%% file: poissondivtype.tex
\documentclass[a4paper,11pt]{amsart}

\usepackage{a4wide}
\usepackage{amsmath, amsfonts, amscd, amssymb, amsthm}
\usepackage{tikz}
\usepackage[UKenglish]{babel}
\usetikzlibrary{matrix,arrows,calc}
\usepackage{xspace}
\usepackage[bottom]{footmisc}
\usepackage{MnSymbol}

\newcommand{\thedocumentname}{Poisson structures of divisor-type}
\newcommand{\theauthor}{Ralph L.\ Klaasse}
\usepackage{aliascnt}
\usepackage[pdftitle={\thedocumentname},pdfauthor={\theauthor},bookmarks,urlcolor=blue,colorlinks=true,linkcolor=black,citecolor=blue]{hyperref}

\theoremstyle{plain}
\newtheorem{thm}{Theorem}[section]

\newaliascnt{lem}{thm}
\newtheorem{lem}[lem]{Lemma}
\aliascntresetthe{lem}

\newaliascnt{prop}{thm}
\newtheorem{prop}[prop]{Proposition}
\aliascntresetthe{prop}               

\newaliascnt{cor}{thm}
\newtheorem{cor}[cor]{Corollary}
\aliascntresetthe{cor}

\newaliascnt{conj}{thm}
\newtheorem{conj}[conj]{Conjecture}
\aliascntresetthe{conj}

\theoremstyle{definition}
\newaliascnt{defn}{thm}
\newtheorem{defn}[defn]{Definition}
\aliascntresetthe{defn}

\newaliascnt{rem}{thm}
\newtheorem{rem}[rem]{Remark}
\aliascntresetthe{rem}

\newaliascnt{exa}{thm}
\newtheorem{exa}[exa]{Example}
\aliascntresetthe{exa}

\theoremstyle{plain}
\newtheorem{thmalpha}{Theorem}

\newcommand{\N}{\mathbb{N}}
\newcommand{\Z}{\mathbb{Z}}
\newcommand{\R}{\mathbb{R}}
\newcommand{\C}{\mathbb{C}}
\newcommand{\bi}{\begin{itemize}}
\newcommand{\ei}{\end{itemize}}
\newcommand{\be}{\begin{equation*}}
\newcommand{\ee}{\end{equation*}}
\newcommand{\bp}[1][]{\begin{proof}[Proof#1.]}
\newcommand{\ep}{\end{proof}}
\newcommand{\wt}{\widetilde}

\newcommand{\mf}{\mathfrak}

\newcommand{\mc}{\mathcal}
\renewcommand{\:}{\text{\rm :}}
\addto\extrasUKenglish{%

}

\numberwithin{equation}{section}
\begin{document}

\author{Ralph L.\ Klaasse}
\address{Department of Mathematics, Utrecht University, 3508 TA Utrecht, The Netherlands}
\address{D\'epartement de Math\'ematique, Universit\'e libre de Bruxelles, Brussels 1050, Belgium}
\subjclass[2010]{53D17}
\email{r.l.klaasse@gmail.com}

\begin{abstract} Poisson structures of divisor-type are those whose degeneracy can be captured by a divisor ideal, which is a locally principal ideal sheaf with nowhere-dense quotient support. This is a large class of Poisson structures which includes all generically-nondegenerate Poisson structures, such as log-, $b^k$-, elliptic, elliptic-log, and scattering Poisson structures.
	
Divisor ideals are used to define almost-injective Lie algebroids of derivations preserving them, to which these Poisson structures can be lifted, often nondegenerately so. The resulting symplectic Lie algebroids can be studied using tools from symplectic geometry. In this paper we develop an effective framework for the study of Poisson structures of divisor-type (also called almost-regular Poisson structures) and their Lie algebroids. We provide lifting criteria for such Poisson structures, develop the language of divisors on smooth manifolds, and discuss residue maps to extract information from Lie algebroid forms along their degeneraci loci.
\end{abstract}

\title{Poisson structures of divisor-type}
\maketitle

\vspace{-2em}
\tableofcontents
\vspace{-3.5em}

\input{introduction}

\input{divisors}

\input{liealgebroidsdivisors}

\input{poissonstructuresdivisors}

\input{residues}

\DeclareRobustCommand{\VAN}[3]{#3}		%for placement of "Van der Leer Duran"
\bibliographystyle{hyperamsplain-nodash}
\bibliography{poissondivtype}
\end{document}

%% file: introduction.tex
\section{Introduction}
\label{sec:introduction}
Poisson geometry can be viewed as the study of singular symplectic foliations. Poisson structures range from the zero Poisson structure, whose symplectic foliation is by points, to nondegenerate Poisson structures having only one symplectic leaf, the entire manifold. Within this wide range of possibilities, one can impose restrictions on the characteristics of the symplectic foliations that may occur, to obtain subclasses of Poisson manifolds.

In this paper we focus on Poisson structures whose symplectic foliation is relatively mild. These are the Poisson structures of divisor-type: those Poisson structures which are generically regular, and whose degeneracy is governed by a divisor ideal. These divisor ideals are locally principal ideals whose local generators have nowhere dense zero set, which we will elaborate on below. Examples of such Poisson structures are those which are generically nondegenerate, i.e.\ those whose union of open symplectic leaves is dense.

The class of Poisson structures of divisor-type coincides with the almost-regular Poisson structures of Androulidakis--Zambon \cite{AndroulidakisZambon17}. They define almost-regular Poisson structures from the point of view of singular foliations, namely as those Poisson structures $\pi$ on a manifold $X$ whose module of Hamiltonian vector fields $\mc{F}_\pi = \pi^\sharp(\Gamma(T^*X)) \subseteq \Gamma(TX)$ is projective. By results of \cite{Debord01}, this class is characterized by having a smooth associated holonomy groupoid.

Our viewpoint -- considering these Poisson structures through their induced divisor -- allows us to  study them using Lie algebroids which are tailor-made to their degeneracy. In this paper we develop a framework for the study of Poisson structures of divisor-type. Below we outline the ideas underlying this framework, and mention our main results.

The techniques and language developed in this paper are used to various establish results for almost-regular Poisson structures, including the computation of Poisson cohomology \cite{KlaasseLanius18}, determining the adjoint symplectic groupoids integrating them \cite{KlaasseLi18}, and developing topological obstructions to their existence \cite{Klaasse18three}. Moreover, this framework was used in \cite{CavalcantiKlaasse18} in the development of Gompf--Thurston methods for symplectic Lie algebroids.
\subsection*{Poisson structures of divisor-type}
In this paper we adapt the language of divisors to the context of smooth manifolds, following their use in \cite{CavalcantiKlaasse18,CavalcantiGualtieri18}. A divisor consist of a line bundle together with a section whose zero set is nowhere dense. Each such pair $(L,s)$ specifies a divisor ideal $I_s := s(\Gamma(L^*)) \subseteq C^\infty(X)$. Abstractly, divisor ideals are locally principal ideals $I$ for which the support $Z(I) \subseteq X$ of the quotient sheaf $C^\infty(X) / \mc{I}$ is nowhere dense, where
\be
	Z(I) = \{x \in X \, | \, f(x) = 0 \text{ for all }f \in \mc{I}_x\}.
\ee
These ideals keep track of degenerative behavior along this support, and are an extension of the concept of a Cartier divisor in algebraic and complex geometry to the smooth setting.

Poisson structures naturally lead to divisors: Given $\pi \in {\rm Poiss}(X)$ such that $\dim(X) = 2n$, consider its Pfaffian $\wedge^n \pi \in \Gamma(\det(TX))$. Note that $\pi$ is nondegenerate exactly at those points $x \in X$ where $\wedge^n \pi_x \neq 0$. Thus $\pi$ is generically nondegenerate if and only if the pair $(\det(TX),\wedge^n \pi)$ is a divisor. In this case we denote this pair by ${\rm div}(\pi)$, with associated divisor ideal $I_\pi \subseteq C^\infty(X)$. This ideal captures the behavior of $\pi$ along $Z(\pi) = Z(I_\pi)$, called the degeneracy locus of $\pi$, consisting of points where the rank of $\pi$ is not maximal. 

In general, given $\pi \in {\rm Poiss}(X)$, the subset $X_{\pi,{\rm reg}} \subseteq X$ of points where $\pi$ is regular is open and dense. In contrast, the subset $X_{\pi, {\rm max}} \subseteq X_{\pi,{\rm reg}}$ of points of maximal rank is not. Given an integer $m \geq 0$, a Poisson structure is of $m$-divisor-type if it is generically of rank $2m$, and its failure to be regular is governed by a divisor. More precisely, for its partial Pfaffians we demand that $\wedge^{m+1} \pi \equiv 0$ and there exists a line subbundle $K \subseteq \wedge^{2m} TX$ such that $(K, \wedge^m \pi)$ is a divisor, which provides a divisor ideal $I_\pi$. A consequence of this definition is that $X_{\pi,{\rm max}} = X_{\pi,{\rm reg}}$ is open and dense, so that $\pi$ is generically regular of rank $2m$. We stress, however, that not every generically regular Poisson structure is of $m$-divisor-type. This definition is closely related to almost-regularity as defined in \cite{AndroulidakisZambon17}, as we discuss in \autoref{sec:almostregularpoisson}.

One class of generically-nondegenerate Poisson structures is that of log-Poisson structures \cite{GuilleminMirandaPires14,MarcutOsornoTorres14,Cavalcanti17,GualtieriLi14}, (also called log- or $b$-symplectic). Another class is given by elliptic Poisson structures \cite{CavalcantiKlaasse18,CavalcantiGualtieri18}, which arise (partially) as the Poisson structures underlying stable generalized complex structures. Further examples include $b^k$-Poisson structures \cite{Scott16,GuilleminMirandaWeitsman17,Lanius16} and scattering Poisson structures \cite{Lanius16}. Moreover, in this paper we introduce a novel class called elliptic-log Poisson structures (a more in-depth study of them is deferred to \cite{KlaasseLanius18}). Each of these Poisson structures of divisor-type has a generically-regular analogue. In particular, there are $m$-log-Poisson structures of generic rank $2m$ (called log-f Poisson structures in \cite{AndroulidakisZambon17}).

By keeping track of the divisor associated to a Poisson structure $\pi$ of divisor-type, we can define almost-injective Lie algebroids $\mc{A} \to X$ whose degeneracy is similar or identical to that of $\pi$. These Lie algebroids are candidates for a lifting procedure, whereby $\pi$ is (partially) desingularized to a Lie algebroid Poisson structure. We first describe such Lie algebroids.
\subsection*{Lie algebroids of divisor-type}
Divisor ideals can also be used to capture the degeneracy of Lie algebroids. This leads to Lie algebroids of divisor-type, which are almost-injective Lie algebroids for which the degeneration of the anchor map is controlled in a precise sense:

Given a Lie algebroid $\mc{A} \to X$, one can view its anchor as a section $\rho_{A} \in \Gamma(\mc{A}^* \otimes TX)$. If ${\rm rank}(\mc{A}) = \dim(X)$, then the zero locus of $\det(\rho_\mc{A})$ consists of those points where $\rho_{A}$ is not an isomorphism. Thus $\mc{A}$ is almost-injective if and only if the pair $(\det(\mc{A}^*) \otimes \det(TX), \det(\rho_{A}))$ is a divisor, which we denote by ${\rm div}(\mc{A})$. In fact, divisors can be associated to any bundle map between bundles of the same rank, so that ${\rm div}(\mc{A}) = {\rm div}(\rho_{A})$. With this in mind, a Lie algebroid $\mc{A} \to X$ is of $m$-divisor-type for some $m \geq 0$ if there exists a Lie subalgebroid $D \subseteq TX$ of rank $m$ so that $\rho_{A}(\mc{A}) \subseteq D$, and for which the induced Lie algebroid morphism $\varphi_{\mc{A}}\colon \mc{A} \to D$ is of divisor-type. This forces that $\mc{A}$ is of rank $m$, and that it is almost-injective.

Above we saw how a Lie algebroid $\mc{A}$ has an associated divisor ${\rm div}(\mc{A})$, and consequently leads to a divisor ideal $I_\mc{A} \subseteq C^\infty(X)$. There are also reverse procedures, where divisor ideals are used to construct Lie algebroids of divisor-type. The first of these proceeds as follows: Given a divisor ideal $I \subseteq C^\infty(X)$, consider the module of derivations preserving the ideal $I$
\be
	\Gamma(TX)_I := \{V \in \Gamma(TX) \, | \, \mc{L}_V I \subseteq I\}.
\ee
If this module is projective, by the Serre--Swan theorem there is a unique almost-injective Lie algebroid $TX_I \to X$ such that $\Gamma(TX_I) \cong \Gamma(TX)_I$. This is called the ideal Lie algebroid associated to $I$, and we call $I$ projective if $\Gamma(TX)_I$ is projective. The Lie algebroid $TX_I$ is of divisor-type, with divisor ideal $I_{TX_I}$. Here $I$ and $I_{TX_I}$ might differ, but are often related. This construction recovers the (complex) log- and elliptic tangent bundles \cite{CavalcantiKlaasse18,CavalcantiGualtieri18,GuilleminMirandaPires14,MarcutOsornoTorres14,Melrose93}.

A second method to construct Lie algebroids  from divisors is the following: Let $\mc{A} \to X$ be a Lie algebroid, and consider a divisor ideal $I \subseteq C^\infty(X)$ whose support $Z(I)$ is a smooth submanifold. Given a Lie subalgebroid $\mc{B}$ of $\mc{A}$ which is supported on $Z(I)$, one can perform elementary modification to construct a new Lie algebroid $[\mc{A}\:(\mc{B},I)]$ which is isomorphic to $\mc{A}$ outside of $Z(I)$. In fact, there is an induced Lie algebroid morphism $\varphi\colon [\mc{A}\:(\mc{B},I)] \to \mc{A}$ which is of divisor-type with divisor ideal $I_\varphi = I^k$, where $k$ is the corank of $\mc{B}$ with respect to $\mc{A}$. This procedure extends the so-called modification procedures (or rescaling) performed in the case where $I = I_Z$, the vanishing ideal of a smooth hypersurface $Z \subseteq X$ (see \cite{Li13,GualtieriLi14,Lanius16,Melrose93}).
\subsection*{Lifting Poisson structures}
Poisson structures and Lie algebroids of divisor-type are connected through the process of lifting. Given a Lie algebroid $\mc{A} \to X$, an $\mc{A}$-Poisson structure $\pi_{\mc{A}} \in {\rm Poiss}(\mc{A})$ is a section $\pi_{A} \in \Gamma(\wedge^2 \mc{A})$ satisfying $[\pi_{\mc{A}},\pi_{\mc{A}}]_{\mc{A}} = 0$, where $[\cdot,\cdot]_\mc{A}$ is the Schouten bracket on $\Gamma(\wedge^\bullet \mc{A})$ extending the Lie bracket on $\Gamma(\mc{A})$. Such a structure turns the tuple $(\mc{A},\mc{A}^*,\pi_{\mc{A}})$ into a triangular Lie bialgebroid (\cite{KosmannSchwarzbach95,KosmannSchwarzbachLaurentGengoux05,LiuWeinsteinXu97,MackenzieXu94}). Given $\pi \in {\rm Poiss}(X)$, we say that $\pi$ is of $\mc{A}$-type if it admits an $\mc{A}$-lift, i.e.\ an $\mc{A}$-Poisson structure $\pi_{\mc{A}}$ such that $\rho_{\mc{A}}(\pi_{\mc{A}}) = \pi$.

If both $\mc{A}$ and $\pi$ are of divisor-type, then so is $\pi_{\mc{A}}$, and their divisor ideals are related:
\be
	I_\pi = I_{\pi_{\mc{A}}} \cdot I_\mc{A}.
\ee
In other words, part of the degeneracy of $\pi$ is absorbed into that of the anchor of $\mc{A}$, making its $\mc{A}$-lift $\pi_{\mc{A}}$ necessarily more regular. Because of the nature of the Lie algebroids considered in this paper, it is possible to use symplectic techniques to study nondegenerate $\mc{A}$-Poisson structures. Indeed, a nondegenerate $\mc{A}$-Poisson structure is dual to an $\mc{A}$-symplectic structure, making $\mc{A}$ into a symplectic Lie algebroid \cite{NestTsygan01}. For example, Moser methods are available for such symplectic Lie algebroids $\mc{A}$ of divisor-type (see \cite[Theorem 5.5]{KlaasseLanius18}), which lead to Darboux-type normal form results for $\pi_{\mc{A}}$, and hence for the underlying $\pi \in {\rm Poiss}(X)$.

Given a Poisson structure $\pi$ of divisor-type, we would like to find a Lie algebroid $\mc{A}$ for which $I_\mc{A} = I_\pi$, so that any $\mc{A}$-lift of $\pi$ will be nondegenerate. Sometimes such Lie algebroids can be constructed from the divisor ideal $I_\pi$. This lifting process necessitates the study of Poisson and symplectic structures in a general Lie algebroid $\mc{A}$, which is the aim of this paper.

While most of this paper is written in terms of $\mc{A}$-Lie algebroids and $\mc{A}$-Poisson structures, our main interest is to study Poisson structures $\pi \in {\rm Poiss}(X)$ in the standard sense. The $\mc{A}$-notions have to be developed for their own sake, but are mainly used as auxiliary tools through the process of lifting. Note that there can be multiple non-isomorphic Lie algebroids to which one can nondegenerately lift a given Poisson structure (see \cite{KlaasseLanius18}).
\subsection*{Results}
The main aim of this paper is to develop an effective framework for the study of Poisson structures of divisor-type. Because of this we introduce and study:
\bi
	\item Divisor ideals keeping track of degeneracies;
	\item Lie algebroids of divisor-type, with emphasis on:
	\bi
		\item Ideal Lie algebroids of derivations preserving a divisor ideal;
		\item Elementary modification of Lie algebroids using divisor ideals;
	\ei
	\item Poisson structures of divisor-type, leading to:
	\bi
		\item The process of lifting Poisson structures to their $\mc{A}$-Poisson counterparts;
	\ei
	\item Residue maps to extract information along the degeneraci loci;
	\item The relation between the lifting Lie algebroid and the modular foliation.
\ei
Here we highlight several aspects of this undertaking that we provide in this paper. 

The first result concerns lifting for Poisson structures of divisor-type to the ideal Lie algebroids they define. This recaptures the lifting results for log- and elliptic Poisson structures of \cite{GuilleminMirandaPires14} and \cite{CavalcantiGualtieri18}. The following is \autoref{thm:ipoissonaitype}, on lifting to the ideal Lie algebroid $TX_I$.
\begin{thmalpha}[Lifting]\label{thm:introlifting} Let $\pi \in {\rm Poiss}(X)$ be of projective $I$-divisor-type with $TX_I^*$ admitting a local basis of closed sections. Then $\pi$ is of $TX_I$-type, and of nondegenerate $TX_I$-type if and only if $I$ is standard, i.e.\ the divisor ideal of $TX_I$ satisfies $I_{TX_I} = I$.
\end{thmalpha}
The next theorem discusses the residue maps of Lie algebroids $\mc{A} \to X$. Given an invariant submanifold $D \subseteq X$ on which the module $(\rho_\mc{A}|_D)(\Gamma(\mc{A}|_D)) \cong \Gamma(\mc{B})$ is projective, there is an induced surjective morphism $\wt{\rho}_\mc{A}|_D\colon \mc{A}|_D \to \mc{B}$ resulting in a Lie algebroid extension sequence
\be
	0 \to E \to \mc{A}|_D \to \mc{B} \to 0.
\ee
The resulting residue maps are used in computing the Lie algebroid cohomology of $\mc{A}$, and to extract information from $\mc{A}$-symplectic structures. In \autoref{sec:residues} we conceptually describe the residue maps found in  \cite{CavalcantiKlaasse18,CavalcantiGualtieri18,GuilleminMirandaPires14,MarcutOsornoTorres14,Scott16}, and are used in \cite{KlaasseLanius18}. The following is \autoref{thm:residues}.
\begin{thmalpha}[Residue maps]\label{thm:introresidues} Let $\mc{A} \to X$ be a Lie algebroid and $D \subseteq X$ be a projective $\mc{A}$-invariant submanifold with induced Lie algebroid $\mc{B} \to D$. Then there exists a residue map
\be
	{\rm Res}_D\colon \Omega^\bullet(X,\mc{A}) \to \Omega^{\bullet-\ell}(D,\mc{B}; \det(E^*)), \qquad \text{where } \ell = {\rm rank}(E).
\ee
This is a cochain morphism for the natural $\mc{B}$-representation on $\det(E^*)$, hence descends to
\be
	[{\rm Res}_D]\colon H^\bullet(X,\mc{A}) \to H^{\bullet-\ell}(D,\mc{B}; \det(E^*)).
\ee
\end{thmalpha}
We further discuss the relation between the lifting process and the modular foliation \cite{GualtieriPym13}, which arises as the sheaf of submodules of $\Gamma(TX)$ spanned by the modular and Hamiltonian vector fields. The symplectic foliation $\mc{F}_\pi$ and the modular foliation $\mc{F}_{\pi,{\rm mod}}$ of a Poisson structure $\pi \in {\rm Poiss}(X)$ interact with any lifting Lie algebroid. The following is \autoref{thm:modfoliation}.
\begin{thmalpha}[Modular foliation]\label{thm:intromodfoliation} Let $\pi \in {\rm Poiss}(X)$ be of $\mc{A}$-type for a Lie algebroid $\mc{A} \to X$. Then $\mc{F}_{\pi} \subseteq \mc{F}_{\pi, {\rm mod}} \subseteq \rho_\mc{A}(\Gamma(\mc{A})) \subseteq \Gamma(TX)$, and the orbits of $\mc{A}$ are $\pi$-Poisson submanifolds.
\end{thmalpha}
This result is most interesting when $\mc{A}$ is almost-injective, i.e.\ when $\rho_\mc{A}$ is injective on sections. In this case \autoref{thm:intromodfoliation} gives an intuitive picture of the lifting procedure. It also helps in determining prospective candidate Lie algebroids to lift Poisson structures to.
\begin{rem} In this paper we do not discuss Dirac geometry nor generalized complex geometry. Their interaction with divisors and the lifting procedure are instead explored in \cite{Klaasse18}, see also \cite{CavalcantiKlaasse18,CavalcantiGualtieri18, Klaasse17}. We do not discuss the Lie groupoids associated to the Lie algebroids we consider, which are all integrable as they are almost-injective \cite{CrainicFernandes03,CrainicFernandes04,Debord01}. See \cite{GualtieriLi14,KlaasseLi18,Li13} for explicit constructions of integrations in the setting of log- and elliptic Poisson manifolds.
\end{rem}
\subsection*{Acknowledgements} This work is partially based on the author's Ph.D.\ thesis \cite{Klaasse17}, and was supported by ERC grant 646649, and by VIDI grant 639.032.221 from NWO, the Netherlands Organisation for Scientific Research. The author would like to thank Melinda Lanius, Rui Loja Fernandes, Ioan M{\u{a}}rcu{\c{t}}, Brent Pym, Marco Zambon, and especially his Ph.D.\ advisor Gil Cavalcanti for useful discussions. Moreover, the author would like to thank the University of Illinois at Urbana--Champaign for its hospitality in early 2017 during a part of this project.
\subsection*{Organization of the paper}
This paper is structured as follows. In \autoref{sec:divonmanifolds} we discuss the theory of divisors on smooth manifolds, emphasizing especially the notion of a divisor ideal. In \autoref{sec:liealgebroidsdivtype} we then turn to Lie algebroids of divisor-type, and how divisor ideals can be used to construct them. Then, in \autoref{sec:poissondivtype} we consider Poisson structures of divisor-type in general Lie algebroids, and discuss the process of lifting Poisson structures along Lie algebroid morphisms. In \autoref{sec:residues} we then discuss residue maps for Lie algebroids and Poisson modules, and how these can be used to study symplectic Lie algebroids. We further study the modular foliation of a Poisson manifold, and how it interacts with the lifting procedure.
%
%END OF TEX FILE

%% file: divisors.tex
\section{Divisors on smooth manifolds}
\label{sec:divonmanifolds}
In this section we introduce real and complex divisors on smooth manifolds (see also \cite{CavalcantiKlaasse18,Klaasse17} and \cite{CavalcantiGualtieri18, VanderLeerDuran16}). These are extensions of the complex geometric notion of divisors to the smooth setting. The theory of divisors permeates much of this paper, as they will be used to keep track of degeneracies of both Lie algebroids and Poisson structures. Let $X$ be a fixed smooth manifold and let $\mc{C}^\infty_X$ be its sheaf of smooth functions, with $C^\infty(X) = \mc{C}^\infty_X(X)$.

In this section we will often not specify whether we are dealing with real or complex objects, and treat them simultaneously and on equal footing. However, in the remainder of the paper we will mostly consider real divisors. Their complex counterparts will play a larger role in the follow-up paper \cite{Klaasse18} (and occur prominently in \cite{CavalcantiKlaasse18,CavalcantiGualtieri18,VanderLeerDuran16}).
\begin{defn} A \emph{divisor} $(U,\sigma)$ on $X$ is a line bundle $U \to X$ together with a section $\sigma \in \Gamma(U)$ with nowhere dense closed zero set $Z_\sigma := \sigma^{-1}(0) \subseteq X$. A \emph{divisor ideal sheaf} is a locally principal ideal sheaf $\mc{I} \subseteq \mc{C}^\infty_X$ generated by nowhere-densely vanishing functions. Its associated \emph{divisor ideal} is the ideal of its global sections, $I := \mc{I}(X) \subseteq C^\infty(X)$.
\end{defn}
\begin{rem} Recall that a nowhere dense subset of a topological space is one whose closure has empty interior. In other words, its intersection with any nonempty open subset is not dense. Often the subset $Z_\sigma$ will be a (union of) submanifold(s) of positive codimension.
\end{rem}
Any divisor specifies a divisor ideal $I_\sigma := \sigma(\Gamma(U^*)) \subseteq C^\infty(X)$ by evaluation, with associated divisor ideal sheaf $\mc{I}_\sigma$. Letting $\alpha$ be a local trivialization of $U^*$, we have $\alpha(\sigma) = g$ for some local function $g$, and hence locally $I_\sigma = \langle \alpha(\sigma) \rangle = \langle g \rangle$. Moreover, note that
\be
	Z(I) := \{x \in X \, | \, f(x) = 0 \text{ for all } f \in \mc{I}_x\}
\ee
is the support of the quotient sheaf $\mc{C}^\infty_X / \mc{I}$, with $Z(I_\sigma) = Z_\sigma$ (here $\mc{I}_x$ denotes the stalk of $\mc{I}$ at $x \in X$). Because of this we will say that $(U,\sigma)$ is \emph{supported} on $Z_\sigma$, and sometimes write $Z_I := Z(I)$. Similarly we say that the divisor ideal $I$ is supported on $Z_I \subseteq X$.

Morphisms of divisors are most conveniently described through their associated divisor ideals. A \emph{morphism of divisors} $(U,\sigma) \to X$ and $(U',\sigma') \to X'$ is a map $f\colon X \to X'$ satisfying $f^* I_{\sigma'} = I_\sigma$, where $f^* I \subseteq C^\infty(X)$ for an ideal $I \subseteq C^\infty(X')$ denotes the ideal generated by the pullback. Divisors are \emph{diffeomorphic} if $f$ is a diffeomorphism, while they are \emph{isomorphic} (written using $=$) if $X = X'$ and the identity map ${\rm id}_X$ is a morphism of divisors.
 
A divisor can be reconstructed up to isomorphism by its divisor ideal. This establishes a bijective correspondence between divisor ideals and isomorphism classes of divisors. This allows us to mostly work with divisor ideals, which carry all essential information. However, often there is a canonical representative of its divisor isomorphism class, for example for Lie algebroid morphisms and Poisson structures.
\begin{prop}[\cite{CavalcantiKlaasse18,VanderLeerDuran16}]\label{prop:locprincideal} Let $I \subseteq C^\infty(X)$ be a divisor ideal. Then there exists a divisor $(U_I,\sigma)$ on $X$ unique up to isomorphism such that $I_\sigma = I$.
\end{prop}
Divisors over a fixed manifold $X$ admit a \emph{product} operation, given by the tensor product $(U,\sigma) \otimes (U',\sigma') := (U \otimes U', \sigma \otimes \sigma')$. The \emph{trivial divisor} $(\underline{\R},\underline{1})$ is the unit up to isomorphism. This operation produces another divisor, because the union of nowhere dense subsets is nowhere dense. The product of divisors corresponds to the product of divisor ideals, i.e.\ $I_{\sigma \otimes \sigma'} = I_\sigma \cdot I_{\sigma'}$, and hence descends to isomorphism classes, with unit $C^\infty(X)$, the trivial divisor ideal. A divisor ideal $I$ is \emph{divisible} by $I'$ if $I = I' \cdot I''$ for some divisor ideal $I''$, and similarly at the level of divisors. We will write $I' < I$ if the divisor ideal $I$ is divisible by $I'$.
\begin{rem} Divisors ideals on $X$ with this notion of morphism and product form an abelian category $\mf{Div}(X)$. Accordingly, we can consider (ir)reducible elements. This leads to asking whether one is able to categorize or classify irreducible divisor ideals. This is akin to the classification of singularities, and related to standardness (\autoref{defn:divprojstandard}).
\end{rem}
\begin{defn}\label{defn:smoothdivisor} A divisor ideal $I \subseteq C^\infty(X)$ is \emph{smooth} if $Z_I \subseteq X$ is a smooth submanifold.
\end{defn}
Here in the notion of submanifold we allow the codimensions of connected components to vary. Any smooth divisor can be uniquely written as the product of divisors on the connected components of its support, thus we can sometimes assume that the support is connected.
\begin{prop}\label{prop:divvanishingideal} Let $(U,\sigma)$ be a divisor. Then $I_\sigma \subseteq I_{Z_\sigma}$, the vanishing ideal of $Z_\sigma$. 
\end{prop}
\bp Locally trivialize $U^*$ by $\alpha^*$ and let $g := \alpha(\sigma) \in \mc{C}^\infty_X$. Then locally $I_\sigma = \sigma(\Gamma(U^*)) = \langle g \rangle$. As $\sigma(Z_\sigma) = 0$, we conclude that $g(Z_\sigma) = \sigma(Z_\sigma) = 0$ so that $g \in I_{Z_\sigma}$, hence $I_\sigma \subseteq I_{Z_\sigma}$. 
\ep
Note that $I_\sigma = I_{Z_\sigma}$ is only possible for $Z_\sigma$ smooth if $Z_\sigma$ is of codimension one, as can be seen by counting local generators.
There is also a relation between morphisms of divisors and their underlying zero sets. We first introduce the following terminology (c.f.\ \cite{CavalcantiKlaasse18}).
\begin{defn}\label{defn:pairsandmaps} A \emph{pair} $(X,Z)$ is a manifold $X$ with a subset $Z \subseteq X$. A \emph{map of pairs} $f\colon (X,Z) \to (X',Z')$ is a smooth map $f\colon X \to X'$ for which $f(Z) \subseteq Z'$. A \emph{strong map of pairs} is a map of pairs $f\colon (X,Z) \to (X,Z')$ for which $f^{-1}(Z') = Z$. A map of pairs $f\colon (X,Z) \to (X',Z')$ where $Z' \subseteq X'$ is a submanifold is \emph{transverse} if $f$ is transverse to $Z'$.
\end{defn}
\begin{prop}\label{prop:divmorphstrong} Let $f\colon (X,U,\sigma) \to (X',U',\sigma')$ be a morphism of divisors. Then $f\colon (X,Z_{\sigma}) \to (X',Z_{\sigma'})$ is a strong map of pairs, equivalently, satisfies the morphism condition $f^* I_{Z_{\sigma'}} = I_{Z_{\sigma}}$.
\end{prop}
In general, for a map $f\colon X \to Y$ and $Z \subseteq Y$ closed, one can only conclude that $f^* I_Z \subseteq I_{f^{-1}(Z)}$.
\bp As $f$ is a morphism of divisors, we have that $U \cong f^* U'$ and $\sigma = g f^* \sigma'$ for a nonvanishing function $g$. As $Z_{\sigma}$ is by definition the zero set of $\sigma$, it is immediate that it equals $f^{-1}(Z_{\sigma'})$.
\ep

Given a smooth map $f\colon X \to X'$ and a divisor $(U',\sigma')$ on $X'$, one can equip $X$ with a divisor by setting $(U,\sigma) := (f^* U', f^* \sigma')$, as long as $f^{-1}(Z_{\sigma'})$ is nowhere dense (e.g.\ by assuming $f$ is a transverse map of pairs onto $(X,Z_{\sigma'})$ if $Z_{\sigma'}$ is smooth). This automatically makes $f$ a morphism of divisors. If $f$ is transverse to $Z_{\sigma'}$, the divisor-type is preserved under pullback: by this we mean that the type of degeneracy is similar, so that the class of divisor is preserved.

Divisors admit another type of product. Let $(X,U,\sigma)$ and $(X',U',\sigma')$ be two divisors and consider the product $X \times X'$ with projections $p_{X}\colon X \times X' \to X$ and $p_{X'}\colon X \times X' \to X'$.
\begin{defn}\label{defn:divdirectstum} The \emph{external tensor product} of $(U,\sigma)$ and $(U',\sigma')$ is the divisor
\be
	(p_X^* U \otimes p_{X'}^* U', p_X^* \sigma \otimes p_{X'}^* \sigma') \to X \times X',
\ee
with zero set $Z_\sigma \times X' \cup X \times Z_{\sigma'}$ and divisor ideal given by the product $p_X^* I_\sigma \cdot p_{X'}^* I_{\sigma'} \subseteq C^\infty(X \times X')$.
\end{defn}
\subsection{Examples}
\label{sec:divexamples}
Divisors of various flavors exist, and can naturally be combined using the product operation and external tensor product to form more involved examples. Here we give several examples.
\begin{exa}[Trivial divisor] Let $U$ be the trivial line bundle with $\sigma \in \Gamma(U)$ nonvanishing. Then $Z_\sigma$ is empty, and $(U,\sigma)$ is called the \emph{trivial divisor} on $X$. In this case, $I_\sigma = C^\infty(X)$.
\end{exa}
\begin{exa}[Log divisors \cite{CavalcantiKlaasse18}]\label{exa:log} A \emph{log divisor} is a divisor for which $\sigma$ is transverse to the zero section in $\Gamma(U)$. Its zero set $Z = Z_\sigma$ is a smooth hypersurface, and the corresponding divisor ideal is the vanishing ideal $I_Z$, locally generated around $Z$ by $\langle z \rangle$, with $z$ a local defining function for $Z$. As shown in \cite{Klaasse17}, this means that $Z$ specifies a unique isomorphism class of log divisors with zero set $Z$. We will henceforth denote log divisors by $Z = (L_Z,s)$.
\end{exa}
\begin{exa}[Normal-crossing log divisors \cite{GualtieriLiPelayoRatiu17}]\label{exa:normalcrossinglog} The above example can be extended to normal-crossing log divisors $\underline{Z} := \bigcup_j Z_j$ associated to a collection $\{Z_j\}$ of hypersurfaces which intersect transversely. Such a divisor is the product of the individual log divisors for each $Z_j$, hence locally $I_{\underline{Z}} = \langle z_1 \cdot \ldots \cdot z_n \rangle$. These specify hyperbolic singularities on two-fold intersections: there we can readily change coordinates so that $I_{\underline{Z}} = \langle z_1 z_2 \rangle = \langle z'^2_1 - z'^2_2 \rangle$. There is also an extension of these ideas to immersed hypersurfaces with transverse self-crossings \cite{MirandaScott18}.
\end{exa}
There is another example resembling \autoref{exa:normalcrossinglog}, namely that of a \emph{star log divisor} \cite{Lanius17}, consisting of the product of log divisors associated to a collection $\{Z_j\}$ of pairwise transversely intersecting hypersurfaces $Z_j \subseteq X$.
\begin{exa}[Elliptic divisors \cite{CavalcantiKlaasse18,CavalcantiGualtieri18}]\label{exa:elliptic} An \emph{elliptic divisor} is one for which $\sigma$ vanishes along a smooth codimension-two critical submanifold along which its normal Hessian
\be
	{\rm Hess}(\sigma) \in \Gamma(Z_\sigma; {\rm Sym}^2 NZ_\sigma \otimes U)
\ee
is definite. We will denote $Z_\sigma$ by $D$ and the elliptic divisor by $|D| = (R,q)$. As shown in \cite{CavalcantiKlaasse18,CavalcantiGualtieri18}, the associated divisor ideal $I_{|D|} := I_q$ is locally generated by $\langle r^2 \rangle$ in polar coordinates $(r,\theta)$ normal to $D$. Note that $I_{|D|}$ is not the vanishing ideal $I_D$ of $D$.
\end{exa}
\begin{rem} Hypersurfaces carry a unique log divisor ideal structure. This does not hold for codimension-two submanifolds and elliptic divisors. A simple example is provided by $X = \R^2$ with $D = \{(0,0)\}$ and coordinates $(x,y)$. Equip $D$ with the elliptic ideals $I = \langle x^2 + y^2 \rangle$ and $I' = \langle x^2 + 2y^2 \rangle$. As these ideals are distinct, they supply $D$ with two non-isomorphic elliptic divisor structures. However, it is easy to see that these elliptic divisors are diffeomorphic.
\end{rem}
We construct a more elaborate class of divisors out of log and elliptic divisors as follows.
\begin{exa}[Elliptic-log divisors]\label{exa:ellipticlog} An \emph{elliptic-log divisor} is a divisor obtained as the product of a log divisor $Z = (L_Z,s)$ and an elliptic divisor $|D| = (R,q)$ such that $D \subseteq Z$ (we will write $|D| \subseteq Z$ as a shorthand for this). We denote this product by $W := Z \otimes |D|$, with divisor ideal $I_W = I_Z \cdot I_{|D|}$. As such, pointwise around $D$ there exist normal coordinates $(x,y)$ to $D$ with $\{x = 0\}$ cutting out $Z$, for which $I_W = \langle x (x^2 + y^2) \rangle$, with $r^2 = x^2 + y^2$.
\end{exa}
It is also possible for a log divisor $Z$ and an elliptic divisor $|D|$ to either not or only partially intersect. In general, one can readily combine divisors with nowhere intersecting supports using the product operation. However, most interesting are cases where the divisors do intersect (yet retain interesting properties), of which \autoref{exa:ellipticlog} is a specific instance.
\begin{exa}[Complex log divisors \cite{CavalcantiGualtieri18}]\label{exa:complexlog} A \emph{complex log divisor} is a complex divisor $(U,\sigma)$ for which $\sigma$ is transverse to zero. This is the complex analogue of \autoref{exa:log}, and locally we have $I_\sigma = \langle w \rangle$ for $w$ a local complex function cutting out $D$. Its zero set $D$ is smooth of codimension two, and by complex conjugation we obtain another complex log divisor $(\overline{U},\overline{\sigma})$ with $I_{\overline{\sigma}} = \langle \overline{w} \rangle$. After taking their product we can extract an elliptic divisor $((U \otimes \overline{U})_\R,\sigma \otimes \overline{\sigma})$ with $I_{\sigma \otimes \overline{\sigma}} = \langle w \overline{w} \rangle = \langle r^2 \rangle$. Given an elliptic divisor and a choice of coorientation for $D$, there is a complex log divisor inducing it which is unique up to diffeomorphism.
\end{exa}
There is a general class of divisors of \emph{Morse--Bott type}, i.e.\ those pairs $(U,\sigma)$ for which $\sigma$ is a \emph{Morse--Bott section}. This means that $Z_\sigma$ consists of nondegenerate critical submanifolds, i.e.\ $Z_\sigma$ is a union of submanifolds $S$ of $\sigma$, where $\ker {\rm Hess}(\sigma)_p = T_p S$ for all $p \in S$. If $U$ is trivial, then a Morse--Bott section for $U$ is not the same thing as a Morse--Bott function on $X$, as we only demand a condition on its zero set. In this sense, an alternate way of phrasing the definition is that $\sigma$ must specify the germ of a Morse--Bott function around $Z_\sigma$.

The elliptic divisors of \autoref{exa:elliptic} are of Morse--Bott type, but normal-crossing log divisors (\autoref{exa:normalcrossinglog}) are not. The analogue of the linearization result for elliptic divisors (see \cite[Proposition 2.20]{CavalcantiKlaasse18}) is readily seen to hold for all divisors of Morse--Bott type (see \cite{Klaasse17}), so that locally the associated divisor ideal $I_\sigma$ will be given by $\langle Q_g \rangle$ with $Q_g$ a homogeneous function in normal bundle coordinates. Note that all Morse--Bott type divisors are smooth.
\begin{rem} Elliptic divisors arise as the real divisors which underlie complex log divisors (c.f.\ \cite{CavalcantiGualtieri18}). There is another as of yet unexplored class of real divisors, namely those underlying \emph{quaternionic log divisors}, i.e.\ coming from a quaternionic line bundle with a transversally vanishing section. These are supported on smooth submanifolds of codimension four, and will be explored further in future work.
\end{rem}
%
%END OF TEX FILE

%% file: liealgebroidsdivisors.tex
\section{Lie algebroids of divisor-type}
\label{sec:liealgebroidsdivtype}
In this section we discuss Lie algebroids and their interaction with divisors, starting in \autoref{sec:liealgebroids} with how Lie algebroids can specify divisors through the degeneracy of their anchor. \autoref{sec:aliealgebroids} through \autoref{sec:singularliealgebroids} discuss several further concepts that we will use.

We describe two processes by which Lie algebroids can be constructed using divisors. The first of these gives rise to ideal Lie algebroids (in \autoref{sec:idealliealgebroids}), which are those Lie algebroids whose sections are those preserving a given ideal of functions. The second is that of elementary modification (in \autoref{sec:elementarymodifications}). Here one uses a injective morphism from, or a surjective comorphism to, a subbundle supported on the zero locus of a smooth divisor ideal. 
\subsection{Lie algebroids}
\label{sec:liealgebroids}
For general information regarding Lie algebroids, see \cite{Mackenzie05}. A vector bundle $\mc{A} \to X$ is \emph{anchored} if it has a bundle map $\rho_{A}\colon \mc{A} \to TX$ called its \emph{anchor}. We refer to morphisms of anchored bundles (those intertwining the anchor maps) by \emph{anchored morphisms}.
\begin{defn} \label{defn:liealgebroid} A \emph{Lie algebroid} $(\mathcal{A},[\cdot,\cdot]_{\mc{A}},\rho_\mc{A})$ is an anchored bundle with a Lie bracket $[\cdot,\cdot]_{\mc{A}}$ on $\Gamma(\mc{A})$ such that $[v,f w]_{\mc{A}} = f [v,w]_{\mc{A}} + \mc{L}_{\rho_{\mc{A}}(v)} f \cdot w$ for all $v, w \in \Gamma(\mc{A})$, $f \in C^\infty(X)$.
\end{defn}
The bracket $[\cdot,\cdot]_\mc{A}$ of a Lie algebroid $\mc{A}$ provides the algebra $\Omega^\bullet(\mc{A}) = \Gamma(\wedge^\bullet \mc{A}^*)$ with a differential $d_\mc{A}$ using the Koszul formula. This allows for a simple definition of morphisms.
\begin{defn} A \emph{Lie algebroid morphism} is a bundle morphism $(\varphi,f)\colon \mc{A} \to \mc{A}'$ for which the induced map $\varphi^*\colon \Omega^*(\mc{A}') \to \Omega^*(\mc{A})$ is a chain map, i.e.\ satisfies $\varphi^* \circ d_{\mc{A}'} = d_\mc{A} \circ \varphi^*$.
\end{defn}
\begin{rem} Given $v \in \Gamma(\mc{A})$, its \emph{$\mc{A}$-Lie derivative} is the operator $\mc{L}_v = \{\iota_v, d_\mc{A}\}$ on $\Omega^\bullet(\mc{A})$, where $\{\cdot,\cdot\}$ is the graded commutator. It satisfies the usual Cartan relations (see e.g.\ \cite{Marle08}).
\end{rem}
Lie algebroid morphisms are in particular anchored morphisms, i.e.\ they satisfy the relation $\rho_{\mc{A}'} \circ \varphi = Tf \circ \rho_\mc{A}$. If $(\varphi,f)\colon \mc{A} \to \mc{A}'$ is a base-preserving anchored morphism, it is a Lie algebroid morphism if and only if it moreover induces a Lie algebra homomorphism $\varphi\colon \Gamma(\mc{A}) \to \Gamma(\mc{A}')$. Note that the anchor map $\rho_\mc{A}$ is a Lie algebroid morphism from $\mc{A}$ to $TX$. For base-preserving Lie algebroid morphisms, we keep track of where it is an isomorphism.
\begin{defn} Let $(\varphi,{\rm id}_X)\colon \mc{A}' \to \mc{A}$ be a Lie algebroid morphism. Its \emph{isomorphism locus} is the open subset $X_\varphi \subseteq X$ where $\varphi$ is an isomorphism.
\end{defn}
The complement of the isomorphism locus is denoted by $Z_\varphi := X \backslash X_\varphi$, and is called the \emph{degeneracy locus} of $\varphi$. For a Lie algebroid $\mc{A}$ we define $X_{\mc{A}} := X_{\rho_\mc{A}}$ and refer to this as the isomorphism locus of $\mc{A}$, and similarly for $Z_\mc{A}$. Note that whenever $X_\varphi$ is nonempty, we have ${\rm rank}(\mc{A}) = {\rm rank}(\mc{A}')$ because $X_\varphi$ is open. Between bundles of the same rank, the condition that $X_\varphi$ is dense is equivalent to injectivity of the induced map on sections $\varphi\colon \Gamma(\mc{A}') \to \Gamma(\mc{A})$.
\begin{rem} Given two composable base-preserving Lie algebroid morphisms $\varphi$ and $\varphi'$, we have $X_{\varphi \circ \varphi'} \supseteq X_\varphi \cap X_{\varphi'}$. If we further assume that all Lie algebroids involved are of the same rank (e.g.\ by each isomorphism locus being nonempty), this becomes an equality. In particular, in this case, if the isomorphism locus $X_{\varphi \circ \varphi'}$ is dense, then so are $X_\varphi$ and $X_{\varphi'}$.
\end{rem}
Given a Lie algebroid morphism $(\varphi,{\rm id}_X)\colon \mc{A}' \to \mc{A}$, we can take determinants to obtain $\det(\varphi)\colon \det(\mc{A}') \to \det(\mc{A})$. This specifies a section $\det(\varphi) \in \Gamma(\det(\mc{A}'^*) \otimes \det(\mc{A}))$. It is now readily verified that $\varphi$ has dense isomorphism locus if and only if this defines a divisor.
\begin{defn}\label{defn:morphismdivisortype} Let $(\varphi,{\rm id}_X)\colon \mc{A}' \to \mc{A}$ be a Lie algebroid morphism with dense isomorphism locus. The \emph{divisor associated to $\varphi$} is ${\rm div}(\varphi) = (\det(\mc{A}'^*) \otimes \det(\mc{A}), \det(\varphi))$, with the associated divisor ideal being denoted by $I_\varphi$. Given a Lie algebroid $\rho_\mc{A}\colon \mc{A} \to TX$ with dense isomorphism locus, we define ${\rm div}(\mc{A}) := {\rm div}(\rho_\mc{A})$ and $I_\mc{A} := I_{\rho_\mc{A}}$, and call $\mc{A}$ of \emph{divisor-type}.
\end{defn}
Given a divisor ideal $I$ we call a base-preserving Lie algebroid morphism $\varphi$ \emph{of $I$-divisor-type} if $I_\varphi = I$, and similarly for Lie algebroids $\mc{A}$. It is immediate that composition of morphisms of divisor-type corresponds to taking products of their divisors (hence their divisors ideals): after noting that $\det(\varphi' \circ \varphi) = \det(\varphi') \circ \det(\varphi)$ we obtain the equality $I_{\varphi' \circ \varphi} = I_{\varphi'} \cdot I_\varphi$.
\begin{prop}\label{prop:divisorcomposition} Let $(\varphi,{\rm id}_X)\colon \mc{A} \to \mc{A}'$ and $(\varphi',{\rm id}_X)\colon \mc{A}' \to \mc{A}''$ be Lie algebroid morphisms of divisor-type. Then $\varphi' \circ \varphi$ is of divisor-type, and ${\rm div}(\varphi' \circ \varphi) \cong {\rm div}(\varphi') \otimes {\rm div}(\varphi)$.
\end{prop}
In the generically-regular setting, we define a Lie algebroid morphism $(\varphi,{\rm id}_X)\colon \mc{A}' \to \mc{A}$ to be of \emph{$m$-divisor-type} for some $m \geq 0$ if there exists a regular Lie subalgebroid $D_\mc{A} \subseteq \mc{A}$ of rank $m$ (see \autoref{sec:liesubalgds}, in particular \autoref{exa:adistribution}) such that $\varphi$ factors through the inclusion $\varphi_{D_\mc{A}}\colon D_{\mc{A}} \to \mc{A}$, i.e.\ $\varphi = \varphi_{D_\mc{A}} \circ \varphi_{\mc{A}'}$ for some Lie algebroid morphism $\varphi_{\mc{A}'}\colon \mc{A}' \to D_{\mc{A}}$, and $\varphi_{\mc{A}'}$ is of divisor-type, with divisor given by ${\rm div}(\varphi_{\mc{A}'}) = (\det(\mc{A}'^*) \otimes \det(D_{\mc{A}}), \det(\varphi_{\mc{A}'}))$. Alternatively, we can demand the existence of a line subbundle $K_\mc{A} \subseteq \wedge^m \mc{A}$ such that $\wedge^m \varphi\colon \det(\mc{A}') \to K_\mc{A}$ and $(\det(\mc{A}'^*)\otimes K_\mc{A}, \wedge^m \varphi)$ is a divisor. These formulations are related via $\det(D_\mc{A}) = K_\mc{A}$. We then define $\mc{A}$ to be of $m$-divisor-type if its anchor $\rho_{A}\colon \mc{A} \to TX$ is.

We introduce several further adjectives to describe the behavior of the anchor.
\begin{defn}\label{defn:lalgebroidadjective} Let $\mc{A} \to X$ be a Lie algebroid. Then $\mc{A}$ is \emph{regular} if $\rho_\mc{A}$ has constant rank, \emph{transitive} if $\rho_\mc{A}$ is surjective, and \emph{totally intransitive} if $\rho_\mc{A} \equiv 0$. Further, $\mc{A}$ is \emph{almost-injective} \cite{Debord01} if $\rho_\mc{A}\colon \Gamma(\mc{A}) \to \Gamma(TX)$ is injective, and \emph{projective} if $\rho_{A}(\Gamma(\mc{A})) \subseteq \Gamma(TX)$ is projective.
\end{defn}
These properties are related: transitivity implies regularity, which in turn implies projectivity. If ${\rm rank}(\mc{A}) = \dim X$, then $\mc{A}$ is almost-injective if and only if it is of divisor-type.
\begin{defn}\label{defn:laisotropy} The \emph{isotropy} of $\mc{A} \to X$ at $x \in X$ is given by the subspace $\ker \rho_{\mc{A},x} \subseteq \mc{A}_x$.
\end{defn}
This subspace makes sense for any anchored vector bundle, but for Lie algebroids the isotropy is a Lie algebra, inheriting a bracket from $\mc{A}$. The isotropy Lie algebras $\ker \rho_{\mc{A},x}$ may vary in dimension as $x \in X$ varies, so that in general these may not form a subbundle of $\mc{A}$. A Lie algebroid (or anchored) morphism $(\varphi,f)\colon (\mc{A},X) \to (\mc{A}',X')$ restricts to a map $\varphi\colon \ker \rho_{\mc{A},x} \to \ker \rho_{\mc{A}',f(x)}$ between pointwise isotropies for all $x \in X$.

The image of any Lie algebroid under the anchor map determines a singular foliation
\be
	\mc{F}_\mc{A} := \rho_{\mc{A}}(\Gamma(\mc{A})) \subseteq \Gamma(TX)
\ee
in the sense of Androulidakis--Skandalis \cite{AndroulidakisSkandalis09}. Its leaves, i.e.\ the associated maximal immersed submanifolds $\mc{O}$ satisfying $T_x \mc{O} = {\rm im}\, \rho_{\mc{A},x}$ for all $x \in \mc{O}$, are called the \emph{orbits} of $\mc{A}$.

For regular Lie algebroids (hence for transitive ones), the kernel of $\rho_\mc{A}$ is of constant dimension, so that it is a subbundle of $\mc{A}$. If $\mc{A}$ is transitive, the isotropy Lie algebras together exhibit $\mc{A}$ as an abelian extension of $TX$. In other words, there is a short exact sequence:
\be
	0 \to \ker \rho_\mc{A} \to \mc{A} \to TX \to 0.
\ee
A similar sequence exists for projective Lie algebroids. In this case (as in \autoref{sec:singularliealgebroids}), there exists an almost-injective Lie algebroid $\mc{B} \to X$ with $\Gamma(\mc{B}) \cong \rho_{A}(\Gamma(\mc{A}))$ which is unique up to isomorphism. This comes with an induced Lie algebroid morphism $\wt{\rho}_\mc{A}\colon \mc{A} \to \mc{B}$ induced by the map on sections $\rho_{A}\colon \Gamma(\mc{A}) \to \Gamma(\mc{B})$. This results in a short exact sequence of Lie algebroids:
\be
	0 \to \ker \wt{\rho}_\mc{A} \to \mc{A} \to \mc{B} \to 0.
\ee
Unlike the transitive case, the Lie algebroid $\ker \wt{\rho}_\mc{A}$ need not be abelian. Following \cite[Definition 1.10]{AndroulidakisZambon17}, we make the following definition in this case.
\begin{defn}\label{defn:germinalisotropy} The \emph{germinal isotropy} of $\mc{A} \to X$ at $x \in X$ is the subspace $\ker \wt{\rho}_{\mc{A},x} \subseteq \mc{A}_x$.
\end{defn}
In general one has $\ker \wt{\rho}_{\mc{A},x} \subseteq \ker \rho_{\mc{A},x}$. This definition still makes sense when the module $\mc{F}_\mc{A} \subseteq \Gamma(TX)$ is projective: if it is not, we can still consider the vector spaces $\mc{F}_\mc{A} / I_x \cdot \mc{F}_\mc{A}$. We will return to the germinal isotropy in \autoref{sec:singularliealgebroids}.\\

Let $\mc{A} \to X$ be a Lie algebroid with dense isomorphism locus $X_\mc{A}$. Then the inclusion $i\colon X_{\mc{A}} \hookrightarrow X$ gives a bijection $\rho_{\mc{A}}^*\colon \Omega^k(X_{\mc{A}}) \to \Omega^k_{\mc{A}}(X_{\mc{A}})$ for all $k$. Consequently, we can view $\mc{A}$-forms as smooth forms with certain ``singular'' behavior at its degeneracy locus $Z_\mc{A}$. The isomorphism given by the restricted anchor map $\rho_\mc{A}\colon \mc{A}|_{X_\mc{A}} \to TX|_{X_\mc{A}}$ also implies the following.
\begin{prop}\label{prop:denseisolocusstrongmap} Let $(\varphi,f)\colon (\mc{A},X) \to (\mc{A}',X')$ be an anchored morphism between Lie algebroids of divisor-type such that $f^{-1}(X_{\mc{A}'}) = X_\mc{A}$. Then $\varphi \equiv Tf$ on sections, and $\varphi$ is a Lie algebroid morphism.
\end{prop}
\bp When restricted to the isomorphism loci, we have that $\varphi \equiv Tf$ as bundle maps after intertwining with the anchors. Because the isomorphism loci are dense, this means that $\varphi \equiv Tf$ on sections, which holds globally on $X$. Because $\varphi$ is a bundle morphism, the map $\varphi^*$ is an algebra morphism. In the isomorphism loci, $\varphi$ must equal $Tf$, and $(Tf,f)$ is a Lie algebroid morphism between $TX$ and $TX'$, which means that $f^*$ is a chain map. By the density of $X_{\mc{A}}$ in $X$, we conclude that the map $\varphi^*$ is a chain map everywhere.
\ep
Thus, for Lie algebroids as \autoref{prop:denseisolocusstrongmap}, to determine whether there is a Lie algebroid morphism $(\varphi, f)$, one can check whether $f^*$ extends to a map $\varphi^*$ on forms. This in turn follows if it holds on generators, so it suffices to verify that $f^*$ extends to a map $\varphi^*\colon \Omega^1(\mc{A}') \to \Omega^1(\mc{A})$.
\subsection{\texorpdfstring{$\mc{A}$}{A}-Lie algebroids}
\label{sec:aliealgebroids}
We now introduce what we call $\mc{A}$-versions of Lie algebroids, where $\mc{A} \to X$ is a fixed Lie algebroid. The main idea is that these objects come equipped with a preferred Lie algebroid morphism onto $\mc{A}$, called its $\mc{A}$-anchor, through which its own anchor factors. Our interest later will be in those examples for which this $\mc{A}$-anchor is of divisor-type.
\begin{defn}\label{defn:aavala} An \emph{$\mc{A}$-anchored vector bundle} is an anchored vector bundle $E_\mc{A} \to X$ with an anchored morphism $\varphi_{E_\mc{A}} \colon E_\mc{A} \to \mc{A}$ called the \emph{$\mc{A}$-anchor} satisfying $\rho_{E_\mc{A}} = \rho_\mc{A} \circ \varphi_{E_\mc{A}}$. In turn, an \emph{$\mc{A}$-Lie algebroid} is an $\mc{A}$-anchored vector bundle $E_\mc{A}$ equipped with a compatible Lie algebroid structure for which $\varphi_{E_\mc{A}}$ is a Lie algebroid morphism.
\end{defn}
This concept is useful especially when dealing with Lie algebroids whose sections form a submodule of $\Gamma(\mc{A})$, i.e.\ those whose $\mc{A}$-anchor $\varphi_{E_\mc{A}}\colon \Gamma(E_{\mc{A}}) \to \Gamma(\mc{A})$ is injective on sections. \autoref{defn:morphismdivisortype} provides any $\mc{A}$-Lie algebroid $\mc{A}'$ with the notion of being of \emph{$\mc{A}$-divisor-type} by looking at the pair $(\det(\mc{A}'^*)\otimes \det(\mc{A}), \det(\varphi_{\mc{A}}))$, with its \emph{$\mc{A}$-divisor ideal} being $I_{\varphi_{\mc{A}}}$. Note also the relevance of \autoref{prop:divisorcomposition} which shows that $I_{\mc{A}'} = I_{\varphi_{\mc{A}}} \cdot I_\mc{A}$ when this makes sense.
\begin{rem} Any Lie algebroid $\mc{A} \to X$ is naturally both an $\mc{A}$- and a $TX$-Lie algebroid.
\end{rem}
Recall that a \emph{comorphism of anchored vector bundles} $(\varphi;f)\colon (E,X) \dashedrightarrow (E',X')$ consists of a base map $f\colon X \to X'$ and a vector bundle morphism $\varphi\colon f^* E' \to E$ such that
\be
	\rho_E (\varphi^*(\sigma)) \sim_f \rho_{E'}(\sigma) \qquad \text{for all } \sigma \in \Gamma(E').
\ee
 We will call these \emph{anchored comorphisms} for short. An anchored comorphism gives a pullback map $\varphi^*\colon \Gamma(E') \to \Gamma(E)$ on sections. A \emph{Lie algebroid comorphism} $(\varphi;f)\colon (E,X) \dashedrightarrow (E',X')$ is an anchored comorphism for which the induced map $\varphi^*\colon \Gamma(E') \to \Gamma(E)$ preserves brackets, i.e.\ satisfies $\varphi^*[\sigma,\tau]_{E'} = [\varphi^*(\sigma),\varphi^*(\tau)]_{E}$ for all $\sigma,\tau \in \Gamma(E')$.
\begin{rem} Using the fact that Lie algebroid structures on $\mc{A}$ are in one-to-one correspondence with linear Poisson structures on $\mc{A}^*$, a Lie algebroid comorphism is alternatively defined as a bundle morphism between Lie algebroids whose dual morphism is a Poisson map.
\end{rem}
The $\mc{A}$-version of the notions of morphism and comorphism are rather straightforward.
\begin{defn} An \emph{$\mc{A}$-anchored morphism} $(\varphi',\varphi,f)\colon (E_\mc{A},\mc{A},X) \to (E_{\mc{A}'}, \mc{A}',X')$ consists of an anchored morphism $(\varphi',f)\colon (E_\mc{A},X) \to (E_{\mc{A}'},X')$ and a compatible Lie algebroid morphism $(\varphi,f)\colon (\mc{A},X) \to (\mc{A}',X')$, i.e.\ such that we have the intertwining relation $\varphi \circ \varphi_{E_\mc{A}} = \varphi_{E_\mc{A}'} \circ \varphi'$.
	
An \emph{$\mc{A}$-Lie algebroid morphism} $(\varphi',\varphi,f)\colon (E_\mc{A},\mc{A},X) \to (E_{\mc{A}'}, \mc{A}',X')$ is an $\mc{A}$-anchored morphism for which $(\varphi,f)\colon (E_\mc{A},X) \to (E_{\mc{A}'}, X')$ is a Lie algebroid morphism.
\end{defn}
\begin{defn} An \emph{$\mc{A}$-anchored comorphism} $(\varphi';\varphi,f)\colon (E_\mc{A},\mc{A},X) \dashedrightarrow (E_{\mc{A}'}, \mc{A}',X')$ consists of an anchored comorphism $(\varphi';f)\colon (E_\mc{A},X) \dashedrightarrow (E_{\mc{A}'},X')$ and a compatible Lie algebroid morphism $(\varphi,f)\colon (\mc{A},X) \to (\mc{A}',X')$, i.e.\ such that $\varphi_{E_\mc{A}}(\varphi^*(\sigma)) \sim_{\varphi'} \varphi_{E_\mc{A}'}(\sigma)$ for all $\sigma \in \Gamma(E_{\mc{A}'})$.  An \emph{$\mc{A}$-Lie algebroid comorphism} $(\varphi';\varphi,f)\colon (E_\mc{A},\mc{A},X) \dashedrightarrow (E_{\mc{A}'}, \mc{A}',X')$ is an $\mc{A}$-anchored comorphism for which $(\varphi;f)\colon (E_\mc{A},X) \to (E_{\mc{A}'}, X')$ is a Lie algebroid comorphism.
\end{defn}
To summarize, $\mc{A}$-Lie algebroid (co)morphisms give rise to diagrams as follows:
\begin{center}
	\begin{tikzpicture}
	\begin{scope}[xshift=-200]
	\matrix (m) [matrix of math nodes, row sep=2.5em, column sep=2.5em,text height=1.5ex, text depth=0.25ex]
	{	\rho_{E_\mc{A}}\colon &[-3.1em] E_\mc{A} & \mc{A} & TX \\ \rho_{E_{\mc{A}'}}\colon &[-3.2em] E_{\mc{A}'} & \mc{A}' & TX'\\};
	\path[-stealth]
	(m-1-2) edge node [above] {$\varphi_{E_\mc{A}}$} (m-1-3)
	(m-1-2) edge node [left] {$\varphi$} (m-2-2)
	(m-2-2) edge node [above] {$\varphi_{E_{\mc{A}'}}$} (m-2-3)
	(m-1-3) edge node [left] {$\varphi'$} (m-2-3)
	(m-1-3) edge node [above] {$\rho_\mc{A}$} (m-1-4)
	(m-2-3) edge node [above] {$\rho_{\mc{A}'}$} (m-2-4)
	(m-1-4) edge node [left] {$T f$} (m-2-4);
	\end{scope}
	\begin{scope}
	\matrix (m) [matrix of math nodes, row sep=2.5em, column sep=2.5em,text height=1.5ex, text depth=0.25ex]
	{	\rho_{E_\mc{A}}\colon &[-3.1em] E_\mc{A} & \mc{A} & TX\\ \rho_{E_{\mc{A}'}}\colon &[-3.2em] E_{\mc{A}'} & \mc{A}' & TX'\\};
	\path[-stealth]
	(m-1-2) edge node [above] {$\varphi_{E_\mc{A}}$} (m-1-3)
	(m-2-2) edge node [above] {$\varphi_{E_{\mc{A}'}}$} (m-2-3)
	(m-1-3) edge node [left] {$\varphi'$} (m-2-3)
	(m-1-3) edge node [above] {$\rho_\mc{A}$} (m-1-4)
	(m-2-3) edge node [above] {$\rho_{\mc{A}'}$} (m-2-4)
	(m-1-4) edge node [left] {$T f$} (m-2-4);
	\path[-stealth, dashed]
	(m-1-2) edge node [left] {$\varphi$} (m-2-2);
	\end{scope}
	\end{tikzpicture}
\end{center}
\subsection{Fiber products}
\label{sec:algoperations}
We introduce the fiber product of Lie algebroids over the same base. This is a special case of the fiber product construction for Lie algebroids as found in \cite{HigginsMackenzie90}.
\begin{defn}\label{def:fiberproduct} Consider a Lie algebroid $\mc{A} \to X$ and $\mc{A}$-Lie algebroids $(\varphi_{\mc{A}'}, {\rm id}_X)\colon \mc{A}' \to \mc{A}$ and $(\varphi_{\mc{A}''}, {\rm id}_X)\colon \mc{A}'' \to \mc{A}$. Assume that $\varphi_{\mc{A}'}$ and $\varphi_{\mc{A}''}$ are transverse in $\mc{A}$, i.e.\
\be
	\varphi_{\mc{A}'}(\mc{A}') + \varphi_{\mc{A}''}(\mc{A}'') = \mc{A}
\ee
at each point in $X$. The \emph{fiber product} of $\mc{A}'$ and $\mc{A}''$ over $\mc{A}$ is the Lie algebroid defined as
	\be
	\mc{A}' \times_{\mc{A}} \mc{A}'' = \{(v,w) \in \mc{A}' \times \mc{A}'' \, | \varphi_{\mc{A}'}(v) = \varphi_{\mc{A}''}(w) \}.
	\ee
	The anchor is inherited from either $\mc{A}'$ or $\mc{A}''$ and its bracket is generated by the expression $[v \oplus v', w \oplus w']_{\mc{A}' \times_{\mc{A}} \mc{A}''} = [v,w]_{\mc{A}'} \oplus [v',w']_{\mc{A}''}$, for $v,w \in \Gamma(\mc{A}')$ and $v',w' \in \Gamma(\mc{A}'')$
\end{defn}
The projections onto $\mc{A}'$ and $\mc{A}''$ turn the fiber product $\mc{A}' \times_{\mc{A}} \mc{A}''$ into both an $\mc{A}'$- and an $\mc{A}''$-Lie algebroid. There is a natural isomorphism $\mc{A}' \times_{\mc{A}} \mc{A} \cong \mc{A}'$. If $\mc{A}$ has dense isomorphism locus, recall that $Z_\mc{A} = X \backslash X_\mc{A}$  is the degeneracy locus of $\mc{A}$, and the support of its divisor.
\begin{prop}\label{prop:fiberprodisolocus} Let $\mc{A}, \mc{A}' \to X$ be Lie algebroids of divisor-type with disjoint supports. Then $\rho_\mc{A}$ and $\rho_{\mc{A}'}$ are transverse, and $\mc{A} \times_{TX} \mc{A}'$ is of divisor-type with support $Z_\mc{A} \cup Z_{\mc{A}'}$.
\end{prop}
\bp As $Z_\mc{A}$ and $Z_{\mc{A}'}$ are disjoint, pointwise either $\rho_\mc{A}$ or $\rho_{\mc{A}'}$ is an isomorphism, so that $\rho_\mc{A}$ and $\rho_{\mc{A}'}$ are everywhere transverse. The rest of the statement readily follows.
\ep
\subsection{Lie subalgebroids}
\label{sec:liesubalgds}
In this section we consider subobjects of Lie algebroids. Let $\mc{A} \to X$ be a Lie algebroid and $(\varphi,f)\colon (\mc{B}, N) \hookrightarrow (\mc{A},X)$ be a vector subbundle supported on $N \subseteq X$. It is an \emph{anchored subbundle} of $\mc{A}$ supported on $N$ if $\varphi$ is an anchored morphism, in which case $\rho_{B} = \rho_{A}|_N$ restricted to $\mc{B}$. It is easy to see that if $(\rho_\mc{A}|_N)(\mc{B}) \subseteq TN$ then $\mc{B}$ carries a unique anchored bundle structure for which $\varphi$ is such an anchored morphism.

Next, consider the submodule $\Gamma(\mc{A},\mc{B}) := \{v \in \Gamma(\mc{A}) \, | \, v|_N \in \Gamma(\mc{B})\}$ of $\Gamma(\mc{A})$, which comes equipped with a surjective restriction map $\Gamma(\mc{A},\mc{B}) \to \Gamma(\mc{B})$ with kernel $\Gamma(\mc{A},0_N)$. Note that $\mc{B}$ is a \emph{Lie subalgebroid} of $\mc{A}$ supported on $N$ if $\varphi$ is a Lie algebroid morphism, in which case the restriction map $\Gamma(\mc{A},\mc{B}) \to \Gamma(\mc{B})$ is a Lie algebra homomorphism. The following proposition is noted in \cite[Proposition 2.15]{Meinrenken17} and in part in \cite[Remark 2.1.21]{Li13}.
\begin{prop}[\cite{Meinrenken17}]\label{prop:liesubalgeroid} Let $(\varphi,f)\colon (\mc{B},N) \hookrightarrow (\mc{A},X)$ be an anchored subbundle supported on $N$. Then the following are equivalent:
	\bi
	\item The submodule $\Gamma(\mc{A},\mc{B}) \subseteq \Gamma(\mc{A})$ is involutive;
	\item $\mc{B}$ carries a unique Lie algebroid structure for which $\varphi$ is a Lie algebroid morphism.
	\ei
\end{prop}
\bp Assume that $\Gamma(\mc{A},\mc{B})$ is involutive. We show that $\Gamma(\mc{A},0_N)$ is an ideal in $\Gamma(\mc{A},\mc{B})$, i.e.\
\be
	[\Gamma(\mc{A},0_N), \Gamma(\mc{A},\mc{B})]_\mc{A} \subseteq \Gamma(\mc{A},\mc{B}).
\ee
Let $v \in \Gamma(\mc{A},0_N)$. Then $v = \sum_i f_i v_i$ for some $v_i \in \Gamma(\mc{A})$ and $f_i \in I_N$. Given $w \in \Gamma(\mc{A},\mc{B})$, we have $\rho_\mc{A}(w)|_N \in \Gamma(TN)$. Consequently $\mc{L}_{\rho_\mc{A}(w)} f_i|_N = 0$ for all $i$, and the Leibniz rule gives that $[v,w]_\mc{A}|_N = 0$, as we have
\be
	[v, w]_\mc{A} = [\sum_i f_i v_i, w]_\mc{A} = \sum_i (f_i [v_i,w]_{\mc{A}} - \mc{L}_{\rho_\mc{A}(w)} f_i \cdot v_i).
\ee
Thus $\Gamma(\mc{B}) \cong \Gamma(\mc{A},\mc{B}) / \Gamma(\mc{A},0_N)$ inherits a Lie bracket, making $\mc{B}$ into a Lie algebroid, since the Leibniz rule for $\rho_{\mc{A}}$ implies the Leibniz rule is satisfied for $\rho_{\mc{B}}$. The converse is immediate.
\ep
\begin{rem}\label{rem:liesubalgd} Sometimes one includes in the definition, apart from the conditions that $(\rho_{A}|_N)(\mc{B}) \subseteq TN$ and involutivity of $\Gamma(\mc{A},\mc{B})$, the condition that $\Gamma(\mc{A},0_N)$ is an ideal in $\Gamma(\mc{A},\mc{B})$. \autoref{prop:liesubalgeroid} shows this latter condition is a consequence of the first two.
\end{rem}
For later reference we recall the following special examples of Lie subalgebroids.
\begin{defn} A Lie subalgebroid $(\varphi,f)\colon (\mc{B},X') \to (\mc{A},X)$ is \emph{wide} if $X' = X$ and $f = {\rm id}_X$.
\end{defn}
\begin{exa}\label{exa:regfoliation} Let $\mc{F}$ be a regular foliation on $X$ and $D = T \mc{F} \subseteq TX$ the associated involutive distribution. Then $\mc{A}_\mc{F} := D$ is a Lie subalgebroid of $TX$. Moreover, if $\mc{A} \to X$ is a regular Lie algebroid with injective anchor, then $\mc{A} = \mc{A}_\mc{F}$ for $\mc{F}$ associated to the distribution $D = \rho_\mc{A}(\mc{A})$.
\end{exa}
\begin{exa}\label{exa:adistribution} Let $\mc{A} \to X$ be a Lie algebroid with an \emph{$\mc{A}$-distribution}, i.e.\ a smooth subbundle $D_\mc{A} \subseteq \mc{A}$. Its image $D := \rho_\mc{A}(D_\mc{A}) \subseteq TX$, need not be a distribution. If $D_\mc{A} \subseteq \mc{A}$ is $\mc{A}$-involutive (closed under $[\cdot,\cdot]_\mc{A}$), then it is an $\mc{A}$-Lie subalgebroid with $\mc{A}$-anchor the natural inclusion.
\end{exa}
\subsection{$\mc{A}$-invariant submanifolds}
\label{sec:degenloci}
In this section we discuss $\mc{A}$-invariant submanifolds, i.e.\ when one can restrict a given Lie algebroid, as a Lie algebroid. We follow in part \cite{Pym13}.
\begin{defn} Let $\mc{A} \to X$ be a Lie algebroid. Then a submanifold $N \subseteq X$ is \emph{$\mc{A}$-invariant} if $\rho_\mc{A}$ is tangent to $N$, i.e.\ $\rho_\mc{A}|_N$ has image inside $TN$ (instead of merely $TX|_N$).
\end{defn}
For closed submanifolds $N \subseteq X$, a vector field $V \in \Gamma(TX)$ preserves its vanishing ideal $I_N$ if and only if it is tangent to $N$. Hence a closed submanifold $N$ is $\mc{A}$-invariant if and only if $I_N$ is preserved by $\mc{A}$, i.e.\ we have $\mc{L}_{\rho_\mc{A}(v)} I_N \subseteq I_N$ for all $v \in \Gamma(\mc{A})$. Thus $\mc{A}$-invariant submanifolds consist of orbits of $\mc{A}$, i.e.\ are unions of (opens of) leaves of the singular foliation given by $\rho_\mc{A}$.
\begin{prop}\label{prop:ainvariantrestr} Let $\mc{A}$ be a Lie algebroid and $N$ be an $\mc{A}$-invariant submanifold of $X$. Then $(\mc{A}|_N, \rho_\mc{A}|_N) \to N$ with restriction of $[\cdot,\cdot]_\mc{A}$ is a Lie subalgebroid of $\mc{A}$ along $N$.
\end{prop}
\bp Using \autoref{prop:liesubalgeroid}, because $\mc{A}|_N \to TN$ is an anchored bundle due to $\mc{A}$-invariance of $N$, we need only show that $\Gamma(\mc{A},\mc{A}|_N)$ is involutive. However, this is immediate.
\ep
\begin{rem}\label{rem:inverseimagesmooth} The proof of \autoref{prop:ainvariantrestr} shows %that
more generally, given
% a Lie algebroid $\mc{A} \to X$ and 
a submanifold $\iota\colon N \hookrightarrow X$, that $\iota^! \mc{A} := \rho_\mc{A}^{-1}(TN)$ is a Lie subalgebroid of $\mc{A}$ along $N$ if it is a smooth subbundle of $\mc{A}|_N$ (see \cite[Proposition 2.17]{Meinrenken17}). Note here that $\iota^! \mc{A} = \mc{A}|_N$ if and only if $N$ is $\mc{A}$-invariant.
\end{rem}
The set of $\mc{A}$-invariant closed submanifolds of $X$ admits unions and intersections. We say a subset of $X$ is \emph{$\mc{A}$-invariant} if it is a union of $\mc{A}$-invariant submanifolds.
\begin{lem}\label{lem:ainvariance} Let $N$ and $N'$ be $\mc{A}$-invariant closed submanifolds of $X$. Then their union $N \cup N'$ and intersection $N \cap N'$ are $\mc{A}$-invariant closed subsets.
\end{lem}
\bp The sets $N \cup N'$ and $N \cap N'$ have vanishing ideals given by $I_{N \cup N'} = I_N \cap I_{N'} = I_N \cdot I_{N'}$ and $I_{N \cap N'} = I_N + I_{N'}$. Using the product rule $\mc{L}_{\rho_\mc{A}(v)}(f g) = (\mc{L}_{\rho_\mc{A}(v)}f) \cdot g + f \cdot (\mc{L}_{\rho_{\mc{A}}(v)}g)$, we conclude that $I_{N \cup N'}$ is $\mc{A}$-invariant. The $\mc{A}$-invariance of $I_{N \cap N'}$ is immediate by additivity. Alternatively, one can note that $\mc{A}$-invariant submanifolds are precisely unions of $\mc{A}$-orbits.
\ep
\begin{rem}\label{rem:preservingideals} The proof of \autoref{lem:ainvariance} shows that more generally, given a collection of vector fields $\mc{F} \subseteq \Gamma(TX)$ preserving two ideals $I$, $I'$, i.e.\ such that $\mc{L}_v I \subseteq I$ and $\mc{L}_v I' \subseteq I'$ for all $v \in \mc{F}$, the collection $\mc{F}$ also preserves $I + I'$ and $I \cap I'$ (see \cite[Proposition 2.1.10]{Pym13}). In \autoref{lem:ainvariance} this is applied to the collection $\mc{F} := \rho_{A}(\Gamma(\mc{A}))$ and ideals $I := I_N$, $I' := I_{N'}$.
\end{rem}
When considering residue maps associated to $\mc{A}$-invariant submanifolds in \autoref{sec:residuemaps}, we will need the following terminology, concerning properties of the restricted anchor map.
\begin{defn} Let $\mc{A} \to X$ be a Lie algebroid and $N \subseteq X$ an $\mc{A}$-invariant submanifold. Then $N$ is \emph{projective} if the restriction $\mc{A}|_N$ is projective, and \emph{transitive} if $\mc{A}|_N$ is transitive.
\end{defn}
Certainly any transitive $\mc{A}$-invariant submanifold is in particular projective. Any orbit $\mc{O}$ of $\mc{A}$ is a transitive $\mc{A}$-invariant submanifold. Further examples will be given later in \autoref{sec:residues}.

Given a projective $\mc{A}$-invariant submanifold $N \subseteq X$, by the definition of projectivity, the Serre--Swan theorem provides a Lie algebroid $\mc{B} \to N$ satisfying $\Gamma(\mc{B}) \cong (\rho_A|_N)(\Gamma(\mc{A}|_N))$. There is an induced Lie algebroid morphism $\wt{\rho}_\mc{A} \colon \mc{A}|_N \to \mc{B}$ obtained from the natural map on sections. This fits in a short exact sequence of Lie algebroids over $N$, namely
\be
	0 \to \ker \wt{\rho}_\mc{A} \to \mc{A}|_N \to \mc{B} \to 0.
\ee
As we did in \autoref{sec:liealgebroids}, we call $\ker \wt{\rho}_\mc{A}$ the (bundle of) \emph{germinal isotropy} of $\mc{A}|_N$, after \cite{AndroulidakisZambon17}. 
\subsection{Degeneracy loci}
In this section we define the degeneracy loci of $\mc{A}$ and show that they are $\mc{A}$-invariant subsets. Given an integer $k \geq 0$, the anchor of $\mc{A}$ defines a map $\wedge^k \rho_\mc{A}\colon \wedge^k \mc{A} \to \wedge^k TX$. This induces a map $\wedge^k \rho_\mc{A}\colon \Gamma(\wedge^k \mc{A} \otimes \wedge^k T^*X) \to C^\infty(X)$, whose image is an ideal $I_{\mc{A},k}$.
\begin{defn} Let $\mc{A} \to X$ be a Lie algebroid. The $k$th \emph{degeneracy locus} of $\mc{A}$ is the subspace $X_{\mc{A},k} \subseteq X$ defined by the ideal $I_{\mc{A},k+1} \subseteq C^\infty(X)$, where the rank of $\rho_\mc{A}$ is $k$ or less.
\end{defn}
\begin{prop}[{\cite[Proposition 2.2.9]{Pym13}}]\label{prop:degenlocus} Let $\mc{A} \to X$ be a Lie algebroid and $k \geq 0$. Then $X_{\mc{A},k}$ is an $\mc{A}$-invariant submanifold whenever it is a smooth submanifold.
\end{prop}
\bp Assume that $X_{\mc{A},k}$ is smooth and denote the map $\wedge^{k+1} \rho_\mc{A}$ by $\varphi_{k}\colon \wedge^{k+1} \mc{A} \to \wedge^{k+1} TX$. Let $v \in \Gamma(\mc{A})$ and $f \in I_{\mc{A},k}$. Then by definition there exists $\xi \in \Gamma(\wedge^{k+1} \mc{A})$ and $\omega \in \Omega^{k+1}(X)$ such that $f = \langle \varphi_k(\xi), \omega \rangle$. Using the compatibility of the bracket and anchor we compute that
\begin{align*}
	\mc{L}_{\rho_\mc{A}(v)} f &= \mc{L}_{\rho_\mc{A}(v)} \langle \varphi_k(\xi), \omega \rangle = \langle \mc{L}_{\rho_\mc{A}(v)} \varphi_k(\xi), \omega \rangle + \langle \varphi_k(\xi), \mc{L}_{\rho_\mc{A}(v)} \omega \rangle\\
	&= \langle \varphi_k(\mc{L}_{v} \xi), \omega\rangle + \langle \varphi_k(\xi), \mc{L}_{\rho_\mc{A}(v)} \omega \rangle \in I_{\mc{A},k}.
\end{align*}
As $I_{\mc{A},k}$ is preserved under the action of $\mc{A}$, we conclude that $X_{\mc{A},k}$ is $\mc{A}$-invariant.
\ep
\begin{rem}\label{exa:denseisolocus} Let $\mc{A}^n \to X^n$ be a Lie algebroid of divisor-type. Then $X \backslash X_\mc{A} = Z_\mc{A} = X_{\mc{A},n-1}$. In fact, we have an equality of ideals $I_\mc{A} = I_{\mc{A},n}$. Both $X_\mc{A}$ and $Z_\mc{A}$ are $\mc{A}$-invariant subsets.
\end{rem}
\subsection{Singular Lie subalgebroids}
\label{sec:singularliealgebroids}
In this section we take a slight step into the world of singular foliations. While those underlying the Lie algebroids of interest to us are relatively mild, it is useful to relate our discussion to the language in existing literature, particularly \cite{AndroulidakisSkandalis09, AndroulidakisZambon13, AndroulidakisZambon16, AndroulidakisZambon17, LaurentGengouxLavauStrobl18, Zambon18}. To extend the results in this paper to more singular Poisson structures, one is forced to work with sheaves or modules which no longer have an underlying bundle. Because we defer this to future work, we will be somewhat brief here.
\begin{rem} As a general comment, projective $C^\infty(X)$-modules are not necessarily locally finitely-generated, i.e.\ need not necessarily correspond to section modules of vector bundles. However, our use of the term `projective' will presuppose that the module in question is locally finitely-generated. Indeed, for us such modules are often obtained as images of bundle maps.
\end{rem}
Let $X$ be a manifold and $\mc{C}^\infty_X$ be the sheaf of smooth functions on $X$, with $C^\infty(X) = \mc{C}^\infty_X(X)$. Throughout this section, a \emph{module} will refer to a $C^\infty(X)$-module without further mention.

Let $E \to X$ be a vector bundle and $\mc{E}' \subseteq \mc{E} := \Gamma(E)$ be a submodule. Given a smooth map $f \colon X' \to X$, let $f^* E \to X'$ be the pullback bundle. The \emph{pullback module} $f^* \mc{E}'$ of $\mc{E}'$ along $f$ is
\be
	f^* \mc{E}' := \langle g f^* \xi \, | \, g \in C^\infty(X'), \xi \in \mc{E}' \rangle \subseteq \Gamma(X'; f^*E).
\ee
If $\iota_{X'}\colon X \hookrightarrow X$ is a submanifold, we call $\iota_{X'}^* \mc{E}'$ the \emph{restriction} of $\mc{E}'$ to $X'$.
A module is \emph{finitely generated} if there exist global sections $\xi_i$ of $\mc{E}$ such that $\mc{E} = C^\infty(X) \langle \xi_i \rangle$, and \emph{locally finitely generated} if there exists an open cover $\{U_i\}$ of $X$ such that each restriction $i_{U_i}^*\mc{E}$ is finitely generated. Call a module $\mc{E}$ of \emph{vector bundle-type} when $\mc{E} \cong \Gamma(E)$ for a vector bundle $E \to X$.
\begin{defn} A module $\mc{E}$ is \emph{projective} if it has a global basis of generating sections, or, equivalently, if it is a direct summand of a finitely generated free module. 
\end{defn}
\begin{rem} There is a one-to-one correspondence between locally finitely generated projective modules $\mc{E}$, and locally free sheaves of $\mc{C}^\infty_X$-modules $\{U \mapsto \mc{E}_U\}$. This extends to singular foliations (adding involutivity and possibly dropping projectiveness), giving a one-to-one correspondence between singular foliations and involutive sheaves of locally finitely generated submodules of $\mf{X}_X$ \cite{Wang17}. Adding projectivity (or local freeness), we obtain that projective singular foliations are exactly locally free sheaves of involutive submodules of $\mf{X}_X$.
\end{rem}
The \emph{Serre--Swan correspondence} \cite{Serre55, Swan62} gives an equivalence of categories between finitely generated projective modules $\mc{E}$, and vector bundles $E \to X$ via the assignment $\mc{E} \mapsto \Gamma(E)$. In other words, locally finitely generated projective modules are of vector bundle-type.

Focusing on involutivity, using a given Lie algebroid $\mc{A} \to X$, results in the following (we give the definition in terms of sheaves, instead of using compactly-supported sections).
\begin{defn}[c.f.\ \cite{Zambon18}]\label{defn:singularliealgebroid} A \emph{singular subbundle} of $\mc{A}$ is a sheaf of locally finitely generated submodules of its sheaf of sections $\Gamma(\mc{A})$.  A \emph{singular Lie subalgebroid} of $\mc{A}$ is an involutive singular subbundle. We also use the associated module of global sections interchangeably.
\end{defn}
\begin{exa}	A wide Lie subalgebroid $\mc{A}' \subseteq \mc{A}$ induces a singular Lie subalgebroid $\Gamma(\mc{A}')$.
\end{exa}
\begin{exa} \autoref{prop:liesubalgeroid} shows that given a Lie subalgebroid $\mc{B}$, the module $\Gamma(\mc{A},\mc{B})$ is a singular Lie subalgebroid of $\mc{A}$. It is not necessarily projective (see also \autoref{sec:elementarymodifications}).
\end{exa}

A projective singular Lie subalgebroid of $\mc{A}$ defines an $\mc{A}$-Lie algebroid $\mc{A}'$ along with a Lie algebroid morphism $\mc{A}' \to \mc{A}$. Note however that $\mc{A}'$ is generally \emph{not} a Lie subalgebroid of $\mc{A}$.

Singular Lie subalgebroids of $TX$ are called \emph{singular foliations} in the sense of Androulidakis--Skandalis \cite{AndroulidakisSkandalis09}. Moreover, projective singular foliations $\mc{F}$ are also called \emph{almost-regular} \cite{Debord01}, and give rise to almost-injective Lie algebroids $\mc{A}_\mc{F}$ with $\Gamma(\mc{A}_\mc{F}) \cong \mc{F}$. Note that underlying any singular Lie subalgebroid $\mc{F}_\mc{A}$ of $\mc{A}$ there is a singular foliation $\mc{F} := \rho_{A}(\mc{F}_\mc{A})$. However, in general there is no relation between projectivity of the modules $\mc{F}_\mc{A} \subseteq \Gamma(\mc{A})$ and $\mc{F} \subseteq \Gamma(TX)$.
\subsubsection{On projectivity}
Each submodule $\mc{E}'$ of a projective module $\mc{E}$ has an evaluation map ${\rm ev}\colon X \times \mc{E}' \to E$ defined by ${\rm ev}(x,\xi) := \xi(x)$. Define ${\rm ev}_x\colon \mc{E}' \to E_x$ by ${\rm ev}_x := {\rm ev}(x, \cdot)$ for $x \in X$. We will denote by $E_x' := {\rm ev}_x(\mc{E}')$ the image of this map, which is a vector subspace of the fiber $E_x$ of $E$. Given a point $x \in X$, let $I_x \subseteq C^\infty(X)$ be its vanishing ideal. Moreover, given a submodule $\mc{E}'$ of $\mc{E} = \Gamma(E)$, consider the quotient vector space $\mc{E}'_x := \mc{E}' / I_x \mc{E}'$ of $\mc{E}'$ by $I_x$, called the \emph{fiber} of $\mc{E}'$ at $x$. Note that $\mc{E}'_x = E'_x$ when $\mc{E}'$ is of vector bundle-type. We have:
\begin{prop} Let $\mc{E}' \subseteq \mc{E}$ be a submodule on $X$ where $\mc{E} = \Gamma(E)$ is of vector bundle-type. Then $x \mapsto \dim E'_x$ is lower semi-continuous while $x \mapsto \dim \mc{E}'_x$ is upper semi-continuous.
\end{prop}
When $\mc{E}$ is locally finitely generated, the vector spaces $\mc{E}_x$ are finite-dimensional for all $x \in X$ (see \cite{AndroulidakisSkandalis09}). The evaluation map ${\rm ev}_x$ vanishes on $I_x \mc{E}'$, so that it descends to a surjective homomorphism $\wt{\rm ev}_x\colon \mc{E}'_x \to E'_x$. When $\mc{E}$ is a Lie module and $\mc{E}' \subseteq \mc{E}$ an involutive submodule, the kernel of ${\rm ev}_x$ is a Lie subalgebra of $\mc{E}'$, and $I_x \mc{E}'$ is an ideal in this subalgebra. This results in the conclusion that $\mf{a}_x := \ker {\rm ev}_x / I_x \mc{E}' = \ker \wt{\rm ev}_x$ is a Lie algebra, called the \emph{isotropy} of $\mc{E}'$ at $x \in X$ when $\mc{E} = \Gamma(TX)$. There is no induced bracket on $\mc{E}'_x$. There is a short exact sequence
\be
0 \to \mf{a}_x \to \mc{E}'_x \stackrel{\wt{\rm ev}_x}{\to} E_x' \to 0.
\ee
\begin{prop}[{\cite[Lemma 1.3]{AndroulidakisZambon13}}] Let $\mc{E}' \subseteq \mc{E}$ be a submodule on $X$ where $\mc{E} = \Gamma(E)$ is of vector bundle-type. Then $\dim \mf{g}_x$ is bounded from above by ${\rm rank}(E)$.
\end{prop}
The Serre--Swan correspondence implies that a locally finitely-generated module $\mc{F}$ is projective if and only if the dimension of $\mc{F}_x$ is constant for all $x \in X$. We can pinpoint the failure of projectivity more precisely, following \cite{AndroulidakisZambon17}. Let $\mc{F}_\mc{A}$ be a singular Lie subalgebroid of $\mc{A}$. Then the module $\Gamma(\mc{A})$ is projective, and $\mc{F}_\mc{A}$ is a submodule of it. As such the discussion above leads to a short exact sequence for each point $x \in X$, namely
\be
	0 \to \mf{a}_x \to \mc{F}_{\mc{A},x} \to \mc{A}_{\mc{F}_\mc{A},x} \to 0,
\ee
where $\mc{A}_{\mc{F}_\mc{A},x} = {\rm ev}_x(\mc{F}_\mc{A}) \subseteq \mc{A}_x$ is the evaluation, and $\mc{F}_{\mc{A},x} = \mc{F}_\mc{A} / I_x \mc{F}_\mc{A}$. Now assume that $\mc{F}_\mc{A}$ is induced by an $\mc{A}$-Lie algebroid. In other words, assume that  we have $\varphi_{\mc{A}}\colon \mc{A}' \to \mc{A}$ with $\mc{F}_\mc{A} = \varphi_{\mc{A}}(\Gamma(\mc{A}')) \subseteq \Gamma(\mc{A})$. We then obtain from $\varphi_{\mc{A}}$ another short exact sequence, namely
\be
	0 \to \mf{g}_x \to \mc{A}'_x \to \mc{A}_{\mc{F}_\mc{A},x} \to 0.
\ee
The $\mc{A}$-anchor $\varphi_{\mc{A}}\colon \Gamma(\mc{A}') \to \Gamma(\mc{A})$ on sections gives another short exact sequence:
\be
	0 \to \mf{h}_x \to \mc{A}'_x \to \mc{F}_{\mc{A},x} \to 0.
\ee
We call the kernel $\mf{h}_x$ the \emph{germinal $\mc{A}$-isotropy} of $\mc{F}_\mc{A}$ at $x \in X$. If the module $\mc{F}_\mc{A}$ is projective, so that $\mc{F}_\mc{A} \cong \Gamma(\mc{B})$ for $\mc{B} \to X$, there is an induced Lie algebroid morphism $\wt{\varphi}_\mc{A}\colon \mc{A}' \to \mc{B}$, and then $\mf{h}_x = \ker(\wt{\varphi}_{\mc{A},x}) \subseteq \mc{A}'_x$ as in \autoref{defn:germinalisotropy} for $\mc{A} = TX$. The above sequences are related, as discussed in \cite[Section 1.3]{AndroulidakisZambon17} for the case of $\mc{A} = TX$. Indeed, we have a short exact sequence
\be
	0 \to \mf{h}_x \to \mf{g}_x \to \mf{a}_x \to 0.
\ee
In conclusion, projectivity of $\mc{F}_\mc{A}$ is equivalent to $\dim(\mf{h}_x)$ being constant for all $x \in X$. This criterion for projectivity also follows from the fact that the above Lie algebroid morphism $\wt{\varphi}_\mc{A}$ is surjective by construction, so that it induces a Lie algebroid extension sequence
\be
	0 \to \ker(\wt{\varphi}_\mc{A}) \to \mc{A}' \to \mc{B} \to 0.
\ee
In \cite{LaurentGengouxLavauStrobl18} the authors discuss the notion of a \emph{projective resolution} of a singular foliation $\mc{F}$, which is an exact sequence of sheaves of $C^\infty(X)$-modules:
\be
	\dots \to \Gamma(E_{-2}) \to \Gamma(E_{-1}) \to \mc{F} \to 0,
\ee
where $E_{-i}$ for $i \geq 1$ are vector bundles over $X$, with anchor map $E_{-1} \to TX$ with image $\mc{F}$.
It is clear that this definition makes sense for singular Lie subalgebroids as well upon replacing $TX$ by $\mc{A}$ (giving \emph{projective $\mc{A}$-resolutions}). The \emph{projective $\mc{A}$-dimension} of such a module (or sheaf of submodules) is then defined as the minimal length of a projective resolution, minus one. Thus projective singular Lie subalgebroids are those with zero projective $\mc{A}$-dimension. In future work we will treat situations of positive projective dimension.
\subsection{Ideal Lie algebroids}
\label{sec:idealliealgebroids}
In this section we discuss what we call ideal Lie algebroids, which are Lie algebroids defined using a given ideal of functions. Of particular importance are those created using divisor ideals. Let $\mc{A} \to X$ be a Lie algebroid and let $I \subseteq C^\infty(X)$ be an ideal. Consider $I \cdot \Gamma(\mc{A}) \subseteq \Gamma(\mc{A})$ consisting of all finite sums of products of elements of $I$ and $\Gamma(\mc{A})$, and let $\Gamma(\mc{A})_I := \{v \in \Gamma(\mc{A}) \, | \, \mc{L}_{\rho_{A}(v)} I \subseteq I\}$ be those sections preserving $I$. Both are $C^\infty(X)$-submodules of $\Gamma(\mc{A})$. In this section we sometimes shorten $\mc{L}_{\rho_{\mc{A}}(v)}$ to $\mc{L}_v$.
\begin{lem}\label{lem:liesubalgs} Let $I \subseteq C^\infty(X)$ be an ideal. Then $I \cdot \Gamma(\mc{A}) \subseteq \Gamma(\mc{A})_I \subseteq \Gamma(\mc{A})$ are Lie subalgebras.
\end{lem}
\bp Let $v, w \in \Gamma(\mc{A})_I$ and $f \in I$. Then $\mc{L}_v f = g$ and $\mc{L}_w f = h$ for $g, h \in I$. Hence
\be
	\mc{L}_{[v,w]} f = (\mc{L}_v \mc{L}_w - \mc{L}_w \mc{L}_v) f = \mc{L}_v h - \mc{L}_w g.
\ee
As $g, h \in I$ we have $\mc{L}_{[v,w]} f \in I$ for all $f \in I$, hence $[v,w] \in \Gamma(\mc{A})_I$.

Let $v, w \in I \cdot \Gamma(\mc{A})$ with $v = \sum_i f_i v_i$ and $w = \sum_j f'_j w_j$ for some $f_i, f'_j \in I$ and $v_i, w_j \in \Gamma(\mc{A})$. Then we have $[v,w] \in I \cdot \Gamma(\mc{A})$ by the Leibniz rule as each term starts with an element of $I$,
\be 
[v,w] = \textstyle{\sum}_{i,j} [f_i v_i, f'_j w_j] = \textstyle{\sum}_{i,j} \left( f_i f'_j [v_i, w_j] + f_i (\mc{L}_{v_i} f'_j) \cdot w_j - f'_j (\mc{L}_{w_j} f_i) \cdot v_i\right).
\ee
Let $v \in I \cdot \Gamma(\mc{A})$. Then $v = \sum_i f_i v_i$ for $f_i \in I$ and $v_i \in \Gamma(\mc{A})$. For $f \in I$ we have $\mc{L}_v f = \sum_i f_i (\mc{L}_{v_i} f)$. As $f_i \in I$ for all $i$ and $I$ is an ideal, we conclude that $\mc{L}_v f \in I$ so that $v \in \Gamma(\mc{A})_I$ as desired.
\ep
It is immediate that if $I$ is locally finitely generated, then both $I \cdot \Gamma(\mc{A})$ and $\Gamma(\mc{A})_I$ are as well. By \autoref{lem:liesubalgs} they are then both singular Lie subalgebroids of $\mc{A}$. Consequently, if they are also projective, they define $\mc{A}$-Lie algebroids by the Serre--Swan correspondence. Due to \autoref{lem:liesubalgs}, there is then a sequence of Lie algebroid morphisms as follows:
\be
	I \cdot \mc{A} \to \mc{A}_I \to \mc{A} \to TX.
\ee
We will mostly focus on the projective case, i.e.\ when $\Gamma(\mc{A})_I$ and $I\cdot \Gamma(\mc{A})$ are projective.
\begin{defn} Let $\mc{A} \to X$ be a Lie algebroid and $I \subseteq C^\infty(X)$ be an ideal. The primary and secondary \emph{ideal Lie algebroids} associated to $(\mc{A},I)$ are the Lie algebroids $\mc{A}_I \to X$ such that $\Gamma(\mc{A}_I) \cong \Gamma(\mc{A})_I$, and $I \cdot \mc{A} \to X$ such that $\Gamma(I \cdot \mc{A}) \cong I \cdot \Gamma(\mc{A})$.
\end{defn}
Both ideal Lie algebroids are unique -- provided they exist -- and their $\mc{A}$-anchors are the natural inclusions on sections. However, their $\mc{A}$-anchors are not injective when viewed as bundle maps, unless $I$ is trivial. The isomorphism locus of the $\mc{A}$-anchor of $\mc{A}_I$, and of $I \cdot \mc{A}$, as well as the morphism between $I\cdot \mc{A}$ and $\mc{A}_I$, is given by $X \backslash Z_I$, the complement of the support of $I$. In other words, these three Lie algebroid morphisms are of divisor-type supported on $Z_I$, which in particular means that the ideal Lie algebroids $\mc{A}_I$ and $I \cdot \mc{A}$ are both of $\mc{A}$-divisor-type.

Our main interest lies in the case where $\mc{A} = TX$ and $I$ is a divisor ideal. Then both ideal Lie algebroids are of divisor-type, with divisors ${\rm div}(I \cdot TX)$ and ${\rm div}(TX_I)$ supported on $Z_I$, and $I_{I \cdot TX} = I^n$ if $\dim X = n$ (and similarly using $\mc{A}$). We introduce the following terminology.
\begin{defn}\label{defn:divprojstandard} A divisor ideal $I \subseteq C^\infty(X)$ is said to be \emph{projective} if $\Gamma(TX)_I$ is projective. A divisor ideal is \emph{standard} if it is projective and the divisor ideal of $TX_I$ satisfies $I_{TX_I} = I$.
\end{defn}
Not every divisor ideal is projective, nor is every projective divisor ideal standard. Projectivity of $I \cdot \Gamma(A)$ is guaranteed for divisor ideals, as is also noted in \cite[Lemma 1.11]{Nistor15}.
\begin{prop}\label{prop:dividealalgebroid} Given a Lie algebroid $\mc{A}^n \to X$ and divisor ideal $I$, $I \cdot \Gamma(\mc{A})$ is projective.
\end{prop}
\bp 
There is a short exact sequence of sheaves $0 \to I \cdot \Gamma(\mc{A}) \to \Gamma(\mc{A}) \to \Gamma(\mc{A}) / I \cdot \Gamma(\mc{A}) \to 0$. As both $\Gamma(\mc{A})$ and $I$ are locally free, so is $\Gamma(\mc{A}) / I \cdot \Gamma(\mc{A}) \cong (I^{-1})^n$ (the inverse of $I$ as a sheaf). The same then holds for $I \cdot \Gamma(\mc{A})$ as any short exact sequence in a projective module splits.
\ep
When $I$ is not necessarily standard, we expect the following relation between $I$ and $I_{TX_I}$.
\begin{conj} Let $I \subseteq C^\infty(X)$ be a projective divisor ideal. Then $I$ is divisible by $I_{TX_I}$.
\end{conj}
Note that in general, given a Lie algebroid morphism $(\varphi,{\rm id}_X)\colon \mc{A} \to \mc{A}'$ and a divisor ideal $I$, there are induced maps $\varphi\colon \Gamma(\mc{A})_I \to \Gamma(\mc{A}')_I$ and $\varphi\colon I \cdot \Gamma(\mc{A}) \to I \cdot \Gamma(\mc{A}')$. In particular this holds upon taking $\varphi = \rho_{A}$ and $\mc{A}' = TX$. The following is based on \cite[Remark 135]{Frejlich11}.
\begin{prop}\label{prop:transLAprojective} Let $\mc{A} \to X$ be a transitive Lie algebroid and $I$ a projective divisor ideal. Then $\Gamma(\mc{A})_I$ is a projective module, and the associated bundle $\mc{A}_I$ is a $TX_I$-Lie algebroid. Moreover, $\mc{A}_I$ is isomorphic to the fiber product $\mc{A} \times_{TX} TX_I$.
\end{prop}
\bp As $\mc{A}$ is transitive, there is a short exact sequence $0 \to K \to \mc{A} \to TX \to 0$ where $K = \ker \rho_{A}$ is totally intransitive. We can always split this sequence using an Ehresmann connection giving an isomorphism $\Gamma(\mc{A}) \cong \Gamma(TX) \oplus \Gamma(K)$. Using this splitting it follows that $\Gamma(\mc{A})_I \cong \Gamma(TX)_I \oplus \Gamma(K)$, showing projectivity of $\Gamma(\mc{A})_I$ and implying that there is an isomorphism $\mc{A}_I \cong \mc{A} \times_{TX} TX_I$, from which the fact that $\mc{A}_I$ is a $TX_I$-Lie algebroid follows.
\ep
Given a surjective Lie algebroid morphism $(\varphi,{\rm id}_X)\colon \mc{A} \to \mc{A}'$ and divisor ideal $I$, if the module $\Gamma(\mc{A}')_I$ is projective, the proof of \autoref{prop:transLAprojective} shows that $\Gamma(\mc{A})_I$ is also projective, and that $\mc{A}_I \cong \mc{A} \times_{\mc{A}'} \mc{A}'_I$. Sadly the transitivity of \autoref{prop:transLAprojective} is never satisfied for Lie algebroids of nontrivial divisor-type, so that this inheritance property does not apply to them.
We next study to what extent $\Gamma(\mc{A})_I$ depends on the ideal $I$. Recall that the \emph{radical} of $I$ is defined as $\sqrt{I} := \{f \in C^\infty(X) \, | \, f^n \in I \text{ for some } n \in \N\}$. Then $I \subseteq \sqrt{I}$ and $Z_I = Z_{\sqrt{I}}$.
\begin{prop}\label{prop:locprincidealmodule} Let $I \subseteq C^\infty(X)$ be a locally principal ideal. Then $\Gamma(\mc{A})_I = \Gamma(\mc{A})_{\sqrt{I}}$.
\end{prop}
\bp This can be verified locally. Let $f \in \sqrt{I}$ be such that $\sqrt{I} = \langle f \rangle$, so that $I = \langle f^n \rangle$ for some $n \in \N$. Take $V \in \Gamma(\mc{A})_I$, so that $\mc{L}_V f^n = g \in I$. Then $g = g' f^n$ for some function $g'$, hence $n f^{n-1} \mc{L}_V f = g' f^n$, so that $\mc{L}_V f = \frac1n g' f$, whence $\mc{L}_V f \in \sqrt{I}$. As $f$ generates $\sqrt{I}$, we conclude that $V \in \Gamma(\mc{A})_{\sqrt{I}}$, and thus $\Gamma(\mc{A})_I \subseteq \Gamma(\mc{A})_{\sqrt{I}}$. For the other inclusion, take $W \in \Gamma(\mc{A})_{\sqrt{I}}$ and $h \in I$. Then $h = h' f^n$ for some function $h'$, hence
\be
	\mc{L}_W h = \mc{L}_W h' f^n = n h' f^{n-1} \mc{L}_W f + f^n \mc{L}_W h'.
\ee
As $\mc{L}_W f \in \sqrt{I}$, write $\mc{L}_W f = f' f$ for some function $f'$. Then $\mc{L}_W h = (n h' f' + \mc{L}_W h') f^n \in I$, so that $W \in \Gamma(\mc{A})_I$. Thus $\Gamma(\mc{A})_{\sqrt{I}} \subseteq \Gamma(\mc{A})_I$, and equality follows.
\ep
As a consequence of \autoref{prop:locprincidealmodule}, the primary ideal Lie algebroids $\mc{A}_I$ associated to divisor ideals $I$ are only sensitive to $I$ up to taking powers; this is not true for the secondary ideal Lie algebroids $I \cdot \mc{A}$.
\begin{rem}\label{rem:nonstandard} There are many divisor ideals which are not standard. For example, using \autoref{prop:locprincidealmodule} we obtain that $\Gamma(\mc{A})_{I^k} = \Gamma(\mc{A})_I$ for all $k > 1$. Consequently, if $I$ is a projective divisor ideal, we have that $TX_{I^k} = TX_I$. Hence, if $I$ is moreover standard, we get that $I_{TX_{I^k}} = I_{TX_I} = I$, which is not equal to $I^k$ unless $I$ is trivial, i.e.\ $I^k$ is not standard.
\end{rem}
\begin{rem}\label{rem:preservingradical} Analogously to \autoref{rem:preservingideals}, a slight adaptation of the proof of \autoref{prop:locprincidealmodule} shows that more generally, given a collection of vector fields $\mc{F} \subseteq \Gamma(\mc{A})$ preserving an ideal $I$, i.e.\ such that $\mc{L}_v I \subseteq I$ for all $v \in \mc{F}$, the collection $\mc{F}$ also preserves $\sqrt{I}$. This is a purely algebraic statement and can be found in \cite[Lemma 3.3.2]{Dixmier96}.
\end{rem}
As $Z_I$ is the degeneracy locus of both $I\cdot \mc{A}$ and $\mc{A}_I$, \autoref{prop:degenlocus} gives the following.
\begin{prop}\label{prop:ziaiinvariant} Let $I \subseteq C^\infty(X)$ be a divisor ideal such that $\mc{A}_I$ exists. Then $Z_I$ is $I \cdot \mc{A}$- and $\mc{A}_I$-invariant if it is smooth, i.e.\ $I \cdot \Gamma(\mc{A})$ and $\Gamma(\mc{A}_I)$ have anchor images tangent to $Z_I$.
\end{prop}
Let $I$ and $I'$ be two divisor ideals. If $\mc{A}_{I \cdot I'}$ exists, then so do $\mc{A}_I$ and $\mc{A}_{I'}$, and we see that $\mc{A}_{I \cdot I'} = \mc{A}_I \times_{TX} \mc{A}_{I'}$ is their fiber product. Moreover, $\mc{A}_{I \cdot I'}$ is both an $\mc{A}_I$- and an $\mc{A}_{I'}$-Lie algebroid (again, if it exists). In favorable cases, the $\mc{A}$-divisor of $\mc{A}_{I \cdot I'}$ is the product of divisors of $\mc{A}_I$ and $\mc{A}_{I'}$. For example this holds if $\mc{A} = TX$ and $I \cdot I'$ is standard.
\begin{rem} While $\Gamma(\mc{A})_{I \cdot I'}$ being projective implies both $\Gamma(\mc{A})_I$ and $\Gamma(\mc{A})_{I'}$ are, the converse is not immediate if $Z_I$ and $Z_{I'}$ are not disjoint. For example, one can take $\mc{A} = TX$, $I = I_Z$ and $I' = I_{Z'}$ for hypersurfaces $Z,Z' \subseteq X$ which have nontransverse intersection.
\end{rem}
\subsubsection{Examples}
\label{exa:ideallalgebroids}
In this section we discuss some examples of ideal Lie algebroids constructed using the divisor examples of \autoref{sec:divexamples}. In general it is quite difficult to determine when a divisor ideal is projective, comparable to finding free divisors in algebraic geometry.
\begin{exa}[Trivial divisors] Let $I = C^\infty(X)$ be the trivial divisor ideal. Then
\be
	I \cdot \Gamma(TX) = \Gamma(TX)_I = \Gamma(TX),
\ee
so that the ideal Lie algebroids satisfy $I \cdot TX = TX_I = TX$. The same holds for any Lie algebroid $\mc{A}$, in that $I \cdot \mc{A} = \mc{A}_I = \mc{A}$. In particular $I$ is both projective and standard. 
\end{exa}
\begin{exa}[Log divisors]\label{exa:zvx} Let $Z \subseteq X$ be a closed submanifold with vanishing ideal $I_Z$. Then $\Gamma(TX)_{I_Z}$ consists of all vector fields tangent to $Z$. It is projective if and only if $Z$ has codimension one. In this case $I_Z$ is a log divisor ideal and $\mc{A}_Z := TX_{I_Z}$ is the \emph{log-tangent bundle}\footnote{It is also called the \emph{$b$-tangent bundle} ${}^b TX$, introduced by Melrose when $Z = \partial X$ (see \cite{Melrose93} and \cite{NestTsygan96, GuilleminMirandaPires14}).} $TX(-\log Z)$. Further, $I_Z \cdot \Gamma(TX)$ consists of all vector fields which vanish on $Z$. As with $\Gamma(TX)_{I_Z}$, it is projective if and only if $Z$ has codimension one, in which case $\mc{B}_Z := I_Z \cdot TX$ is the \emph{zero tangent bundle} (also denoted by ${}^0 TX$). Log divisors are standard (\cite[Proposition 3.4.6]{Klaasse17}), as locally we have
\be
	\Gamma(\mc{A}_Z) = \Gamma(TX)_{I_Z} = \langle z \partial_z, \partial_{x_i} \rangle \quad \text{and} \quad \Gamma(\mc{B}_Z) = I_Z \cdot \Gamma(TX) = \langle z \partial_z, z \partial_{x_i} \rangle.
\ee
\end{exa}
\autoref{exa:zvx} extends to normal-crossing and star log divisors $I_{\underline{Z}}$. In both cases, there is a projective module $\Gamma(TX)_{I_{\underline{Z}}}$ defining the log-tangent bundle $TX(-\log \underline{Z})$, which locally is generated by $\langle z_j \partial_{z_j}, \partial_{x_i} \rangle$ when $\underline{Z} = \bigcup_j Z_j$ with $Z_j = \{z_j = 0 \}$ (see \cite{GualtieriLiPelayoRatiu17} and \cite{Lanius17} respectively). In the former case, the log-tangent bundle $TX(-\log \underline{Z})$ is the fiber product of the individual log-tangent bundles $TX(-\log Z_j)$ associated to the log divisors $I_{Z_j}$ (c.f.\ below \autoref{prop:ziaiinvariant}).
\begin{exa}[Elliptic divisors]\label{exa:ellvx} Let $|D| = (R,q)$ be an elliptic divisor with elliptic divisor ideal $I_{|D|}$. Then the module $\Gamma(TX)_{I_{|D|}}$ is projective, defining the \emph{elliptic tangent bundle} $\mc{A}_{|D|} = TX(-\log |D|)$ \cite{CavalcantiGualtieri18}, which is locally given by $\Gamma(\mc{A}_{|D|}) = \Gamma(TX)_{I_{|D|}} = \langle r \partial_r, \partial_{\theta}, \partial_{x_i} \rangle$ in normal polar coordinates $(r,\theta)$ to $D$ in which $I_{|D|} = \langle r^2 \rangle$. From the fact that $r \partial_r \wedge \partial_{\theta} = r^2 \partial_x \wedge \partial_y$ where $r^2 = x^2 + y^2$ we see that elliptic divisors are standard (see \cite[Proposition 3.4.17]{Klaasse17}).
\end{exa}
\begin{rem} Morphisms of divisors often give rise to Lie algebroid morphisms between their primary ideal Lie algebroids. Indeed, in \cite{Klaasse17,CavalcantiKlaasse18} it is shown that a morphism between log divisor ideals leads to an induced Lie algebroid morphism between their resulting log-tangent bundles. Similarly for morphisms between elliptic divisor ideals. Morphisms between elliptic and log divisors are considered in \cite{CavalcantiKlaasse18} as boundary maps, and lead to Lie algebroid morphism from the elliptic tangent bundle to the log-tangent bundle. Finally, regarding morphisms between log and elliptic divisors one has to be more careful, see \cite{KlaasseLi18}.
\end{rem}
\begin{exa}[Complex log divisors]\label{exa:complexlogtgnt} Let $D = (U,\sigma)$ be a complex log divisor with ideal $I_\sigma$. Then $\Gamma(TX_\C)_{I_D}$ is a projective module of complex vector fields, defining the \emph{complex log-tangent bundle} $TX_{\C}(- \log D)$, locally given by $\langle w \partial_w, \partial_{\overline{w}}, \partial_{x_i} \rangle$ in coordinates where $I_\sigma = \langle w \rangle$ (see \cite{CavalcantiGualtieri18}). Because $w \partial_w \wedge \partial_{\overline{w}} = - 2i w \partial_x \wedge \partial_y$, we see that complex log divisors are standard. This complex Lie algebroid is also denoted by $\mc{A}_D \to X$, and we will not use it in this paper.
\end{exa}
We discuss one further example, obtained from an elliptic-log divisor. Let $W = Z \otimes |D|$ be an elliptic-log divisor with divisor ideal $I_W$. Then the primary ideal Lie algebroid $\mc{A}_W := \mc{A}_{I_W}$, the \emph{elliptic-log tangent bundle} $TX(-\log W)$, exists and is both an $\mc{A}_Z$ and an $\mc{A}_{|D|}$-Lie algebroid, as it is the fiber product of the two. The following was pointed out to us by Gil Cavalcanti.
\begin{prop}\label{prop:elllogalgebroid} Let $(X,Z)$ be a log pair and $|D| \subseteq Z$ an elliptic divisor. Then the module $\Gamma(TX)_{I_W}$ is projective, so that $\mc{A}_{W}$ is a Lie algebroid. Its divisor ideal is given by $I_{\mc{A}_W} = I_W$.
\end{prop}
\bp By definition $\Gamma(TX)_{I_W}$ consists of vector fields preserving $I_W = I_Z \cdot I_{|D|}$. It is immediate (see \cite{Klaasse17}) that all vector fields tangent to $D$ belong to $\Gamma(TX)_{I_W}$. Hence the associated sheaf is locally free if around points in $D$ it contains only two independent vector fields normal to $D$, as $D$ has codimension two. Away from $D$ we have $I_W = I_Z$, so that there $\Gamma(TX)_{I_W} = \Gamma(TX)_{I_Z}$ which we know to be locally free. Let $(x,y)$ be normal coordinates around $p \in D$ such that $I_W = \langle x (x^2 + y^2) \rangle$. Set $w(x,y) := x (x^2 + y^2)$ and take $V = a \partial_x + b \partial_y$. Note that $V$ belongs to $\Gamma(TX)_{I_W}$ if and only if the function $f = \mc{L}_V (\log w)$ is smooth. We compute that
\be
f = \mc{L}_V(\log w) = \mc{L}_V(\log x) + \mc{L}_V(\log(x^2 + y^2)) = a x^{-1} + (2a x + 2 b y)(x^2+y^2)^{-1},
\ee
so by multiplying by $x^2 + y^2$ and rearranging we obtain $(x^2 + y^2) f - 2a x - 2b y = a x^{-1} (x^2 + y^2)$. The left-hand side is smooth, hence so must be the right-hand side, i.e.\ $a$ must be divisible by $x$. Hence $a = x \alpha$ for some smooth $\alpha$. This implies that $f = \alpha + (2 x^2 \alpha + 2 b y)(x^2 + y^2)^{-1}$, so that $\lambda_1 := (x^2 \alpha + b y)(x^2 + y^2)^{-1}$ must be smooth. Rewriting this we obtain that
\be
	b = \lambda_1 y + (\lambda_1 - \alpha)x^2 y^{-1},
\ee
and as $b$ is smooth, $y$ divides $\lambda_1 - \alpha$. Set $\lambda_2 := (\lambda_1 - \alpha) y^{-1}$. Given smooth functions $\lambda_1$ and $\lambda_2$ we can now express $a$ and $b$ as $a = x (\lambda_1 - y \lambda_2)$ and $b = \lambda_1 y + \lambda_2 x^2$. This implies that $V$ is given by $V = \lambda_1 (x \partial_x + y \partial_y) + \lambda_2 (-x y \partial_x + x^2 \partial_y)$. This shows that $V$ must lie in the two-dimensional span $\langle x \partial_x + y \partial_y, x(y \partial_x - x \partial_y) \rangle$. To finish the proof of local freeness, we check that these generators preserve $I_W$. This follows, as we readily compute:
\bi
	\item \makebox[3cm][l]{$\mc{L}_{x \partial_x + y \partial_y}(\log w)$} $= x \partial_x \log x + (x \partial_x + y \partial_y) \log(x^2 + y^2) = 1 + 2 = 3$,
	\item \makebox[3cm][l]{$\mc{L}_{x(y \partial_x - x \partial_y)}(\log w)$} $= x y \partial_x \log x + x y \partial_x \log(x^2 + y^2) - x^2 \partial_y \log(x^2 + y^2)$
	\item[] \makebox[3cm][l]{} $= y + (2 x^2 y - 2 x^2 y)(x^2 + y^2)^{-1} = y$.
\ei
We conclude that $\mc{A}_W$ exists and that locally we have (with slight abuse of notation)
\be
\Gamma(\mc{A}_W) = \Gamma(TX)_{I_W} = \langle x \partial_x + y \partial_y, x(y \partial_x - x \partial_y) \rangle \oplus \Gamma(TD).
\ee
The final statement of this proposition follows from the realization that
\be
(x \partial_x + y \partial_y) \wedge x(y \partial_x - x \partial_y) = x (x^2 + y^2) \partial_x \wedge \partial_y = w \partial_x \wedge \partial_y. \qedhere
\ee
\ep
We can alternatively write the generators of $\mc{A}_W$ using a polar coordinate system $(r,\theta)$ normal to $D$ supplied by $I_{W} = \langle r^3 \cos \theta \rangle$. In such a coordinate system, we have that
\be
\Gamma(\mc{A}_W) = \langle r \partial_r, r \cos \theta \partial_\theta \rangle + \Gamma(TD).
\ee
In particular, \autoref{prop:elllogalgebroid} says that $I_W$ is standard. The elliptic-log tangent bundle is further explored in \autoref{exa:ellipticlogresidue} and \cite{KlaasseLanius18}, the latter computing its Lie algebroid cohomology.

A slightly different example is provided by the $b^k$-tangent bundles \cite{Scott16}. The idea here is that instead of considering derivations preserving a certain divisor ideal $I$ (which would give its primary ideal Lie algebroid), we consider those which send $I$ to $I^k$ in a certain sense.

Let $(X,Z)$ be a log pair with $\iota_Z\colon Z \hookrightarrow X$ the inclusion, $k \geq 1$ be a given integer. The \emph{sheaf of $k$-jets} at $Z$ is $\mc{J}_Z^k := \iota_Z^{-1}(C^\infty(X) / I_Z^{k+1})$. A \emph{$k$-jet} at $Z$ is a global section of $\mc{J}_Z^k$. Denote the set of $k$-jets at $Z$ by $J_Z^k$.
Given a function $f$ defined in a neighbourhood of $Z$, let $[f]^k_Z \in J_Z^k$ denote the $k$-jet at $Z$ represented by $f$. For a given $k$-jet at $Z$, $j \in J_Z^k$, write $f \in j$ if $[f]_Z^k = j$. 
\begin{exa}[$b^k$-tangent bundles, \cite{Scott16}]\label{exa:bktangentbundle} Let $j_{k-1} \in J_Z^{k-1}$ be a choice of $(k-1)$-jet and define the sheaf of vector fields $\Gamma(TX)_{I^k_Z, j_{k-1}} := \{v \in \Gamma(TX) \, | \, \mc{L}_v f \in I_Z^k \text{ for all } f \in j_{k-1}\}$. One verifies that this is a locally free sheaf of involutive submodules of $\Gamma(TX)$, and that given $v \in \mc{V}_X(I_Z)$ and $f \in C^\infty(X)$, the jet $[\mc{L}_v f]_Z^{k-1}$ depends only on $[f]_Z^{k-1}$. This defines the \emph{$b^k$-tangent bundle} $\mc{A}_Z^k \to X$ by $\Gamma(\mc{A}_Z^k) \cong \Gamma(TX)_{I^k_Z, j_{k-1}}$, which is locally given by $\Gamma(\mc{A}_Z^k) = \langle z^k \partial_z, \partial_{x_i} \rangle$ for $z \in j$.
\end{exa}
Note that when $k = 1$ the jet data is vacuous, so that $\mc{A}_Z^1 = \mc{A}_Z$, the log-tangent bundle. Further discussion of this process can be found in \cite{Scott16,Klaasse17}. The same construction described above can also be performed with smooth divisor ideals other than $I_Z$, where one uses a choice of $k$-jet at $Z_I$, i.e.\ a global section of $\mc{J}^k_I := \iota_{Z_I}^{-1}(C^\infty(X) / I^{k+1})$, for $\iota_{Z_I}\colon Z_I \hookrightarrow X$ the inclusion. The sheaf of vector fields then becomes $\Gamma(TX)_{I^k, j_{k-1}} := \{v \in \Gamma(TX) \, | \, \mc{L}_v f \in I^k \text{ for all } f \in j_{k-1}\}$.
\begin{rem} Let $I \subseteq C^\infty(X)$ be a projective smooth divisor ideal. Then its associated primary ideal Lie algebroid $TX_I$ can be restricted to $Z_I$, because it is the degeneracy locus of $TX_I$. Assuming that $Z_I$ is a transitive invariant submanifold for $TX_I$, the restriction $TX_I|_{Z_I}$ is then a Lie subalgebroid of ${\rm At}(NZ_I) \to Z_I$, the Atiyah algebroid of the normal bundle. The latter is a transitive Lie algebroid with rank equal to ${\rm rank}({\rm At}(NZ_I)) = \dim(Z_I) + {\rm codim}(Z_I)^2$. If $Z_I$ is not transitive, then instead there is only an almost-injective Lie algebroid morphism from $TX_I|_{Z_I}$ to ${\rm At}(NZ_I)$. As ${\rm rank}(TX_I|_{Z_I}) = {\rm rank}(TX_I) = \dim(X)$, this means that in order for $I$ to be projective, the vector fields normal to $Z_I$ which preserve $I$ must cut down the isotropy of ${\rm At}(NZ_I)$ in dimension from ${\rm codim}(Z_I)^2$ to ${\rm codim}(Z_I)$. Related to this is the discussion in \autoref{sec:divexamples} on Morse--Bott divisor ideals being locally homogeneous around $Z_I$.
	
	For example, if $I = I_Z$ is a log divisor ideal, then the above implies that $TX_{I_Z}|_Z \cong {\rm At}(NZ)$. For elliptic divisor ideals, instead $TX_{I_{|D|}}|_D$ is a corank-two Lie subalgebroid of ${\rm At}(ND)$.
\end{rem}
\subsection{Elementary modifications}
\label{sec:elementarymodifications}
In this section we discuss a process called \emph{elementary modification} \cite{GualtieriLi14, Li13} (partially also \emph{rescaling} \cite{Lanius16, Melrose93}), by which, given a Lie algebroid, new Lie algebroids are constructed by specifying a projective singular Lie subalgebroid using extra data. In loc.\ cit.\ this is done using Lie algebroids supported on hypersurfaces, but we consider a generalization of this procedure by using Lie algebroids supported on the support of a smooth divisor ideal. When using a log divisor ideal (see \autoref{exa:log}), one recovers the procedures of loc.\ cit. There are two procedures which are in a sense dual to another, referred to as lower and upper elementary modification respectively. As noted in \cite{Li13}, these processes are known in algebraic geometry under the name of Hecke modifications. Recall from \autoref{defn:smoothdivisor} that a divisor ideal $I \subseteq C^\infty(X)$ is smooth if its support $Z_I \subseteq X$ is a smooth submanifold.
\subsubsection{Upper elementary modification}
We first describe the process of upper modification (see \cite{Li13} for the case when $I = I_Z$), which uses surjective comorphisms. Let $I$ be a smooth divisor ideal with ideal sheaf $\mc{I}$, and let $\mc{A} \to X$ and $\mc{B} \to Z_I$ be vector bundles. Assume the existence of a surjective bundle comorphism $(\varphi;i_{Z_I})\colon (\mc{B},Z_I) \dashedrightarrow (\mc{A},X)$, i.e.\ the map $\varphi^*\colon i_{Z_I}^* \mc{A} \to \mc{B}$ is surjective, and has kernel $\mc{K} := \ker \varphi \to Z_I$. From this data, we define a $C^\infty(X)$-module
\be
\Gamma(\mc{A};(\mc{B},I)) := \{v \in \Gamma(\mc{A}) \otimes \mc{I}^{-1} \, | \, \mc{I} \cdot v \in \Gamma(\mc{A},\mc{K})\},
\ee
where it is immediate that $\Gamma(\mc{A}) \subseteq \Gamma(\mc{A};(\mc{B},I))$. Here $\mc{I}^{-1}$ is the inverse of $\mc{I}$ as a sheaf, and we make use of the canonical isomorphism $\mc{I} \otimes \mc{I}^{-1} \cong \mc{C}^\infty_X$. Alternatively, we can choose local generators for the sheaf $\mc{I}$ and its inverse and check that the definition is independent of such a choice, or first realize $\mc{I}$ by a divisor using \autoref{prop:locprincideal}. It is immediate that $\Gamma(\mc{A};(\mc{B},I))$ is locally finitely generated because $I$ and $\Gamma(\mc{A})$ are. The main point of this construction is:
\begin{lem}\label{lem:uppermodprojective} Under the above assumptions, the module $\Gamma(\mc{A};(\mc{B},I))$ is projective.
\end{lem}
\bp Trivialize $\mc{B}^m \to Z_I$ around a point in $Z_I$ with basis of sections $(w_1,\dots,w_m)$, and let $f \in \mc{I}$ be a local generator. Choose a complementary basis of sections $(v_{m+1},\dots,v_n)$ for $\mc{K}$, so that $(w_1,\dots,w_m,v_{m+1},\dots,v_n)$ is a local basis of $\Gamma(\mc{A}^n)$. Then $(w_1,\dots,w_m, f^{-1} v_{m+1},\dots, f^{-1} v_n)$ is a local basis of $\Gamma(\mc{A};(\mc{B},I))$, showing projectivity, as $\Gamma(\mc{A};(\mc{B},I)) \cong \Gamma(\mc{A})$ away from $Z_I$.
\ep
By the Serre--Swan correspondence there exists a vector bundle
\be
	\{\mc{A}\:(\mc{B},I)\} \to X \quad \text{defined by} \quad \Gamma(\{\mc{A}\:(\mc{B},I)\}) \cong \Gamma(\mc{A};(\mc{B},I)).
\ee
Next, assume that $\mc{A} \to X$ and $\mc{B} \to Z_I$ are anchored and $\varphi$ is an anchored comorphism, i.e.\ $i_{Z_I}^* (\mc{L}_{\rho_{A}(v)} f) = \mc{L}_{\rho_\mc{B}(\varphi(v))}(i_{Z_I}^*(f))$ for all $v \in \Gamma(\mc{A})$ and $f \in C^\infty(X)$. Alternatively, the anchor maps should satisfy $\rho_{\mc{B}}(\varphi^*(v)) \sim_\varphi \rho_{\mc{A}}(v)$ for all $v \in \Gamma(\mc{A})$.
\begin{lem}[c.f.\ \cite{Li13}]\label{lem:uppermodanchored} Under the above assumptions, the bundle $\{\mc{A}\:(\mc{B},I)\}$ is anchored.
\end{lem}
\bp The anchored comorphism condition implies $\mc{K} \subseteq \ker(\rho_{A}|_{Z_I})$, so that for $v \in \Gamma(\mc{A},K)$ we have $\rho_{A}(v) \in \Gamma(TX,TZ_I)$, which shows that $\rho_{A}$ lifts to an anchor for $\{\mc{A}\:(\mc{B},I)\}$.
\ep
There is an induced (anchored) morphism $\mc{A} \to \{\mc{A}\:(\mc{B},I)\}$ with isomorphism locus $Z_I$. Finally, assume that $\mc{A} \to X$ and $\mc{B} \to Z_I$ are Lie algebroids and $\varphi$ is a surjective Lie algebroid comorphism, i.e.\ $\varphi$ further satisfies $\varphi^*([v,w]) = [\varphi^*(v),\varphi^*(w)]$ for all $v,w \in \Gamma(\mc{A})$.
\begin{lem}\label{lem:uppermodinvolutive} Under the above assumptions, the module $\Gamma(\mc{A};(\mc{B},I))$ has an induced bracket.
\end{lem}
\bp Because $\varphi^*$ preserves brackets, there is an induced Lie bracket on $\Gamma(\mc{A};(\mc{B},I))$ inherited from the module $\Gamma(\mc{A})$ through the natural map $\Gamma(\mc{A}) \to \Gamma(\mc{A};(\mc{B},I))$.
\ep
\autoref{lem:uppermodprojective}, \autoref{lem:uppermodanchored} and \autoref{lem:uppermodinvolutive} combine to the following definition (c.f.\ \cite{Li13}).
\begin{defn}\label{defn:uppermod} Let $\mc{A} \to X$ be a Lie algebroid, $I \subseteq C^\infty(X)$ a smooth divisor ideal and $\mc{B} \to Z_I$ a Lie algebroid with a surjective Lie algebroid comorphism onto $\mc{A}$. The \emph{$(\mc{B},I)$-upper modification} of $\mc{A}$ is the 
	%$\mc{A}$-coLie algebroid  $\mc{A} \to
	Lie algebroid $\{\mc{A}\:(\mc{B},I)\}$ for which $\Gamma(\{\mc{A}\:(\mc{B},I)\}) \cong \Gamma(\mc{A};(\mc{B},I))$.
\end{defn}
In this case there is a Lie algebroid morphism $\mc{A} \to \{\mc{A}\:(\mc{B},I)\}$ making $\mc{A}$ into an $\{\mc{A}\:(\mc{B},I)\}$-Lie algebroid, and the $\{\mc{A}\:(\mc{B},I)\}$-divisor ideal of $\mc{A}$ is given by $I^k$, where $k = \dim(\mc{B})$.
\begin{rem} From the divisor ideal of the induced morphism $\mc{A} \to \{\mc{A}\:(\mc{B},I)\}$ we see that it is only possible for $\{\mc{A}\:(\mc{B},I)\}$ to be anchored if the divisor of $\mc{A}$ is divisible by $I^k$. This places an often nontrivial restriction on when upper modifications can be Lie algebroids.
\end{rem}
There is a commutativity property for upper modifications, which is readily verified.
\begin{prop} Let $\mc{A} \to X$ be a Lie algebroid, $I,I' \subseteq C^\infty(X)$ smooth divisor ideals and $(\mc{B},Z_I)$, $(\mc{B}',Z_{I'})$ Lie algebroids equipped with surjective Lie algebroid comorphisms onto $\mc{A}$. Then upper modification at $Z_I$ and $Z_{I'}$ commute, i.e. we have isomorphisms
	\be
	\{\{\mc{A}\:(\mc{B},I)\}\:(\mc{B}',I')\} \cong \{\{\mc{A}\:(\mc{B}',I')\}\:(\mc{B},I)\}.
	\ee
\end{prop}
\subsubsection{Lower modification}
We turn to the dual procedure to upper modification, namely lower modification (or rescaling). This uses injective morphisms instead of surjective comorphisms.
 
Let $I$ be a smooth divisor ideal and $\mc{B} \to Z_I$ be a vector subbundle of a vector bundle $\mc{A} \to X$. To define an appropriate $C^\infty(X)$-module in this setting, we remark that injective morphisms dualize to surjective comorphisms. That is, the bundle inclusion $(\varphi,i_{Z_I})\colon (\mc{B},Z_I) \to (\mc{A},X)$ dualizes to the surjection $(\varphi;i_{Z_I})\colon (\mc{B}^*,Z_I) \dashrightarrow (\mc{A}^*,X)$. We can thus perform upper modification to obtain the vector bundle
\be
	\{\mc{A}^*\:(\mc{B}^*,I)\} \to X \quad \text{with} \quad \Gamma(\{\mc{A}^*\:(\mc{B}^*,I)\}) \cong \Gamma(\mc{A}^*\:(\mc{B}^*,I)).
\ee
This comes equipped with a natural bundle morphism $\mc{A}^* \to \{\mc{A}^*\:(\mc{B}^*,I)\}$. Dualizing again, using $(\mc{A}^*)^* \cong \mc{A}$ we obtain a comorphism $\mc{A} \dashrightarrow \{\mc{A}^*\:(\mc{B}^*,I)\}^*$. Reinterpreted, this says that there is a bundle morphism $\{\mc{A}^*\:(\mc{B}^*,I)\}^* \to \mc{A}$. This motivates the following definitions:
\be
	[\mc{A}\:(\mc{B},I)] := \{\mc{A}^*\:(\mc{B}^*,I)\}^* \quad \text{and} \quad \Gamma(\mc{A},(\mc{B},I)) := \Gamma([\mc{A}\:(\mc{B},I)]).
\ee
\begin{rem} We had to proceed by dualizing twice because of the use of a smooth divisor ideal $I$. When $I = I_Z$, one can merely define $\Gamma(\mc{A},(\mc{B},I_Z)) := \Gamma(\mc{A},\mc{B})$ as in \cite{Melrose93,Lanius16,GualtieriLi14}. Note here that the `dual definition' $\{v \in \Gamma(\mc{A}) \otimes \mc{I} \, | \, \mc{I}^{-1} \cdot v \in \Gamma(\mc{A},\mc{B})\}$ does not make sense.
\end{rem}
Now that we have the module of interest, it can be described locally, as in \autoref{lem:uppermodprojective}. Trivialize $\mc{B}^m \to Z_I$ around a point in $Z_I$ with basis of sections $(w_1, \dots w_m)$, and let $f \in \mc{I}$ be a local generator. Extend this basis to a local basis $(w_1,\dots,w_m,v_{m+1},\dots,v_n)$ of $\Gamma(\mc{A}^n)$. Then $(w_1,\dots,w_m, f v_{m+1}, \dots, f v_n)$ is a local basis for $\Gamma(\mc{A},(\mc{B},I))$. This also shows projectivity.

When $\mc{A} \to X$ is an anchored bundle and $\mc{B} \to Z_I$ an anchored subbundle of $\mc{A}$, the bundle $[\mc{A}\:(\mc{B},I)]$ becomes anchored by inheriting an anchor from $\rho_{A}\colon \mc{A} \to TX$. Note that there is an induced (anchored) morphism $[\mc{A}\:(\mc{B},I)] \to \mc{A}$ with isomorphism locus $Z_I$. Note that $[\mc{A}\:(\mc{B},I)]$ is \emph{not} an anchored subbundle of $\mc{A}$ (except in the trivial case when $I = C^\infty(X)$).

Next assume that $\mc{A} \to X$ is a Lie algebroid and $\mc{B} \to Z_I$ is a Lie subalgebroid.
\begin{lem}\label{lem:lowermodinvolutive} Under the above assumptions, the submodule $\Gamma(\mc{A},(\mc{B},I)) \subseteq \Gamma(\mc{A})$ is involutive.
\end{lem}
\bp The submodule $\Gamma(\mc{A},\mc{B}) \subseteq \Gamma(\mc{A})$ is involutive by \autoref{prop:liesubalgeroid}. To extend this property to $\Gamma(\mc{A},(\mc{B},I))$ we can make use of the local description above. Let $f \in \mc{I}$ be a local generator and consider two sections $f s, f s' \in \Gamma(\mc{A},(\mc{B},I))$ with $s \in \Gamma(\mc{A})$ but not in $\Gamma(\mc{A},\mc{B})$ (here $g \in \mc{I}$ can be written as $g = g' f$, but then $g'$ can be absorbed into $s$). Then we have:
\be
	[f s, f s']_\mc{A} = f^2 [s,s']_\mc{A} + f \left((\mc{L}_{\rho_{A}(s)} f) s' - (\mc{L}_{\rho_{A}(s')} f) s\right).
\ee
We see from this that $[fs,f s']_\mc{A} \in \mc{I} \otimes \Gamma(\mc{A})$, as in \autoref{lem:liesubalgs}. The other cases are similar.
\ep
Thus, making use of \autoref{lem:lowermodinvolutive} we see that $\Gamma(\mc{A},(\mc{B},I))$ is a projective singular Lie subalgebroid of $\mc{A}$ in the sense of \autoref{defn:singularliealgebroid}. This leads to the following.
\begin{defn}[{\cite{GualtieriLi14, Lanius16, Melrose93}}]\label{defn:lowermod} Let $\mc{A} \to X$ be a Lie algebroid, $I \subseteq C^\infty(X)$ a smooth divisor ideal and $\mc{B} \to Z_I$ a Lie subalgebroid of $\mc{A}$. Then the \emph{$(\mc{B},I)$-lower modification} of $\mc{A}$ is the $\mc{A}$-Lie algebroid $[\mc{A}\:(\mc{B},I)] \to \mc{A}$ for which $\Gamma([\mc{A}\:(\mc{B},I)]) \cong \Gamma(\mc{A},(\mc{B},I))$.
\end{defn}
It is immediate that the $\mc{A}$-divisor ideal of $[\mc{A}\:(\mc{B},I)]$ is given by $I^k$, where $k = {\rm codim}(\mc{B})$.
\begin{rem} When $I = I_Z$, the above procedure is also called \emph{rescaling} \cite{Lanius16, Melrose93}; we then denote $[\mc{A}\:(\mc{B},I_Z)]$ by $[\mc{A}\:\mc{B}]$ or ${}^{\mc{B}} \mc{A}$ and refer to it as the \emph{$(\mc{B},Z)$-rescaling} of $\mc{A}$.
\end{rem}
\begin{rem}\label{rem:bbprimerescaling} Let $\mc{A} \to X$ be a Lie algebroid, $I$ a smooth divisor ideal and $\mc{B}, \mc{B}' \to Z_I$ two Lie subalgebroids of $\mc{A}$ such that $\mc{B} \subseteq \mc{B}'$. Then $[\mc{A}\:(\mc{B},I)]$ is a $[\mc{A}\:(\mc{B}',I)]$-Lie algebroid.
\end{rem}
One can always modify using $0_{Z_I} \subseteq \mc{A}|_{Z_I}$, the trivial Lie subalgebroid, with
\be
	[\mc{A}\:(0_{Z_I},I)] = I \cdot \mc{A},
\ee
which is the secondary ideal Lie algebroid with section module $I \cdot \Gamma(\mc{A})$. Thus lower modification by $0_{Z_I}$ is the harshest, while lower modification using $\mc{B} = \rho_\mc{A}^{-1}(T(Z_I))$ (if it exists, c.f.\ \autoref{rem:inverseimagesmooth}) is the mildest. Note that $0_{Z_I}$-modification is not idempotent.
The process of lower modification is commutative when performed relative to disjoint submanifolds.
\begin{prop}\label{prop:lowermodifycommute} Let $\mc{A} \to X$ be a Lie algebroid, $I,I' \subseteq C^\infty(X)$ smooth divisor ideals and $(\mc{B}, Z_I)$, $(\mc{B}',Z_{I'})$ Lie subalgebroids of $\mc{A}$ such that $Z_I \cap Z_{I'} = \emptyset$. Then lower modification at $Z_I$ and $Z_{I'}$ commute and coincides with taking the fiber product, i.e. we have isomorphisms
\be
	[[\mc{A}\:(\mc{B},I)]:(\mc{B}',I')] \cong [\mc{A}\:(\mc{B},I)] \oplus_{\mc{A}} [\mc{A}\:(\mc{B}',I')] \cong [[\mc{A}\:(\mc{B}',I')]:(\mc{B},I)].
\ee
\end{prop}
The proof of the above is immediate upon noticing we can view $(\mc{B}',Z_{I'})$ as a Lie subalgebroid of $[\mc{A}\:(\mc{B},I)]$ because the latter is isomorphic to $\mc{A}$ outside $Z_I$, and $Z_I \cap Z_{I'} = \emptyset$.
\begin{rem} Elementary modification can also be performed in an even slightly more general setting. Namely, instead of a smooth divisor ideal one can use a $C^\infty(X)$-module $\mc{M}$ that is locally free of rank one, and a suitable vector bundle supported on its smooth support.
\end{rem}
Upper modification is a left inverse to lower modification, and vice versa regarding a right inverse when this makes sense. Often the morphism used is the anchor map on sections.
\begin{prop} Let $\mc{A} \to X$ be a Lie algebroid, $I \subseteq C^\infty(X)$ a smooth divisor ideal and $\mc{B} \to Z_I$ a Lie subalgebroid of $\mc{A}$. Then $\{[\mc{A}\:(\mc{B},I)]\:(\mc{B},I)\} \cong \mc{A}$, and similarly when there is a surjective Lie algebroid comorphism from $\mc{B}$ onto $\mc{A}$ we have $[\{\mc{A}\:(\mc{B},I)\}\:(\mc{B},I)] \cong \mc{A}$.
\end{prop}
\subsubsection{Examples}
\label{sec:examplesmodification}
We finish the discussion of modifications by considering some examples. See also \cite{Melrose95} for further examples when $I = I_Z$, which are of interest in geometric analysis.
\begin{exa} Let $(X,Z)$ be a log pair and consider $\mc{A} = TX$ and $\mc{B} = 0_Z$. Then $[\mc{A}\:\mc{B}] = \mc{B}_Z$, the zero tangent bundle associated to $Z$ (see \autoref{exa:zvx}).
\end{exa}
\begin{exa} Let $(X,Z)$ be a log pair and consider $\mc{A} = TX$ and $\mc{B} = TZ$. Then $[\mc{A}\:\mc{B}] = \mc{A}_Z$, the log-tangent bundle associated to $Z$ (see \autoref{exa:zvx}).
\end{exa}
The above example can be extended to cases where $\mc{B} = TF$ for $F$ an involutive distribution on $Z$, for example the kernel of a closed one-form on $Z$. The resulting Lie algebroid $TX(-\log(Z,F)) := [TX\:TF]$ occurs in the context of log-Poisson structures (c.f.\ \cite{KlaasseLanius18,GualtieriLiPelayoRatiu17}).
\begin{exa} Let $(X,Z)$ be a log pair and consider $\mc{A} = \mc{A}_Z$ and $\mc{B} = TZ$ with surjective bundle morphism $\rho_\mc{A}|_Z\colon \mc{A}|_Z \to \mc{B}$. Then $\rho_\mc{A}$ is a Lie algebroid comorphism and $\{\mc{A}_Z\:\mc{B}\} = TX$.
\end{exa}
\begin{exa} Let $(X,Z)$ be a log pair and consider $\mc{A} = \mc{A}_Z$ and $\mc{B} = 0_Z$. Then $[\mc{A}\:\mc{B}] = \mc{C}_Z$, the \emph{scattering tangent bundle} of \cite{Lanius16,Melrose95}.
\end{exa}
Further examples obtained using elliptic divisor ideals can be found in \cite{KlaasseLanius18}.
%
%END OF TEX FILE

%% file: poissonstructuresdivisors.tex
\section{Poisson structures of divisor-type}
\label{sec:poissondivtype}
In this section we discuss the notion of an $\mc{A}$-Poisson structure (c.f.\ \cite{CannasDaSilvaWeinstein99}), their interaction with divisors, and the process of lifting using base-preserving Lie algebroid morphisms. Most $\mc{A}$-Poisson structures we will encounter will be generically nondegenerate, or generically regular in a precise sense. At first reading, the reader can let $\mc{A} = TX$ in most of what follows.
\begin{defn} An \emph{$\mathcal{A}$-Poisson structure} is an $\mathcal{A}$-bivector $\pi_\mathcal{A}\in \Gamma(\wedge^2\mathcal{A})$ with $[\pi_\mathcal{A},\pi_\mathcal{A}]_{\mc{A}} = 0$.
\end{defn}
We denote the space of such $\mc{A}$-bivectors by ${\rm Poiss}(\mc{A})$, and set ${\rm Poiss}(X) := {\rm Poiss}(TX)$. Of particular interest to us are those Poisson bivectors that are \emph{generically nondegenerate}, i.e.\ for which $\pi^\sharp_{\mc{A}}\colon \mc{A}^* \to \mc{A}$ is nondegenerate on an open dense subset of $X$. As for Poisson structures, an $\mc{A}$-Poisson structure induces an $\mc{A}$-Lie algebroid structure on the dual bundle $\mc{A}^* \to X$.
\begin{defn}\label{defn:apoissonalgebroid} Let $\pi_\mc{A} \in {\rm Poiss}(\mc{A})$ be given. The \emph{$\mc{A}$-Poisson algebroid} $\mc{A}^*_{\pi_\mc{A}}$ consists of the bundle $\mc{A}^*$ with $\mc{A}$-anchor map $\pi_\mc{A}^\sharp\colon \mc{A}^* \to \mc{A}$, whose Lie bracket on $\Gamma(\mc{A}^*_{\pi_{\mc{A}}})$ is given by
	\be
	[v,w]_{\mc{A}^*_{\pi_\mc{A}}} = \mc{L}_{\pi_\mc{A}^\sharp v} w - \mc{L}_{\pi_\mc{A}^\sharp w} v - d_{\mc{A}} \pi_{\mc{A}}(v,w), \qquad v,w \in \Gamma(\mc{A}^*_{\pi_{\mc{A}}}).
	\ee
\end{defn}
Poisson structures can be related to the divisors of \autoref{sec:divonmanifolds} as follows. Denote the space of $\mc{A}$-bivectors by $\mf{X}^2(\mc{A}) = \Gamma(\wedge^2 \mc{A})$ (\emph{not} $\Gamma(\wedge^2 T\mc{A})$), so that we have $\mf{X}^2(X) = \mf{X}^2(TX)$, although no confusion should arise. The \emph{Pfaffian} of $\pi_\mc{A}\in \mf{X}^2(\mc{A}^{2n})$ is the section $\wedge^n \pi_\mc{A} \in \Gamma(\det(\mc{A}))$. An $\mc{A}$-Poisson structure $\pi_{\mc{A}}$ is called \emph{nondegenerate} if $\pi_{\mc{A}}^\sharp\colon \mc{A}^* \to \mc{A}$ is an isomorphism, whose existence forces the rank of $\mc{A}$ to be even. Alternatively, one can demand the Pfaffian $\wedge^n \pi_\mc{A}$ to be nowhere vanishing, and similarly for $\wedge^m \pi_{\mathcal{A}}$ regarding $2m$-regularity.
\begin{defn}\label{defn:adivisortype} An $\mc{A}$-Poisson structure $\pi_\mc{A} \in {\rm Poiss}(\mc{A})$ is of \emph{$m$-divisor-type} for some $m \geq 0$ if $\wedge^{m+1} \pi_{\mc{A}} \equiv 0$ and there exists a line subbundle $K \subseteq \wedge^{2m} \mc{A}$ such that $(K, \wedge^{m} \pi_\mc{A})$ is a divisor. Denote the divisor of $\pi_\mc{A}$ by ${\rm div}(\pi_\mc{A})$ and the associated divisor ideal by $I_{\pi_\mc{A}} \subseteq C^\infty(X)$.
\end{defn}
Note that if $\pi_\mc{A}$ is of $m$-divisor-type, it generically has rank $2m$, and if $\pi_{\mathcal{A}}$ specifies the trivial divisor, then $\pi_{\mathcal{A}}$ is $2m$-regular. In the above definition, when $2m = {\rm rank}(\mc{A})$, we will simply say that $\pi_\mc{A}$ is of \emph{divisor-type}, and then the line bundle $K = \det(\mc{A})$ can be used.
\begin{rem}\label{rem:apoissontriangbialgebroid} An $\mc{A}$-Poisson structure defines a triangular Lie bialgebroid $(\mc{A}, \mc{A}^*_{\pi_\mc{A}})$ \cite{KosmannSchwarzbach95,KosmannSchwarzbachLaurentGengoux05,LiuWeinsteinXu97,MackenzieXu94}: in general, a pair of Lie algebroids $(\mc{A},\mc{A}^*)$ is a Lie bialgebroid if for all $v, w \in \Gamma(\mc{A})$
	\be
	d_{\mc{A}^*} [v,w]_\mc{A} = \mc{L}_{v} d_{\mc{A}^*} w - \mc{L}_w d_{\mc{A}^*} v.
	\ee
	Any Lie bialgebroid $(\mc{A},\mc{A}^*)$ induces a Poisson structure $\pi \in {\rm Poiss}(X)$ with $\pi^\sharp = \rho_{\mc{A}} \circ \rho_{\mc{A}^*}^*$. Not every Lie bialgebroid is triangular (yet $\pi = \rho_{\mc{A}}(\pi_{\mc{A}})$ if it is), and the maps $\pi_{\mc{A}}^\sharp \colon \mc{A}^* \to \mc{A}$ and $\rho_{\mc{A}^*}^* \colon T^*X \to \mc{A}^*$ are Lie algebroid morphisms. These types of structures are also considered in \cite{Lanius16,Lanius18}, for Lie algebroids whose isomorphism locus is the complement of a hypersurface.
\end{rem}
\begin{rem} Given $\pi_{\mathcal{A}} \in {\rm Poiss}(\mc{A})$ we obtain a (singular) $\mc{A}$-distribution $D_{\pi_\mc{A}} = \pi_\mc{A}^\sharp(\mc{A}^*)$. Regularity of $\pi_\mc{A}$ is the same as regularity of $D_{\pi_\mc{A}}$, and then $\det(D_{\pi_\mc{A}}) = K$, with $K \subseteq \wedge^{2m} \mc{A}$.
\end{rem}
\begin{rem} If $\pi_\mc{A}$ is of $m$-divisor-type, it is generically regular of rank $2m$. This is not an if and only if unless $2m = \dim(X)$: being of $m$-divisor-type for $2m < \dim(X)$ further demands that the way in which $\pi_{\mathcal{A}}$ degenerates is controlled in a precise sense. See also \autoref{sec:almostregularpoisson}.
\end{rem}
Given a divisor ideal $I \subseteq C^\infty(X)$, we say $\pi_\mc{A}$ is of \emph{($m$-)$I$-divisor-type} if $I_{\pi_\mc{A}} = I$, noting that we must specify both the divisor ideal $I$, and the generic rank $2m$. Given a divisor ideal $I$, we denote by $\mf{X}^2_I(\mc{A}) \subseteq \mf{X}^2(\mc{A})$ the space of divisor $\mc{A}$-bivectors $\pi_\mc{A}$ such that $I_{\pi_\mc{A}} = I$. We denote by ${\rm Poiss}_I(\mc{A})$ the space of $\mc{A}$-Poisson structures of $I$-divisor-type, and by ${\rm Poiss}_{I,m}(\mc{A})$ the space of $\mc{A}$-Poisson structures of $m$-$I$-divisor-type.
\begin{rem} One of our main interests lies in nondegenerate $\mc{A}$-Poisson structures, as we then have access to symplectic techniques (see \autoref{sec:symplecticliealgebroids}). Indeed, the goal of the lifting process described in Section \ref{sec:liftingpoisson} is to find a Lie algebroid $\mc{A} \to X$ for which a given $\pi \in {\rm Poiss}(X)$ can be viewed as a nondegenerate $\mc{A}$-Poisson structure.
\end{rem}
A function $f \in C^\infty(X)$ is said to be \emph{$\pi_{\mathcal{A}}$-Casimir} if $\pi_\mc{A}^\sharp(d_\mc{A} f) = 0$ (alternatively, $[f,\pi_{\mathcal{A}}]_\mc{A} = 0$ using the $\mc{A}$-Schouten bracket). Given an $\mc{A}$-Poisson structure $\pi_\mc{A}$ of $m$-divisor-type and a $\pi_\mc{A}$-Casimir function $f$, we can form a new $\mc{A}$-Poisson structure $\pi'_\mc{A} = f \cdot \pi_\mc{A}$. If the zero set of $f$ is nowhere dense, it specifies a divisor $(\underline{\R},f)$, so as $\wedge^{m} \pi'_\mc{A} = f^m \cdot \wedge^m \pi_\mc{A}$ we conclude that $\pi'_\mc{A}$ is again of $m$-divisor-type. Moreover, $I_{\pi'_\mc{A}} =  \langle f^m \rangle \cdot I_{\pi_\mc{A}}$. See also \cite[Lemma 3.3]{AndroulidakisZambon17} and \autoref{sec:almostregularpoisson}.
\begin{rem} When $\pi_\mc{A} \in {\rm Poiss}(\mc{A})$ has maximal rank equal to $2$ it follows that $f \cdot \pi_\mc{A}$ is $\mc{A}$-Poisson for \emph{any} function $f \in C^\infty(X)$ (as is known for $\mc{A} = TX$): we have $[f\cdot \pi_\mc{A}, f \cdot \pi_\mc{A}]_\mc{A} = f^2 [\pi_\mc{A}, \pi_\mc{A}]_\mc{A} + 2 f \pi_\mc{A} \wedge \pi_\mc{A}^\sharp(d_\mc{A} f)$, where the latter term vanishes as $\pi_\mc{A} \wedge v = 0$ for all $v \in {\rm im}(\pi^\sharp_\mc{A})$. This can be used to construct examples of $\mc{A}$-Poisson structures of higher divisor-type.
\end{rem}
Given an $\mc{A}$-Poisson structure, we obtain its \emph{$\mc{A}$-Hamiltonian vector fields} as the image of the map $f \mapsto \pi_\mc{A}^\sharp(d_\mc{A} f) =: v_{\mc{A}, f} \in \mf{X}(\mc{A})$. Similarly we have the \emph{$\mc{A}$-Poisson vector fields}, i.e.\ those $v_\mc{A} \in \mf{X}(\mc{A})$ such that $\mc{L}_{v_\mc{A}} \pi_{\mc{A}} = 0$, using the $\mc{A}$-Lie derivative. Define the spaces ${\rm Ham}(\pi_\mc{A}) \subseteq {\rm Poiss}(\pi_\mc{A}) \subseteq \mf{X}(\mc{A})$ as for Poisson structures. Given $f \in C^\infty(X)$ we have $\pi_{\mc{A}}^\sharp(d_{\mc{A}} f) = \pi_{\mc{A}}^\sharp(\rho_{\mc{A}}^*(df))$, so that by setting $\pi := \rho_\mc{A}(\pi_\mc{A})$, we see from the relation $\pi^\sharp = \rho_{\mc{A}} \circ \pi_{\mc{A}}^\sharp \circ \rho_{\mc{A}}^*$ that $\rho_\mc{A}(v_{\mc{A},f}) = V_f$. Thus the $\mc{A}$-Hamiltonian vector fields for $\pi_\mc{A}$ surject onto the Hamiltonian vector fields for $\pi$.
\begin{rem}\label{rem:productapoisson} Let $(X,\mc{A},\pi_\mc{A})$ and $(X',\mc{A}',\pi_{\mc{A}'})$ be two $\mc{A}$-Poisson manifolds. Then their product $(X \times X', \mc{A} \oplus \mc{A}', \pi_{\mc{A}} + \pi_{\mc{A}'})$ is also $\mc{A}$-Poisson, where $\mc{A} \oplus \mc{A}' \to X \times X'$ is the direct product of Lie algebroids (see \cite{HigginsMackenzie90}), and both projections are $\mc{A}$-Poisson maps. Thus the class of $\mc{A}$-Poisson structures of divisor-type is closed under products, using the external tensor product of divisors (see \autoref{defn:divdirectstum}). This is noted in \cite{AndroulidakisZambon17} for $\mc{A} = TX$ using \autoref{cor:divtypealmostreg}.
\end{rem}
\subsection{Degeneraci loci}
We discuss degeneraci loci of $\mc{A}$-Poisson structures $\pi_\mc{A}$, as in \cite{Pym13}. Given $\pi_{\mc{A}} \in \mf{X}^2(\mc{A}^n)$ with bundle map $\pi_\mc{A}^\sharp\colon \mc{A}^* \to \mc{A}$, we have both $\det(\pi_\mc{A}^\sharp) \in \Gamma(\det(\mc{A})^2)$ and $\wedge^n \pi_\mc{A} \in \Gamma(\det(\mc{A}))$, which are related by $\det(\pi_\mc{A}^\sharp) = \wedge^n \pi_\mc{A} \otimes \wedge^n \pi_\mc{A}$. Considering the definition for Lie algebroids, we could define the degeneracy loci of $\pi_\mc{A}$ using the map $\pi_\mc{A}^\sharp$, i.e.\ as the degeneracy locus of its $\mc{A}$-Poisson algebroid $\mc{A}^*_{\pi_{\mathcal{A}}}$. However, $\pi^\sharp_\mc{A}$ is skew-symmetric, so its vanishing degree is too high. Instead we will use the Pfaffians $\wedge^k\pi_\mc{A} \in \Gamma(\wedge^{2k} \mc{A})$ for $k \geq 0$. These give rise to maps $\wedge^k \pi_\mc{A}\colon \Gamma(\wedge^{2k} \mc{A}^*) \to C^\infty(X)$ on sections, and thus define ideals $I_{\pi_\mc{A},2k}$.
\begin{defn} Let $\pi_\mc{A} \in {\rm Poiss}(\mc{A})$. The \emph{$2k$th degeneracy locus} of $\pi_\mc{A}$ for $k \geq 0$ is the subspace $X_{\pi_\mc{A},2k} \subseteq X$ of points where ${\rm rank}(\pi_{\mathcal{A}}) \leq 2k$, defined by the ideal $I_{\pi_\mc{A},2k+2} \subseteq C^\infty(X)$.
\end{defn}
Thus the degeneraci loci of $\mc{A}_{\pi_\mc{A}}^*$ and $\pi_\mc{A}$ agree as subspaces of $X$, but their ideals do not. Recall that a submanifold $N$ is \emph{$\pi$-Poisson} if $\pi^\sharp(T_x^* X) \subseteq T_x N$ for all $x \in N$. When $N$ is closed this is equivalent to its vanishing ideal $I_N$ being a \emph{Poisson ideal}, i.e.\ $\{I_N,C^\infty(X)\}_\pi \subseteq I_N$, or to $N$ being $T^*_\pi X$-invariant. Moreover, note that $I$ being a $\pi$-Poisson ideal is equivalent to $\pi^\sharp(d I) = 0$. When $N$ is $\pi$-Poisson, it carries a uniquely induced Poisson structure $\pi_N$ for which the inclusion is a Poisson map. The following is then a consequence of \autoref{prop:degenlocus}.
\begin{prop}\label{prop:pidegenpipoisson} Each degeneracy locus $X_{\pi,2k}$ is a $\pi$-Poisson submanifold if it is smooth.
\end{prop}
For $\mc{A}$-Poisson structures the notion of a \emph{$\pi_{\mc{A}}$-Poisson submanifold} is more involved: it instead is a Lie subalgebroid $(\varphi,i)\colon (\mc{B},N) \hookrightarrow (\mc{A},X)$ with $\pi_\mc{B} \in {\rm Poiss}(\mc{B})$ with $\varphi_*(\pi_\mc{B}) = \pi_\mc{A}$.
A point $x \in X$ is $\pi_{\mathcal{A}}$-\emph{regular} if the rank of $\pi_\mc{A}$ is constant in some open neighbourhood of $x$, and singular otherwise. We denote by $X_{\pi_\mc{A}, {\rm reg}} \subseteq X$ the space of $\pi_\mc{A}$-regular points, also called the \emph{regular locus} of $\pi_\mc{A}$. This is an open dense subspace of $X$, with its complement $X_{\pi_\mc{A}, {\rm sing}}$, the \emph{singular locus} of $\pi_\mc{A}$, being closed and nowhere dense. Denote by $X_{\pi_\mc{A}, {\rm max}} \subseteq X$ the open subspace where $\pi_\mc{A}$ has maximal rank. It immediately follows that $X \backslash X_{\pi_\mc{A}, {\rm max}} = X_{\pi_\mc{A}, 2k-2} = (\wedge^k \pi_\mc{A})^{-1}(0)$, where $2k$ is the maximal rank of $\pi_\mc{A}$. In general we have $X_{\pi_\mc{A}, {\rm max}} \subsetneq X_{\pi_\mc{A}, {\rm reg}}$.
\subsection{Relation to almost-regularity}
\label{sec:almostregularpoisson}
A Poisson structure $\pi \in {\rm Poiss}(X)$ specifies a singular Lie subalgebroid $\mc{F}_\pi := \pi^\sharp(\Gamma(T_\pi^*X))$ of $TX$ consisting of $\pi$-Hamiltonian vector fields. Following Androulidakis--Zambon \cite{AndroulidakisZambon17}, we call the Poisson structure $\pi$ \emph{almost-regular} if $\mc{F}_\pi$ is projective. There is a nice criterion for almost-regularity in terms of the $\pi$-symplectic leaves.
\begin{prop}[{\cite[Theorem 2.8]{AndroulidakisZambon17}}]\label{prop:almostregulardistr} Let $\pi \in {\rm Poiss}(X)$. Then $\pi$ is almost-regular if and only if $X_{\pi, {\rm max}} \subseteq X$ is dense and there exists a distribution $D_\pi \subseteq TX$ such that $D_{\pi,x} = T_x L$ for all $x \in X_{\pi, {\rm max}}$, where $L$ is the $\pi$-symplectic leaf through $x$.
\end{prop}
Hence almost-regularity of $\pi$ depends only on the partition of $X$ into immersed leaves.
\begin{rem} As is stressed in \cite{AndroulidakisZambon17}, in general projectivity of a singular Lie subalgebroid $\rho_{A}(\Gamma(\mc{A})) \subseteq \Gamma(TX)$ cannot be tested from just the orbits of $\mc{A}$. However, this is possible for those singular Lie subalgebroids arising from Poisson structures due to skew-symmetry of the bivector $\pi$, which implies that the anchor $\pi^\sharp$ of $T^*_\pi X$ is skew, so that $\ker(\pi^\sharp_x) = {\rm im}(\pi^\sharp_x)^\circ$.
\end{rem}
By continuity and density of $X_{\pi, {\rm max}}$, the distribution $D_\pi$ in the above proposition is unique and involutive, and determines a $2k$-dimensional regular foliation by $\pi$-Poisson submanifolds. Phrased differently, $D_\pi$ specifies a $2k$-regular Lie algebroid whose orbits are $\pi$-Poisson submanifolds. Note that $X_{\pi, {\rm max}} = X_{\pi, {\rm reg}}$ for almost-regular Poisson structures. There is another characterization of almost-regularity, phrased in terms of $\pi$ itself. Let ${\rm Poiss}_{2k}(X) \subseteq {\rm Poiss}(X)$ be the space of Poisson structures of maximal rank $2k$, i.e.\ for which $X_{\pi,{\rm max}} = X \backslash X_{2k-2}$.
\begin{prop}[{\cite[Proposition 2.11]{AndroulidakisZambon17}}]\label{prop:almostregularline} Let $\pi \in {\rm Poiss}_{2k}(X)$. Then $\pi$ is almost-regular if and only if $X_{\pi, {\rm max}}$ is dense and there exists a line bundle $K \subseteq \wedge^{2k} TX$ such that $\wedge^k \pi \in \Gamma(K)$.
\end{prop}
The above distribution $D_\pi$ and line bundle $K$ are related via $\det(D_\pi) = \wedge^k D_\pi = K$. Note that because for $\pi$ as in \autoref{prop:almostregularline} we have $X_{\pi,{\rm max}} = X_{\pi,{\rm reg}} = X \backslash X_{2k-2} = (\wedge^{k} \pi)^{-1}(0)$, so that density of this subspace amounts to saying that the pair $(K,\wedge^k \pi)$ is a divisor on $X$.
\begin{rem} In work of Lanius \cite{Lanius16two,Lanius16,Lanius17}, specifically for computing Poisson cohomology, the \emph{rigged algebroid} $\mc{R}_\pi$ of a Poisson structure $\pi$ plays a large role. It is defined as the almost-injective Lie algebroid inducing $\mc{F}_\pi = \pi^\sharp(\Gamma(T^*X))$. In the joint paper \cite{KlaasseLanius18} we systematically develop its use in computing the Poisson cohomology of almost-regular Poisson structures.
\end{rem}
Together with \autoref{defn:adivisortype} we obtain the following dictionary between concepts.
\begin{cor}\label{cor:divtypealmostreg} Let $\pi \in {\rm Poiss}(X)$ be given. Then the following are equivalent:
	\bi
	\item $\pi$ is of $m$-divisor-type for some $m \geq 0$;
	\item $\pi$ is almost-regular;
	\item The rigged algebroid $\mc{R}_\pi$ exists.
	\ei
\end{cor}
By defining almost-regularity for $\mc{A}$-Poisson structures as projectivity of $\pi^\sharp_{\mc{A}}(\Gamma(\mc{A}^*))$, and its \emph{$\mc{A}$-rigged algebroid} $\mc{R}_{\pi_{\mathcal{A}}}$ via $\Gamma(\mc{R}_{\pi_{\mathcal{A}}}) \cong \pi^\sharp_{\mc{A}}(\Gamma(\mc{A}^*))$ \cite{KlaasseLanius18}, \autoref{cor:divtypealmostreg} is seen to hold for $\mc{A}$-Poisson structures for any Lie algebroid $\mc{A}$, after adapting the proof of \cite[Proposition 2.11]{AndroulidakisZambon17}.
\subsection{Lifting Poisson structures}
\label{sec:liftingpoisson}
We next introduce the process of lifting $\mc{A}$-Poisson structures through base-preserving Lie algebroid morphisms (see also the discussion in \cite{CavalcantiKlaasse18}).
\begin{defn}\label{defn:poissonalift} Let $(\varphi,{\rm id}_X)\colon \mathcal{A} \to \mathcal{A}'$ be a Lie algebroid morphism. An $\mc{A}'$-bivector $\pi_{\mathcal{A}'} \in \Gamma(\wedge^2 \mc{A}')$ is \emph{$\mc{A}$-liftable} if there exists an $\mc{A}$-bivector $\pi_{\mc{A}} \in \Gamma(\wedge^2 \mc{A})$ such that $\varphi(\pi_{\mc{A}}) = \pi_{\mc{A}'}$.
\end{defn}
In the above situation we again call $\pi_{\mc{A}}$ the \emph{$\mc{A}$-lift} of $\pi_{\mc{A}'}$, and say that $\pi_{\mc{A}'}$ is of \emph{$\mc{A}$-type}. We will denote the space of $\mc{A}$-liftable bivectors by $\mf{X}_\mc{A}^2(X) := \rho_\mc{A}(\mf{X}^2(\mc{A})) \subset \mf{X}^2(X)$, and the space of $\mc{A}$-liftable Poisson structures by ${\rm Poiss}_{\mc{A}}(X) \subseteq {\rm Poiss}(X)$.
\begin{rem} In \cite[Definition 2.16]{Lanius16} this notion is considered using slightly different terminology: relating it to ours, there $\pi \in {\rm Poiss}(X)$ is called $\mc{A}$-Poisson if it is of $\mc{A}$-type. It seems preferable to describe Poisson structures by their divisor (c.f.\ \autoref{defn:adivisortype}) instead of to which Lie algebroid they can be (nondegenerately) lifted, as the latter need not be unique.
\end{rem}
The $\mc{A}$-lift of a $\mc{A}'$-Poisson structure of $\mc{A}$-type need not necessarily be $\mc{A}$-Poisson (consider $\mc{A}$ with trivial anchor and $\pi \equiv 0$). However, this is true when $\varphi$ is of divisor-type, as will be the case for us (moreover, the $\mc{A}$-lift will then be unique). More generally this holds for almost-injective $\varphi$. We have $\pi_{\mc{A}'}^\sharp = \varphi \circ \pi_\mc{A}^\sharp \circ \varphi^*$ as maps, summarized in the following diagram.
\begin{center}
	\begin{tikzpicture}
	\matrix (m) [matrix of math nodes, row sep=2.5em, column sep=2.5em,text height=1.5ex, text depth=0.25ex]
	{	\mc{A}^* & \mc{A} \\ \mc{A}'{}^* & \mc{A}' \\};
	\path[-stealth]
	(m-1-1) edge node [above] {$\pi_\mc{A}^\sharp$} (m-1-2)
	(m-2-1) edge node [left] {$\varphi^*$} (m-1-1)
	(m-2-1) edge node [above] {$\pi_{\mc{A}'}^\sharp$} (m-2-2)
	(m-1-2) edge node [right] {$\varphi$} (m-2-2);
	\end{tikzpicture}
\end{center}
\begin{defn}\label{defn:apoissonmap} An \emph{$\mc{A}$-Poisson map} $(\varphi,f)\colon (\mc{A},X,\pi_\mc{A}) \to (\mc{B},X',\pi_\mc{B})$ is a Lie algebroid morphism between Lie algebroids with $\mc{A}$-Poisson structures such that $\varphi(\pi_\mc{A}) = \pi_\mc{B}$.
\end{defn}
As for Poisson maps, it is not true that any $\mc{A}$-Poisson structure can be pushed forward to an $\mc{A}'$-Poisson structure along a given Lie algebroid morphism $(\varphi,f)\colon \mc{A} \to \mc{A}'$. However, this is always possible along base-preserving morphisms. In particular, given an $\mc{A}$-Poisson structure, the bivector $\pi := \rho_{A}(\pi_{\mc{A}})$ is always Poisson. The lifting condition of an $\mc{A}'$-Poisson structure to an $\mc{A}$-Poisson structure is exactly that $(\varphi,{\rm id}_X)\colon (\mc{A},\pi_\mc{A}) \to (\mc{A}',\pi_{\mc{A}'})$ is an $\mc{A}$-Poisson map. If $\pi_{\mc{A}'}$ has a nondegenerate $\mc{A}$-Poisson lift $\pi_\mc{A}$, we say $\pi_{\mc{A}'}$ is of \emph{nondegenerate $\mc{A}$-type}, and similarly for other types (e.g.\ for regular ones). We often omit writing the Lie algebroid morphism $\varphi$ (as typically it is induced by an inclusion of modules), and always use the anchor map $\rho_{\mc{A}}$ when considering $\mc{A}$-lifts of bivectors on $TX$. When $\pi_{\mc{A}}$ nondegenerately lifts $\pi_{\mathcal{A}'}$, we will call its inverse $\omega_{\mc{A}} := \pi_\mc{A}^{-1} \in \Omega^2(\mc{A})$ the \emph{dual $\mc{A}$-form} to $\pi_{\mathcal{A}'}$. This is an $\mc{A}$-symplectic structure in the sense of \autoref{sec:symplecticliealgebroids}, making $\mc{A}$ into a symplectic Lie algebroid \cite{NestTsygan01}.

Thus the process of lifting constitutes the existence of a base-preserving $\mc{A}$-Poisson map. In terms of $\mc{A}$-Poisson algebroids, we obtain the following analogue of the well-known fact that a Poisson map induces a Lie algebroid comorphism between the respective Poisson algebroids.
\begin{prop}\label{prop:apoissoncomorphism} Let $(\varphi,f)\colon (\mc{A},X,\pi_\mc{A}) \to (\mc{B},X',\pi_\mc{B})$ be an $\mc{A}$-Poisson map. Then there is an induced $\mc{A}$-Lie algebroid comorphism $(\varphi;\varphi,f)\colon (\mc{A}^*_{\pi_\mc{A}},\mc{A},X) \dashedrightarrow (\mc{B}^*_{\pi_\mc{B}},\mc{B},X')$.
\end{prop}
The lifting process interacts with Poisson structures of divisor-type as follows.
\begin{prop}\label{prop:apoissondivtypelifts} Let $(\varphi,{\rm id}_X)\colon \mathcal{A} \to \mathcal{A}'$ be a Lie algebroid morphism of divisor-type, and let $\pi_{\mc{A}'}$ be an $\mc{A}'$-Poisson structure of $m$-divisor-type for some $m \geq 0$ that is also of $\mc{A}$-type. Then its unique $\mc{A}$-lift $\pi_\mc{A}$ is of $m$-divisor-type, and ${\rm div}(\pi_{\mc{A}'}) = {\rm div}(\pi_\mc{A}) \cdot {\rm div}(\varphi)$ so that $I_{\pi_{\mc{A}'}} = I_{\pi_{\mc{A}}} \cdot I_\varphi$.
\end{prop}
\bp Let $\pi_\mc{A}$ be an $\mc{A}$-lift of $\pi_{\mc{A}'}$, so that $(\wedge^2 \varphi)(\pi_\mc{A}) = \pi_{\mc{A}'}$. We have ${\rm div}(\pi_{\mathcal{A}'}) = (K, \wedge^m \pi_{\mc{A}'})$ for some line bundle $K \subseteq \wedge^{2m} \mc{A}'$. As $\wedge^{2m} \varphi\colon \wedge^{2m} \mc{A} \to \wedge^{2m} \mc{A}'$, we find a line bundle $K_\mc{A} \subseteq \wedge^{2m} \mc{A}$ for which $\wedge^{2m}\varphi\colon K_\mc{A} \to K$ and $\wedge^m \pi_\mc{A} \in \Gamma(K_\mc{A})$. If $\pi_\mc{A}$ vanishes then so must $\pi_{\mc{A}'}$, so that $(K_\mc{A}, \wedge^m \pi_\mc{A})$ is a divisor, ${\rm div}(\pi_{\mathcal{A}})$, as a subset of a nowhere dense set is nowhere dense. The statement that ${\rm div}(\pi_{\mc{A}'}) = {\rm div}(\pi_\mc{A}) \cdot {\rm div}(\varphi)$ then follows from the relation $\pi_{\mc{A}'}^\sharp = \varphi \circ \pi_\mc{A}^\sharp \circ \varphi^*$.
\ep
This proposition makes precise that the process of lifting makes $\mc{A}$-Poisson structures less degenerate, as their degeneracies are absorbed in that of the Lie algebroid morphism (hence, in the Lie algebroid it is lifted to). Moreover, it says that $\pi_{\mc{A}'}$ is of $m$-regular $\mc{A}$-type if and only if ${\rm div}(\pi_{\mc{A}'}) = {\rm div}(\varphi)$, once $\pi_{\mc{A}'}$ is of $\mc{A}$-type (this includes being of nondegenerate $\mc{A}$-type. There is a similar statement for not necessarily base-preserving $\mc{A}$-Poisson maps.
\begin{lem} Let $(\varphi,f)\colon (\mc{A}^{2n},X,\pi_{\mc{A}}) \to (\mc{B}^{2n},X',\pi_{\mathcal{B}})$ be an $\mc{A}$-Poisson map between Lie algebroids with $\mc{A}$-Poisson structures of divisor-type. Then $I_{\pi_{\mathcal{A}}}$ divides $f^*I_{\pi_{\mathcal{B}}}$.
\end{lem}
\bp As ${\rm rank}(\mc{A}) = {\rm rank}(\mc{B})$, the fact that $\varphi(\pi) = \sigma$ implies that $\det(\varphi)(\pi_{\mc{A}}^n) = \pi_{\mathcal{B}}^n$. Take local volume forms $\mu_\mc{A}, \mu_\mc{B}$ for $\mc{A}$ and $\mc{B}$. For certain functions $g \in C^\infty(X)$, $h \in C^\infty(X')$:
\be
	\langle\pi_\mc{A}^n,\mu_\mc{A}\rangle = g \text{, so } I_{\pi_{\mathcal{A}}} = \langle g \rangle, \qquad \text{and} \qquad \langle\pi_{\mc{B}}^n,\mu_\mc{B}\rangle = h \text{, so } I_{\pi_{\mathcal{B}}} = \langle h \rangle.
\ee
Note that $\varphi^* \mu_\mc{B} = w \mu_\mc{A}$ for some nonnegative function $w \in C^\infty(X)$. Thus $I_{\pi_{\mathcal{A}}} \subseteq f^*I_{\pi_{\mathcal{B}}}$, as
\be
f^*(h) = f^*\langle\pi_{\mc{B}},\mu_\mc{B}\rangle = f^*\langle\det(\varphi)(\pi_{\mc{A}}^n),\mu_\mc{B}\rangle = \langle\pi_\mc{A}^n, \varphi^*\mu_\mc{B}\rangle = \langle\pi_\mc{A}^n, w \mu_\mc{A}\rangle = w \langle \pi_\mc{A}^n, \mu_\mc{A}\rangle = w g.\qedhere
\ee
\ep
The above lemma is of particular interest in the standard case when $\mc{A} = TX$ and $\mc{B} = TX'$.
For clarity, in the specific case that $\mc{A}' = TX$, \autoref{prop:apoissondivtypelifts} says that a Poisson structure $\pi \in {\rm Poiss}(X)$ of divisor-type $I$ can only admit nondegenerate $\mc{A}$-lifts to Lie algebroids satisfying $I_\mc{A} = I$.
We can iterate the lifting procedure as follows, whose proof is immediate.
\begin{prop}\label{prop:liftiterate} Let $(\varphi,{\rm id}_X)\colon \mc{A} \to \mc{A}'$ and $(\varphi',{\rm id}_X)\colon \mc{A}' \to \mc{A}''$ be Lie algebroid morphisms, and $\pi_{\mc{A}''} \in \rm{Poiss}(\mc{A}'')$. Then if $\pi_{\mc{A}''}$ is of $\mc{A}$-type, it is of $\mc{A}'$-type. If $\pi_\mc{A}$ is an $\mc{A}$-Poisson lift of $\pi_{\mc{A}''}$, then $\pi_{\mc{A}'} := \varphi(\pi_\mc{A})$ is an $\mc{A}'$-Poisson lift of $\pi_{\mc{A}''}$ which itself is of $\mc{A}$-type.
\end{prop}
Consequently, an $\mc{A}'$-Poisson structure $\pi_{\mathcal{A}'}$ being of $\mathcal{A}$-type is the same thing as its underlying Poisson structure $\pi := \rho_{\mathcal{A}'}(\pi_{\mathcal{A}'})$ being of $\mathcal{A}$-type. Let $\pi_{\mc{A}'}$ be a nondegenerate $\mc{A}'$-Poisson structure. Then any $\mc{A}$-Poisson map $(\varphi,f)\colon (\mc{A},\pi_\mc{A}) \to (\mc{A}',\pi_{\mc{A}'})$ has to be a Lie algebroid submersion, i.e.\ the morphism $\varphi$ must be fiberwise surjective. This implies the following.
\begin{prop}\label{prop:liftingnondegiso} Let $\pi \in {\rm Poiss}(X)$ be of nondegenerate $\mc{A}' $-type and let  $(\varphi,{\rm id}_X)\colon \mc{A} \to \mc{A}'$ be a morphism between Lie algebroids of the same rank. Then the following are equivalent:
	\bi
	\item $\varphi$ is an isomorphism;
	\item $\pi$ is of $\mc{A}$-type;
	\item $\pi$ is of nondegenerate $\mc{A}$-type.
	\ei
\end{prop}
As a consequence, once one has lifted $\pi$ to being nondegenerate, one cannot (meaningfully) lift further. This makes sense, as the lifting process is meant to desingularize the Poisson structure, and nondegenerate ones are maximally nonsingular. This is consistent with \autoref{prop:apoissondivtypelifts}, as isomorphisms specify trivial divisors. Moreover, we see that $\pi_{\mc{A}'}$ is of nondegenerate $\mc{A}$-type if and only if the underlying Poisson structure $\pi = \rho_{\mc{A}'}(\pi_{\mc{A}'})$ is.
\begin{rem}\label{rem:multiplelifts} A Poisson structure can lift nondegenerately to multiple non-isomorphic Lie algebroids of (necessarily) the same divisor-type. The orbits of these Lie algebroids will specify the same partition of $X$ into leaves. For some consequences of this phenomenon, see \cite{KlaasseLanius18}.
\end{rem}
\begin{exa}\label{exa:mregularapoisson} Let $\pi_\mc{A}$ be an $m$-regular $\mc{A}$-Poisson structure. Then its image $D_{\pi_\mc{A}} = \pi_\mc{A}^\sharp(\mc{A}^*)$ is a regular $\mc{A}$-distribution. By \autoref{exa:adistribution}, $D_{\pi_{\mc{A}}}$ is a Lie subalgebroid of $\mc{A}$ and thus carries an $\mc{A}$-Lie algebroid structure $\varphi\colon D_{\pi_\mc{A}} \to \mc{A}$. It follows that $\pi_{\mc{A}}$ is of nondegenerate $D_{\pi_\mc{A}}$-type.
\end{exa}
Let us give another consequence of lifting. Let $\mc{A} \to X$ be a Lie algebroid and recall that $\mc{F}_\mc{A} = \rho_{\mc{A}}(\Gamma(\mc{A}))$ is its induced singular foliation whose leaves are the orbits of $\mc{A}$.
\begin{prop}\label{prop:atypepoissonsubmfd} Let $\pi \in {\rm Poiss}_\mc{A}(X)$. Then the orbits of $\mc{A}$ are $\pi$-Poisson submanifolds.
\end{prop}
\bp The lifting condition $\rho_{\mathcal{A}}(\pi_{\mc{A}}) = \pi$ equivalently reads $\pi^\sharp = \rho_{\mc{A}} \circ \pi_\mc{A}^\sharp \circ \rho_{\mc{A}}^*$. As a submanifold $N$ is $\pi$-Poisson if $\pi^\sharp(T^*_x X) \subseteq T_x N$ for all $x \in N$, we note that this condition pointwise becomes $\rho_{\mc{A},x}(\mc{B}_x) \subseteq T_x N$ for the subspace $\mc{B}_x := \pi_{\mc{A},x}^\sharp(\rho_{\mc{A},x}^*(T^*_x X)) \subseteq \mc{A}_x$. If $N$ is an orbit of $\mc{A}$, then by definition $T_x N = {\rm im}(\rho_{\mc{A},x}) \subseteq T_x X$, so that this condition is certainly satisfied.
\ep
This gives an intuitive description of the lifting process when $\mc{A}$ is almost-injective. Namely, a choice of lifting algebroid $\mc{A}$ amounts to a ``grouping'' of the $\pi$-symplectic leaves into orbits of $\mc{A}$. Unless $\mc{A}$ is regular (with injective anchor), there may be multiple non-isomorphic Lie algebroids with the same orbits (cf.\ \autoref{rem:multiplelifts}). There is a partial converse to \autoref{prop:atypepoissonsubmfd} in the regular injective case, in which case $\mc{A}$ is just a regular foliation $T\mc{F} \subseteq TX$.
\begin{prop}\label{prop:regularpipoissonlift} Let $\pi \in {\rm Poiss}(X)$ and $\mc{F}$ be a regular foliation by $\pi$-Poisson submanifolds, with $D = T\mc{F}$ the associated regular Lie algebroid. Then $\pi$ is of $D$-type.
\end{prop}
\bp Let $L$ be a leaf of $\mc{F}$. As it is $\pi$-Poisson, we have $\pi_x \in \wedge^2 T_x L = \wedge^2 D_x$ for all $x \in L$. By regularity there exists $\pi_D \in {\rm Poiss}(D)$ agreeing with $\pi$ on each leaf, which is a $D$-lift of $\pi$.
\ep
This result recovers the fact that a regular Poisson structure is equivalently described as a nondegenerate Poisson bivector for the Lie algebroid of tangencies to its symplectic leaves. The consequences of \autoref{prop:regularpipoissonlift} for almost-regular Poisson structures (i.e., those of $m$-divisor-type), in light of \autoref{prop:almostregulardistr}, are explored in \autoref{sec:liftingalmostreg}.
\begin{rem} \autoref{prop:regularpipoissonlift} is false when $\mc{F}$ is no longer regular, as the lifting property cannot be tested purely by looking at the underlying singular foliation of the Lie algebroid. For example, on $\R^2$ the Lie algebroids $\mc{A}_Z^k = \langle x^k \partial_x, \partial_y \rangle$ associated to the jets of $x \in C^\infty(\R^2)$ for $Z = \{x = 0\}$ all have the same orbits, but $\pi = x \partial_x \wedge \partial_y$ only admits $\mc{A}_Z^k$-lifts for $k = 0$ or $1$. Essentially what is used in the regular case is that a decomposition into leaves determines a unique module of vector fields if the decomposition is regular (see the discussion in \cite{AndroulidakisZambon16}).
\end{rem}
\subsection{Lifting Poisson structures of divisor-type}
\label{sec:liftingdivtype}
In this section we discuss the process of lifting Poisson structures of divisor-type. Let $I \subseteq C^\infty(X)$ be a divisor ideal, and assume that the associated ideal Lie algebroid $TX_I \to X$ exists. In other words, assume that the divisor ideal $I$ is projective. Recall that given a Lie algebroid $\mc{A}$ we denote the spaces of $\mc{A}$-liftable and $I$-divisor-bivectors by $\mf{X}^2_{\mc{A}}(X)$ and $\mf{X}^2_{I}(X)$ respectively. Similarly, we have their Poisson counterparts ${\rm Poiss}_{\mc{A}}(X) \subseteq \mf{X}^2_{\mc{A}}(X)$ and ${\rm Poiss}_{I}(X) \subseteq \mf{X}^2_I(X)$.

We wish to study the relation between the spaces ${\rm Poiss}_{TX_I}(X)$ and ${\rm Poiss}_{I}(X)$. Note that for $\pi \in {\rm Poiss}_I(X)$ to be liftable to $\mc{A}$, by \autoref{prop:apoissondivtypelifts} $I$ must be $I$ divisible by $I_\mc{A}$. Consequently whether $I$ is standard (\autoref{defn:divprojstandard}) is relevant, i.e.\ whether $I_{TX_I} = I$.

Divisor bivectors are generally not liftable to their ideal Lie algebroid, i.e.\ $\mf{X}^2_I(X) \not\subseteq \mf{X}^2_{TX_I}(X)$. The other inclusion is also false, $\mf{X}^2_{TX_I}(X) \not\subseteq \mf{X}^2_I(X)$, as the following examples show.
\begin{exa}\label{exa:logbivector} The bivector $\pi = x \partial_x \wedge \partial_y + \partial_z \wedge \partial_w + \partial_x \wedge \partial_w$ on $\R^4$ with coordinates $(x,y,z,w)$ satisfies $\wedge^2 \pi = x \partial_x \wedge \partial_y \wedge \partial_z \wedge \partial_w$, which vanishes transversally on $Z_\pi = \{x = 0\}$. Hence $\pi$ is of log divisor-type, but does not lift to $\mc{A}_{Z_\pi} = \langle x \partial_x, \partial_y, \partial_z, \partial_w \rangle$ due to the presence of $\partial_x \wedge \partial_w$, showing that $\mf{X}^2_I(X) \not\subseteq \mf{X}^2_{TX_I}(X)$ in this case. Note that $\pi$ is not Poisson, as $[\pi,\pi] = -\partial_x \wedge \partial_y \wedge \partial_w \neq 0$.
\end{exa}
\begin{exa} The Poisson bivector $\pi = x^2 \partial_x \wedge \partial_y$ on $\R^2$ with coordinates $(x,y)$ lifts to $\mc{A}_Z$, where $Z = \{x = 0\}$, with $\mc{A}_Z$-lift $\pi_{\mc{A}_Z} = x (x \partial_x) \wedge \partial_y$, but it does not lift nondegenerately, nor does $\pi$ specify a log divisor structure on $Z$. Thus we see that $\mf{X}^2_{TX_I}(X) \not\subseteq \mf{X}^2_I(X)$ in this case.
\end{exa}
The Poisson condition is important for liftability as the above examples show. To study this in more detail, we first discuss which vector fields are liftable to $TX_I$.
\begin{prop}\label{prop:liftifpreserve} Let $\pi \in \mf{X}^2_I(X)$. A vector field $V \in \mf{X}(X)$ lies in $\mf{X}_{TX_I}(X)$ if $\mc{L}_V (\wedge^n \pi) = 0$.
\end{prop}
\bp Denote $\wedge^n \pi$ by $\pi^n$. By definition of $\Gamma(TX_I) = \Gamma(TX)_I$, we must show that $\mc{L}_V I \subset I$, which is a local statement. Let ${\rm vol}_X$ be a local volume form and define $f := \pi^n({\rm vol}_X)$ so that $I = \langle f \rangle$. Then $\mc{L}_V f = \mc{L}_V \pi^n({\rm vol}_X) = (\mc{L}_V \pi^n)({\rm vol}_X) + \pi^n \mc{L}_V \rm{vol}_X$. Note that $\mc{L}_V {\rm vol}_X = g {\rm vol}_X$ for some local function $g$. From this we conclude that
\be
	\mc{L}_V f= \mc{L}_V \pi^n ({\rm vol}_X) + \pi^n(g {\rm vol}_X) = \mc{L}_V \pi^n ({\rm vol}_X) + g \pi^n({\rm vol}_X) = \mc{L}_V \pi^n ({\rm vol}_X) + g f.
\ee
As $g f \in I$, we see that $V$ lifts if $\mc{L}_V \pi^n = 0$.
\ep
The assumption in \autoref{prop:liftifpreserve} seems strong, yet allows for the following conclusion.
\begin{cor}\label{cor:poissonhamiltonian} Let $\pi \in {\rm Poiss}_I(X)$. Then ${\rm Ham}_\pi(X) \subseteq {\rm Poiss}_{\pi}(X) \subseteq \mf{X}_{TX_I}(X)$.
\end{cor}
\bp That ${\rm Ham}_\pi(X) \subseteq {\rm Poiss}_{\pi}(X)$ is clear. Given $V \in {\rm Poiss}_\pi(X)$ we have $\mc{L}_{V} \pi = 0$. But then we compute that $\mc{L}_{V} \pi^n = n (\mc{L}_V \pi) \wedge \pi^{n-1} = 0$, so that $V \in \mf{X}_{TX_I}(X)$ by \autoref{prop:liftifpreserve}.
\ep
\begin{rem} For a general bivector $\pi \in \mf{X}^2(X)$, the definition of its `Hamiltonian vector fields' $V_f = \pi^\sharp(df)$ for $f \in C^\infty(X)$ still makes sense. However, the fact that all such vector fields preserve $\pi$, i.e.\ $\mc{L}_{V_f} \pi = 0$, is equivalent to $\pi$ being Poisson.
\end{rem}
Due to \autoref{cor:poissonhamiltonian}, we see that for $\pi \in {\rm Poiss}_I(X)$ there exists a map $\wt{\pi}^\sharp\colon T^*X \to TX_I$ fitting in the following diagram, which is not yet the existence of an $TX_I$-lift.
\begin{center}
	\begin{tikzpicture}
	\matrix (m) [matrix of math nodes, row sep=2.5em, column sep=2.5em,text height=1.5ex, text depth=0.25ex]
	{	TX^*_I & TX_I \\ T^*X & TX \\};
	\path[-stealth]
	(m-2-1) edge node [left] {$\rho_{TX_I}^*$} (m-1-1)
	(m-2-1) edge node [above] {$\pi^\sharp$} (m-2-2)
	(m-2-1) edge node [above] {$\wt{\pi}^\sharp$} (m-1-2)
	(m-1-2) edge node [right] {$\rho_{TX_I}$} (m-2-2);
	\draw[dotted, ->] (m-1-1) to (m-1-2);
	\end{tikzpicture}
\end{center}
One can obtain an actual $TX_I$-lift under the assumption that $\Gamma(TX_I^*)$ admits local bases of closed sections. Namely, we can now prove \autoref{thm:introlifting} from the introduction.
\begin{thm}\label{thm:ipoissonaitype} Let $I \subseteq C^\infty(X)$ be a projective divisor ideal with $TX^*_I$ admitting local bases of closed sections. Then ${\rm Poiss}_I(X) \subseteq {\rm Poiss}_{TX_I}(X)$, i.e.\ those of $I$-divisor-type are of $TX_I$-type.
\end{thm}
\bp We obtain from the discussion above a map $\wt{\pi}^\sharp\colon T^*X \to TX_I$. We can dualize the bundle morphism $\wt{\pi}^\sharp$ to a map $(\wt{\pi}^\sharp)^*\colon TX_I^* \to TX$. To test whether this lifts to a map $\pi_{TX_I}^\sharp\colon TX_I^* \to TX_I$, we proceed as in \autoref{cor:poissonhamiltonian} by checking whether $(\wt{\pi}^\sharp)^*$ maps to $\Gamma(TX)_I \cong \Gamma(TX_I)$ on sections. Assuming that $\Gamma(TX_I^*)$ admits local bases of closed sections, we can by continuity answer this in the isomorphism locus $X \backslash Z_I$, where such local sections are exact under the isomorphism with $TX$, and where $(\wt{\pi}^\sharp)^* = (\pi^\sharp)^*$ using the isomorphism given by $\rho_{TX_I}$. Here the lifting property follows as in \autoref{cor:poissonhamiltonian}, using that $(\pi^\sharp)^* = - \pi^\sharp$ by skew-symmetry. The lifting property then holds in the entirety of $X$ by density of the isomorphism locus.
\ep
\begin{rem} By inspecting the proof of \autoref{thm:ipoissonaitype}, we readily see that it still holds for Poisson structures of $m$-$I$-divisor-type, that is, any almost-regular Poisson structure.
\end{rem}
The converse to the above theorem is false unless $\pi \in {\rm Poiss}(X)$ is of nondegenerate $TX_I$-type and $I$ is standard (see the discussion below \autoref{prop:apoissondivtypelifts}). Further, we know that if $I$ is standard and $\pi \in {\rm Poiss}_I(X)$ is of $TX_I$-type, then $\pi$ is in fact of nondegenerate $TX_I$-type.
\begin{rem} Following up on \autoref{rem:nonstandard}, let $I \subseteq C^\infty(X)$ be a nontrivial projective divisor ideal and $\pi \in {\rm Poiss}_{I^k}(X)$ for some $k > 1$. While $\pi$ can be of $TX_I$-type using \autoref{thm:ipoissonaitype}, it can never be of nondegenerate $TX_I$-type by comparing the divisor ideals of $\pi$ and $TX_I$.
\end{rem}
Note that if $\pi$ is of $TX_I$-type, the degeneracy locus $Z_I$ must be a $\pi$-Poisson subset, i.e.\ $I_{Z_I}$ is a $\pi$-Poisson ideal. Similarly we have the following (recall that $I \subseteq I_{Z_I}$ by \autoref{prop:divvanishingideal}).
\begin{prop}\label{prop:itypeipoissonideal} Let $I$ be a divisor ideal and $\pi \in {\rm Poiss}_I(X)$. Then $I$ is a $\pi$-Poisson ideal.
\end{prop}
\bp Let ${\rm vol}_X$ be a local volume form, and $g \in C^\infty(X)$. Then $I = \langle \pi^n({\rm vol}_X) \rangle = \langle f \rangle$, and
\be
	\{g, f\}_\pi = \{g, \pi^n({\rm vol}_X)\}_\pi = \pi(dg, d(\pi^n{\rm vol}_X)) = \pm \mc{L}_{V_g} (d \pi^n({\rm vol}_X)) = \pi^n \mc{L}_{V_g} {\rm vol}_X,
\ee
using that as $V_g$ is Hamiltonian, it preserves $\pi^n$. We have $\mc{L}_{V_g} {\rm vol}_X = h {\rm vol}_X$ for some function $h \in C^\infty(X)$, so that $\pi^n \mc{L}_{V_g} {\rm vol}_X = h f \in I$. This implies that $\{g,f\}_\pi \in I$ as desired.
\ep
\subsection{Examples}
In this section we discuss several interesting examples of Poisson structures of divisor-type, and of the lifting procedure to their associated ideal Lie algebroids. We will only do so here for Poisson structures in the usual sense, as this is the main motivation for the development of the theory. However, it is quite possible, and certainly interesting, to consider $\mc{A}$-Poisson structures of divisor-type for given, well-understood, Lie algebroids $\mc{A} \to X$, in cases where it is difficult to intrinsically capture when a Poisson structure is of $\mc{A}$-type. Recall the examples of divisors and ideal Lie algebroids given in \autoref{sec:divexamples} and \autoref{exa:ideallalgebroids}.

Throughout, this section let $X$ be a given $2n$-dimensional manifold.
\begin{exa}[Nondegenerate] Let $\pi \in {\rm Poiss}(X)$ be a nondegenerate Poisson structure. Then $\pi$ is of $C^\infty(X)$-divisor-type, i.e.\ $I_\pi$ is the trivial divisor ideal. By \autoref{prop:liftingnondegiso}, $\pi$ is only of $\mc{A}$-type for the Lie algebroid $\mc{A} = TX$, and certainly lifts to it nondegenerately.
\end{exa}
\begin{exa}[Log-Poisson, \cite{GuilleminMirandaPires14}]\label{exa:logpoissonlift} Let $\pi \in {\rm Poiss}(X)$ be a \emph{log-Poisson structure}, i.e.\ such that it is of log divisor-type. Then $I_\pi = I_Z$ is a log divisor ideal. Regarding liftability of $\pi$, recall from \autoref{exa:zvx} that the log-tangent bundle $\mc{A}_Z$ is the (primary) ideal Lie algebroid of $I_Z$. By \autoref{thm:ipoissonaitype} we see that $\pi$ is of $\mc{A}_Z$-type. As $I_Z$ is standard, it is in fact of nondegenerate $\mc{A}_Z$-type, so that it is dual to a \emph{log-symplectic structure} \cite{GuilleminMirandaPires14,GualtieriLi14,Cavalcanti17,MarcutOsornoTorres14}.
\end{exa}
We can determine the lifting property to log-tangent bundles $\mc{A}_Z \to X$ more generally.
\begin{prop}\label{prop:azliftpoisson} Let $\pi \in {\rm Poiss}(X)$ be given and let $Z \subseteq X$ be a hypersurface. Then the following are equivalent:
\bi
	\item $Z$ is a $\pi$-Poisson submanifold;
	\item $I_Z$ is a $\pi$-Poisson ideal;
	\item $\pi$ of $\mc{A}_Z$-type.
\ei
\end{prop}
\bp This is a rephrasing of \cite[Proposition 4.4.1]{Pym13}. The equivalence of the first two points is standard, as $Z$ is a closed submanifold. If $\pi$ of $\mc{A}_Z$-type then $Z$ is $\pi$-Poisson, because $Z$ is the degeneracy locus of $\mc{A}_Z$ (c.f.\ \autoref{prop:pidegenpipoisson}), or because it is an orbit of $\mc{A}_Z$ (c.f.\ \autoref{prop:atypepoissonsubmfd}). To see that if $Z$ is $\pi$-Poisson, then $\pi$ is of $\mc{A}_Z$-type, we turn to the proof of \autoref{thm:ipoissonaitype}: there we showed first that each $\pi$-Hamiltonian vector field is liftable to $\mc{A}_Z$, before dualizing to obtain liftability of $\pi$. However, $Z$ being $\pi$-Poisson is equivalent to all $\pi$-Hamiltonian vector fields being tangent to $Z$ (or preserving the vanishing ideal $I_Z$).
\ep
From \autoref{prop:azliftpoisson} one can wonder when a Poisson structure $\pi$ is of nondegenerate $\mc{A}_Z$-type. However, this is already answered by \autoref{prop:apoissondivtypelifts}: for this $\pi$ must be of $I_Z$-divisor-type (because $I_{\mc{A}_Z} = I_Z$), which is also sufficient by, for example, \autoref{prop:itypeipoissonideal}. Hence, $\pi$ is of nondegenerate $\mc{A}_Z$-type if and only if it of $I_Z$-divisor-type, i.e.\ is log-Poisson.
\begin{rem} Because $I_Z$ is a locally principal (divisor) ideal, \autoref{prop:azliftpoisson} says that Poisson structures of $I$-divisor-type are liftable to $\mc{A}_Z = TX_{I_Z}$ if $I_Z$ divides $I$. The converse statement is the content of \autoref{prop:apoissondivtypelifts}. It would be interesting to determine whether there are other Lie algebroids of divisor-type for which a similar statement is true, or to show that this only happens in this case. We seem to crucially use that $I_Z$ is not only a divisor ideal, but in fact the vanishing ideal of its support. This only happens in the log setting.
\end{rem}
We can determine the lifting property for log-bivectors following \cite[Lemma 138]{Frejlich11}.
\begin{prop}[\cite{Frejlich11}]\label{prop:liftazbivectors} Let $\pi \in \mf{X}^2_{I_Z}(X)$ be a bivector of log-divisor-type. Then $\pi$ is of $\mc{A}_Z$-type if and only if $\pi^\sharp(T^*_x X) \subseteq T_x Z$ for some $x \in Z_i$ for each connected component $Z_i \subseteq Z$. 
\end{prop}
\bp Given $x \in Z$ with local generator $z$ of $I_Z$ we have $\pi = \partial_z \wedge \pi^\sharp(dz) + \nu$ with $\nu^\sharp(dz) \equiv 0$. Then $\wedge^n \pi \sim \partial_z \wedge \pi^\sharp(dz) \wedge \nu^{n-1}$. This decomposition also shows that for all $x \in Z$, the condition $\pi^\sharp(T^*_x X) \subseteq T_x Z$ is equivalent to the condition $\pi^\sharp(d_x I_Z) = 0$. If $\pi^\sharp(d_x I_Z) = 0$, then $\pi^\sharp(dz)|_Z = 0$, so that in turn $\pi^\sharp(dz) = x V$ for some $V \in \Gamma(TX)$. Thus $\pi^\sharp(d\log z) = V$ which shows that $\pi$ is of $\mc{A}_Z$-type (pointwise, that $\pi_x \in \wedge^2 \mc{A}_{Z,x}$). The converse statement is immediate, because $Z$ is an $\mc{A}_Z$-invariant submanifold. Next, the space $\mc{Z}(\pi) = \{x \in Z \, : \, \pi^\sharp(T^*_x x) \subseteq T_x Z\} \subseteq Z$ is both open and closed, where for openness we use that if $x \in \mc{Z}(\pi)$, then $\pi^\sharp(d_x z) = 0$, whence $\nu_x^{n-1} \neq 0$ by transversality. But then there exists some open $U$ around $x$ on which $\nu^{n-1}$ is nonvanishing. As $\wedge^n \pi$ vanishes on $U \cap Z$, this implies that $\pi^\sharp(dz)$ divides $\nu^{n-1}$ on $U \cap Z$. To not contradict transversality of $\wedge^n \pi$, we must have that $\pi^\sharp(dz) = 0$ on $U \cap Z$. This in turn implies that $U \cap Z \subseteq \mc{Z}(\pi)$, which shows openness. The rest follows immediately.
\ep
Essentially, the above proposition uses the fact that transversality is an open condition.
\begin{rem}\label{rem:examplelogbivectorlift} Comparing this to \autoref{exa:logbivector}, we see that there $\pi^\sharp(dw) = \partial_z + \partial_x \not\in TZ$, and indeed the bivector $\pi$ of that example does not lift to $\mc{A}_Z$, as said by \autoref{prop:liftazbivectors}.
\end{rem}
\begin{exa}[Normal-crossing log-Poisson] Let $\pi\in {\rm Poiss}(X)$ be a \emph{normal-crossing log-Poisson structure}, i.e.\ such that it is of normal-crossing log divisor-type, with divisor ideal $I_{\underline{Z}}$. The above discussion (that is, \autoref{prop:azliftpoisson}) also holds for normal-crossing log-Poisson structures: a Poisson structure is of $\mc{A}_{\underline{Z}}$-type if and only if $I_{\underline{Z}}$ is a $\pi$-Poisson ideal, which is true in this case by \autoref{prop:itypeipoissonideal}. Being of nondegenerate $\mc{A}_{\underline{Z}}$-type is again equivalent to being of $I_{\underline{Z}}$-divisor-type. See also \cite[Proposition 4.4.1; Proposition 4.4.2]{Pym13}, and \cite{Klaasse17}.
\end{exa}
\begin{exa}[$b^k$-Poisson, \cite{Scott16}]\label{exa:bkpoissonlift} Another interesting class is given by the \emph{$b^k$-Poisson structures} $\pi_k \in {\rm Poiss}(X)$, which are of $I^k_Z$-divisor-type for $k \geq 0$, where $I_Z$ is a log divisor ideal. This does not characterize $b^k$-Poisson structures. To do so, demand that $\{\cdot,\cdot\}_{\pi_k}$ satisfies
\be
	\{\cdot,\cdot\}_{\pi_k}\colon I_Z \times C^\infty(X) \to I_Z^k.
\ee
If this holds, then $\wedge^n \pi_k$ provides a local generator of $I_Z^k$, and its $(k-1)$-jet can be used to define the $b^k$-tangent bundle $\mc{A}_Z^k$ as in \autoref{exa:bktangentbundle}. It then follows that $\pi_k$ is of nondegenerate $\mc{A}_Z^k$-type, or of nondegenerate \emph{$b^k$-type}. Consequently they are dual to $b^k$-symplectic structures. Note here that $b^1$-Poisson structures are just the log-Poisson structures of \autoref{exa:logpoissonlift}.
\end{exa}
\begin{rem} Our notion of nondegenerate $b^k$-type is almost the same as that of $b^k$-type in \cite{Scott16}. There is one difference, namely that there is a class of (non-canonically) isomorphic Lie algebroids, namely $\mc{A}_Z^k$ with respect to different jet data, and Scott's $b^k$-type demands liftability to any one of these Lie algebroids. We instead extract the jet data from $\pi_k$ itself.
\end{rem}
\begin{exa}[Scattering Poisson, \cite{Lanius16}]\label{exa:scatteringpoissonlift} Let $(X,Z)$ be a log pair. Then $\pi \in {\rm Poiss}(X)$ is a \emph{scattering Poisson structure} if its associated Poisson bracket satisfies
	\be
	\{\cdot,\cdot\}_\pi\colon C^\infty(X) \times C^\infty(X) \to I_Z, \qquad \text{and} 
	\qquad \{\cdot,\cdot\}_\pi\colon I_Z \times C^\infty(X) \to I_Z^2.
	\ee
	Then $\pi$ is of nondegenerate $\mc{C}_Z$-type, and also of $\mc{A}_Z$-type (because $\mc{C}_Z$ is an $\mc{A}_Z$-Lie algebroid).
\end{exa}
\begin{rem} A Poisson structure $\pi \in {\rm Poiss}(X)$ is of $\mc{C}_Z$-type if $Z$ is $\pi$-Poisson, and the first jet of its $\mc{A}_Z$-lift $\pi_{\mc{A}_Z}$ vanishes at $Z$ (using \autoref{prop:azliftpoisson}). Moreover, $\pi$ is of nondegenerate $\mc{C}_Z$-type if and only if $Z = X_{\pi_{\mc{A}_Z}, 2n-2}$ as $\pi_{\mc{A}_Z}$-Poisson divisors (c.f.\ \cite[Lemma 5.2]{Lanius16}), which provides an alternative definition of $\pi$ being scattering Poisson for the log pair $(X,Z)$.
\end{rem}
\begin{exa}[Elliptic Poisson, \cite{CavalcantiGualtieri18}]\label{exa:ellpoissonlift} Let $\pi \in {\rm Poiss}(X)$ be an \emph{elliptic Poisson structure}, i.e.\ such that it is of elliptic divisor-type. Then $I_\pi = I_{|D|}$ for an elliptic pair $(X,|D|)$. The (primary) ideal Lie algebroid of $I_{|D|}$ is the elliptic tangent bundle $\mc{A}_{|D|}$ (see \autoref{exa:ellvx}). By \autoref{thm:ipoissonaitype}, the bivector $\pi$ is of nondegenerate $\mc{A}_{|D|}$-type because $I_{|D|}$ is standard. This recovers \cite[Lemma 3.4]{CavalcantiGualtieri18}. More discussion of these structures can be found in \cite{KlaasseLanius18, KlaasseLi18}, determining their Darboux models, Poisson cohomology and adjoint symplectic groupoids.
\end{exa}
\begin{exa}[Elliptic-log Poisson]\label{exa:elllogpoissonlift} Let $\pi \in {\rm Poiss}(X)$ be an \emph{elliptic-log Poisson structure}, i.e.\ such that it is of elliptic-log divisor-type. Then $I_\pi = I_W$ for an elliptic-log ideal $I_W = I_Z \cdot I_{|D|}$. The degeneracy locus of $\pi$ is then given by $Z$, with $D$ being a $\pi$-symplectic leaf. Consequently, by \autoref{prop:azliftpoisson} we know that $\pi$ is of $\mc{A}_Z$-type. Naturally one wonders whether $\pi$ is in fact of $\mc{A}_W$-type (noting that $\mc{A}_W = TX_{I_W}$ is an $\mc{A}_Z$-Lie algebroid, c.f.\ \autoref{exa:ideallalgebroids}. However, this cannot be immediately answered by \autoref{thm:ipoissonaitype} because $\mc{A}_W^*$ does not admit local closed generators. However, the first part of the lifting procedure (the map $\wt{\pi}^\sharp$) does not use this, so that certainly $\pi$ lifts to a map $\wt{\pi}^\sharp\colon T^*X \to \mc{A}_W$. As $I_W$ is standard, we further know that if $\pi$ is of $\mc{A}_W$-type, then it is of nondegenerate $\mc{A}_W$-type. These types of Poisson structures are studied in more detail in \cite{KlaasseLanius18}, where in particular their Poisson cohomology is determined.
\end{exa}
It is very interesting to extend the above list of examples using a new class of divisor ideals.
\subsection{Lifting Poisson structures of $m$-divisor-type}
\label{sec:liftingalmostreg}
We next turn to lifting Poisson structures of $m$-$I$-divisor-type. To this end, let $\pi \in {\rm Poiss}(X)$ be of $m$-divisor-type. Then $\pi$ is almost-regular by \autoref{cor:divtypealmostreg}, so that by \autoref{prop:almostregulardistr} it has a unique associated involutive regular distribution $D_\pi \subseteq TX$ of $\pi$-Poisson submanifolds. We can readily view this distribution as an injective Lie algebroid $\mc{D}_\pi$, and can wonder whether $\pi$ can be lifted to it.
\begin{prop}\label{prop:mdivtypelift} Let $\pi \in {\rm Poiss}(X)$ be of $m$-divisor-type with its associated injective Lie algebroid $\mc{D}_\pi \to X$. Then $\pi$ is of $I_\pi$-$\mc{D}_\pi$-type, where $I_\pi$ is the divisor ideal of $\pi$.
\end{prop}
\bp As discussed in \autoref{sec:almostregularpoisson}, the regular distribution $D_\pi$ obtained from $\pi$ consists of $\pi$-Poisson submanifolds. \autoref{prop:regularpipoissonlift} thus immediately implies that $\pi$ is of $\mc{D}_\pi$-type. If we consider the definition of $\pi$, we note that its associated divisor $(K, \wedge^m \pi)$ for $K \subseteq \wedge^{2m} TX$ is such that $\det D_\pi = K$. From this we see that its $\mc{D}_\pi$-lift is of $I_\pi$-divisor-type as desired.
\ep
This result is noticed without the above precise language in \cite[Section 4.2]{AndroulidakisZambon17}. Note that the $\mc{D}_\pi$-lift $\pi_{\mc{D}} \in {\rm Poiss}(\mc{D}_\pi)$ is automatically generically-nondegenerate by \autoref{prop:apoissondivtypelifts}.
\begin{rem} \autoref{prop:mdivtypelift} allows for an alternative definition of almost-regularity (or of being of $m$-divisor-type). Namely, an almost-regular Poisson structure is uniquely determined by a regular involutive distribution $D$, and a generically-nondegenerate $\mc{D}$-Poisson structure (with $\mc{D}$ the injective Lie algebroid of $D$). Indeed, if $\pi_{\mc{D}} \in {\rm Poiss}(\mc{D})$ is generically-nondegenerate, then $\pi := \rho_D(\pi_{\mc{D}})$ is almost-regular. \autoref{prop:mdivtypelift} shows the converse.
\end{rem}
The above remark shows a useful viewpoint on almost-regular Poisson structures, or those of $m$-divisor-type. Namely, they are just Poisson structures of divisor-type, except using the injective Lie algebroid $\mc{D}$ instead of $TX$. We further remark here that a Poisson structure being almost-regular is different from saying that it is generically regular ($X_{\pi,{\rm reg}} = X_{\pi,{\rm max}}$), as its degeneracy must further be governed by a divisor (instead of more involved degeneration).

Given an almost-regular Poisson structure $\pi$ with its $\mc{D}$-lift $\pi_{D}$, in light of \autoref{prop:mdivtypelift} and the discussion in \autoref{sec:liftingdivtype}, one can try to lift the $\mc{D}$-Poisson structure $\pi_{\mc{D}}$ further, to the ideal Lie algebroid $\mc{D}_{I_\pi} \to \mc{D} \to TX$. We next consider this for log divisors.
\subsubsection{The case of $m$-log-Poisson structures}
Let $\pi \in {\rm Poiss}(X)$ be an \emph{$m$-log-Poisson structure}, i.e.\ such that it is of $m$-log-divisor-type. Such structures are called \emph{log-f Poisson structure} in \cite{AndroulidakisZambon17}, but we favor our name due to its consistency with other types of almost-regular Poisson structures. Our discussion here is partially parallel to \cite[Section 3.4]{AndroulidakisZambon17}. By definition $I_\pi$ is a log divisor ideal, and the degeneracy locus $Z_\pi$ of $\pi$ is a hypersurface. For ease of notation we set $Z := Z_\pi$ so that $I_\pi = I_Z$, and moreover let $\mc{D} := \mc{D}_\pi$ be the associated injective Lie algebroid.

By \autoref{prop:mdivtypelift} we know that $\pi$ is of log-$\mc{D}$-type. A natural question to ask is whether each leaf of $D$ is log-Poisson. However, this need not necessarily be the case.
\begin{prop}[{\cite[Lemma 3.7]{AndroulidakisZambon17}}]\label{prop:mlogpoissonleaflogpoisson} Let $(X,Z,\mc{D},\pi)$ be an $m$-log-Poisson manifold. Then a leaf $P$ of $D$ is log-Poisson if and only if for all $x \in P \cap Z$, we have that $D_x + T_x Z = T_x X$.
\end{prop}
To find a Lie algebroid to which $\pi$ lifts nondegenerately, we can try to form the primary ideal Lie algebroid $\mc{D}_{I_Z}$, which would be the $(TZ \cap TD,Z)$-rescaling of $\mc{D}$. However, it is not clear a priori that the involutive submodule $\Gamma(\mc{D})_{I_Z} \subseteq \Gamma(\mc{D})$ is projective. For this, we have:
\begin{prop}\label{prop:injectivelogalgebroid} Let $(X,Z)$ be a log pair and $\mc{B} \to X$ be an injective Lie algebroid. Then the $C^\infty(X)$-module $\Gamma(\mc{B})_{I_Z}$ is projective if $Z$ is transverse to the orbits of $\mc{B}$.
\end{prop}
\bp The orbits of $\mc{B}$ are the leaves of the associated distribution $T\mc{B} := \rho_{\mc{B}}(\mc{B}) \subseteq TX$. If $Z$ is transverse to this regular foliation, we can find local bases of generating vector fields $\partial_{x_i}$ for $\mc{B}$ and $\partial_{y_j}$ for $TZ$ so that $\langle \partial_{x_1}, \partial_{y_j} \rangle = \Gamma(TX)$. Then if locally $I_Z = \langle z \rangle$, we have that $\mc{B}_{I_Z} = \langle z \partial_{z}, \partial_{x_i} \, : \, i \geq 2 \rangle$, with $x_1 = z$, which is thus projective (note that $\Gamma(TX)_{I_Z} = \langle z \partial_{z}, \partial_{y_j} \rangle$).
\ep
\begin{rem}\label{rem:dtangentZ} When $Z$ is everywhere tangent to the orbits of $\mc{B}$, the module $\Gamma(\mc{B})_{I_Z}$ is also projective, but just equals $\Gamma(\mc{B})$. Thus in this case we have that $\mc{B}_{I_Z} = \mc{B}$, i.e.\ $\mc{B}$ is unchanged.
\end{rem}
We can apply \autoref{prop:injectivelogalgebroid} to the case where $\mc{B} = \mc{D}$. This shows that if $Z$ is transverse to the orbits of $\mc{D}$, then $\mc{D}_{I_Z}$ is an almost-injective Lie algebroid for which the natural map $\mc{D}_{I_Z} \to \mc{D}$ is of $I_Z$-divisor-type. Moreover it follows from the discussion in \autoref{sec:liftingdivtype} that the $m$-log-Poisson structure $\pi$ is not only of $\mc{D}$-type, but also of $\mc{D}_{I_Z}$-type. We conclude that when $Z$ is transverse to the orbits of $\mc{D}$, then $\pi$ is of nondegenerate $\mc{D}_{I_Z}$-type by \autoref{prop:apoissondivtypelifts}.

An alternate approach is to consider $\pi$ and whether it lifts to the usual ideal Lie algebroid $TX_{I_Z} = \mc{A}_Z$, the log-tangent bundle. As $Z$ is by definition the degeneracy locus of $\pi$, it is a $\pi$-Poisson submanifold by \autoref{prop:pidegenpipoisson}. From this we conclude immediately using \autoref{prop:azliftpoisson} that $\pi$ is of $\mc{A}_Z$-type with lift $\pi_{\mc{A}_Z}$. In fact, it follows that $\pi_{\mc{A}_Z}$ is $2m$-regular when $Z$ is transverse to the orbits of $\mc{D}$. Similarly, when $Z$ is everywhere tangent to the orbits of $\mc{D}$, then $\pi_{\mc{A}_Z}$ is of $m$-log-divisor-type. By construction $\mc{D}_{I_Z}$ is a $\mc{D}$-Lie algebroid, but from the proof of \autoref{prop:injectivelogalgebroid} we see that $\mc{D}_{I_Z}$ is also an $\mc{A}_Z$-Lie algebroid (also \autoref{rem:dtangentZ}). This summarizes into the following diagram.
\begin{center}
	\begin{tikzpicture}
	\matrix (m) [matrix of math nodes, row sep=2.5em, column sep=2.5em,text height=1.5ex, text depth=0.25ex]
	{	(\mc{D}_{I_Z}, \pi_{\mc{D}_{I_Z}}) & (\mc{D}, \pi_{\mc{D}}) \\ (\mc{A}_Z, \pi_{\mc{A}_Z}) & (TX, \pi) \\};
	\path[-stealth]
	(m-2-1) edge node [above] {$\rho_{\mc{A}_Z}$} (m-2-2)
	(m-1-1) edge node [above] {$\rho_{\mc{A}_Z}$} (m-1-2);
	\draw[right hook-latex]
	(m-1-1) edge (m-2-1)
	(m-1-2) edge (m-2-2);
	\end{tikzpicture}
\end{center}
Thus, an $m$-log-Poisson structure $\pi \in {\rm Poiss}(X)$ with divisor ideal $I_Z$ and distribution $D$:
\bi
	\item Is always of log-$\mc{D}$-type;
	\item Is always of $\mc{A}_Z$-type, with $2m$-regular lift when $Z$ is transverse to the $\mc{D}$ orbits;
	\item Is of nondegenerate $\mc{D}_{I_Z}$-type when $Z$ is transverse to the $\mc{D}$-orbits.
\ei
At this point it is clear that it is important to understand the relation between transversality of $Z$ and the orbits of $\mc{D}$, and Poisson geometric properties of $\pi$. The answer to this is provided by the behavior of the modular foliation $\mc{F}_{\pi,{\rm mod}}$ of $\pi$, which in general satisfies $\mc{F}_\pi \subseteq \mc{F}_{\pi,{\rm mod}}$ (discussed in \autoref{sec:modularfoliation}, see \autoref{exa:modfolmlog} in particular). Following \cite{AndroulidakisZambon17}, we denote by
\be
	Z_{\rm sing} = \{x \in Z \, | \, D_x \subseteq T_x Z\}
\ee
the points on $Z$ where $Z$ does not intersect the orbits of $\mc{D}$ transversely. Due to \autoref{prop:mlogpoissonleaflogpoisson}, these are the points on $Z$ whose associated leaf $P$ of $D$ is not log-Poisson. The two extreme cases are where $Z_{\rm sing} = \emptyset$, so that $Z$ is everywhere transverse to the orbits of $\mc{D}$, and the case $Z_{\rm sing} = Z$, where $\mc{D}$ is tangent to $Z$. The general setting is more complicated.

We next follow the above discussion for the three examples contained in \cite[Example 3.8]{AndroulidakisZambon17}.
\begin{exa}\label{rem:mlogpoissonlocal} Consider the following three $1$-log-Poisson structures on $(X =\R^3,(x,y,z))$:
\bi
	\item $\pi_1 = x \partial_x \wedge \partial_y$;
	\item $\pi_2 = z \partial_x \wedge \partial_y$;
	\item $\pi_3 = (z - x^2) \partial_x \wedge \partial_y$.
\ei
We then have that $\Gamma(TX) = \langle \partial_x, \partial_y, \partial_z \rangle$ and moreover $\mc{D} := \mc{D}_1 = \mc{D}_2 = \mc{D}_3 = \langle \partial_x, \partial_y \rangle$. To proceed we list some of the relevant objects for each of these (where $I_i := I_{\pi_i}$ for $i = 1,2,3$):
\bi
	\item $I_1 = \langle x \rangle$, $\Gamma(TZ_1) = \langle \partial_y, \partial_z \rangle$, $Z_{1,{\rm sing}} = \emptyset$, $\Gamma(\mc{A}_{Z_1}) = \langle x \partial_x, \partial_y, \partial_z \rangle$, $\mc{D}_{I_1} = \langle x \partial_x, \partial_y \rangle$;
	\item $I_2 = \langle z \rangle$, $\Gamma(TZ_2) = \langle \partial_x, \partial_y \rangle$, $Z_{2,{\rm sing}} = Z_2$, $\Gamma(\mc{A}_{Z_2})  = \langle \partial_x, \partial_y, z \partial_z \rangle$, $\mc{D}_{I_2} = \langle \partial_x, \partial_y \rangle = \mc{D}$.
	\item $I_3 = \langle z - x^2 \rangle$, $\Gamma(TZ_3) = \langle \partial_x + 2x \partial_z, \partial_y \rangle$, $Z_{3,{\rm sing}} = \{x_1 = 0\} \cap Z$.
	\ei
For $\pi_1$ and $\pi_2$, $Z$ is either everywhere transverse, or everywhere tangent to the orbits of $\mc{D}_i$. The case of $\pi_3$ is more involved, as is the behavior of the modules $\Gamma(TX)_{I_{Z_3}}$ and $\Gamma(\mc{D})_{I_{Z_3}}$.

This shows that $\pi_i$ lifts to be of log-type to $\mc{D}$ for $i = 1,2,3$, as in \autoref{prop:mdivtypelift}. Moreover, the $\mc{A}_{Z_i}$-lifts of the Poisson structures $\pi_i$ for $i = 1,2$ are given by
\be
	\pi_{\mc{A}_{Z_1}} = (x \partial_x) \wedge \partial_y, \quad \pi_{\mc{A}_{Z_2}} = z \partial_x \wedge \partial_y.
\ee
Thus $\pi_{\mc{A}_{Z_1}}$ is $2$-regular, $\pi_{\mc{A}_{Z_2}}$ is of $1$-log-divisor-type, and $\pi_1$ is of nondegenerate $\mc{D}_{I_1}$-type.
\end{exa}
\subsection{Symplectic Lie algebroids}
\label{sec:symplecticliealgebroids}
In this section we briefly discuss symplectic Lie algebroids \cite{NestTsygan01}, the associated study of $\mc{A}$-symplectic geometry, and its relation with $\mc{A}$-Poisson geometry.
\begin{defn}[\cite{NestTsygan01}] A \emph{symplectic Lie algebroid} is a Lie algebroid $\mc{A} \to X$ equipped with an \emph{$\mc{A}$-symplectic structure}, i.e.\ a closed nondegenerate $\mc{A}$-two-form ${\omega_{\mc{A}} \in \Omega^2_{\mc{A}}(X)}$.
\end{defn}
In other words, the $\mc{A}$-two-form $\omega_\mc{A}$ satisfies $d_\mc{A} \omega_\mc{A} = 0$ and $\omega_\mc{A}^n \neq 0$, where ${\rm rank}(\mc{A}) = 2n$. The nondegeneracy condition for $\omega_\mc{A}$ is equivalent to the map $\omega_\mc{A}^\flat\colon \mc{A}^* \to \mc{A}$, given by $v \mapsto \iota_v \omega_\mc{A}$ for $v \in \mc{A}$, being an isomorphism. The space of $\mc{A}$-symplectic forms is denoted by ${\rm Symp}(\mc{A})$, and we say that $X$ is \emph{$\mc{A}$-symplectic} if it is equipped with a symplectic Lie algebroid $(\mc{A},\omega_{\mc{A}}) \to X$.

The nondegeneracy condition for $\omega_{\mc{A}}$ implies that the rank of $\mc{A}$ must be even (but note that $\dim(X)$ need not necessarily be even). Any symplectic Lie algebroid $(\mc{A},\omega_\mc{A})$ defines an $\mc{A}$-cohomology class $[\omega_\mc{A}] \in H^2(\mc{A})$, is naturally oriented by $\omega_\mc{A}^n \in \Gamma(\det(\mc{A}^*))$, and always admits an \emph{$\mc{A}$-almost-complex structure}, i.e.\ a complex structure $J_\mc{A} \in {\rm End}(\mc{A})$ for $\mc{A}$ (see \cite{Klaasse18three}). Lie algebroids $\mc{A}$ of rank two admit $\mc{A}$-symplectic structures if and only if they are orientable, i.e.\ if $w_1(\mc{A}) = 0$.
There is a standard notion of morphism between $\mc{A}$-symplectic manifolds.
\begin{defn} Let $(\varphi,f)\colon (X, \mc{A}, \omega_\mc{A}) \to (X', \mc{A}', \omega_{\mc{A}'})$ be a Lie algebroid morphism between symplectic Lie algebroids. Then the map $(\varphi,f)$ is
\bi
	\item \emph{$\mc{A}$-symplectic} if $\varphi^* \omega_{\mc{A}'} = \omega_\mc{A}$;
	\item an \emph{$\mc{A}$-symplectomorphism} if it $\mc{A}$-symplectic, and $\varphi$ is a fiberwise isomorphism.
\ei
\end{defn}
The $\mc{A}$-symplectic condition $\varphi^* \omega_{\mc{A}'} = \omega_\mc{A}$ is equivalent to demanding $\omega_\mc{A}^\flat = \varphi \circ \omega_\mc{B}^\flat \circ \varphi^*$ as maps. As both $\omega_\mc{A}$ and $\omega_{\mc{A}'}$ are nondegenerate, any $\mc{A}$-symplectic map $(\varphi,f)$ must in particular be fiberwise injective. Consequently, if ${\rm rank}(\mc{A}) = {\rm rank}(\mc{A}')$, any $\mc{A}$-symplectic map is necessarily an $\mc{A}$-symplectomorphism, but the base map $f\colon X \to X'$ is not necessarily a diffeomorphism. Our main interest however will be in Lie algebroid isomorphisms which are $\mc{A}$-symplectic.
\begin{rem} If $\mc{A} \to X$ is a Lie algebroid of divisor-type, one can perform any symplectic operation to a given $\mc{A}$-symplectic structure in the isomorphism locus $X_\mc{A}$. These include for example the symplectic fiber sum \cite{Gompf95} or symplectic blow-up \cite{McDuffSalamon17} procedures.
\end{rem}
Let $(\mc{A},\omega_{\mc{A}}) \to X$ be a symplectic Lie algebroid. Then the map $\omega_{\mc{A}}^\flat\colon \mc{A} \to \mc{A}^*$ can be inverted to a map $(\omega_{\mc{A}}^\flat)^{-1}\colon \mc{A}^* \to \mc{A}$, which we will denote by $\pi_{\mathcal{A}}^\flat$. This construction defines an $\mc{A}$-Poisson structure $\pi_{\mathcal{A}} \in {\rm Poiss}(\mc{A})$, which we will also denote by $\pi_{\mathcal{A}} = \omega_{\mc{A}}^{-1}$. Indeed, the conditions $d\omega_\mc{A} = 0$ and $[\pi_{\mathcal{A}},\pi_{\mathcal{A}}]_\mc{A} = 0$ are equivalent. As discussed below \autoref{defn:apoissonmap}, the two-form $\omega_{\mc{A}}$ is called the \emph{dual $\mc{A}$-form} to the underlying Poisson structure $\pi := \rho_{A}(\pi_{\mc{A}})$.

The following is the Moser theorem for Lie algebroids $\mc{A}$ of smooth divisor-type \cite{KlaasseLanius18}. It covers several results in the literature (\cite{GuilleminMirandaPires14,Lanius16,MarcutOsornoTorres14,MirandaPlanas18,Moser65,NestTsygan96,Radko02,Scott16}), and is similar to \cite{MirandaScott18}.
\begin{thm}[\cite{KlaasseLanius18}]\label{thm:amoser} Let $\mc{A} \to X$ be a Lie algebroid of smooth divisor-type, and let $k \in \{1, 2, \dim X\}$. Let $\omega, \omega' \in \Omega^k(\mc{A})$ be $d_\mc{A}$-closed $\mc{A}$-$k$-forms that are nondegenerate on $Z_\mc{A}$. Assume that either
	\bi
	\item $[\omega] = [\omega'] \in H^k(\mc{A})$, or
	\item $\wt{\omega} - \omega = \rho_{A}^* \tau'$ for $\tau' \in \Omega^k(X)$ satisfying $\tau'|_{Z_\mc{A}} = 0$.
	\ei
	Then there exists a Lie algebroid isomorphism $(\varphi,f)\colon (\mc{A}|_{U},\omega) \to (\mc{A}|_{U'},\omega')$ on neighbourhoods $U$ and $U'$ of $Z_\mc{A}$ for which $\varphi^* \omega' = \omega$, which in the second case can be chosen such that $f|_{Z_\mc{A}} = {\rm id}$.
\end{thm}
\begin{rem} As a consequence of \autoref{thm:amoser}, to establish an $\mc{A}$-Darboux theorem providing a pointwise normal form for $\mc{A}$-symplectic structures, one need only establish what an $\mc{A}$-symplectic structure must look like locally at a point in $Z_\mc{A}$.
\end{rem}
This result implies that $\mc{A}$-Nambu structures of top degree (nonvanishing sections of $\det(\mc{A})$), and hence $\mc{A}$-symplectic structures on surfaces, specifying the same $\mc{A}$-orientation are classified by their $\mc{A}$-cohomology class. Namely, letting $n = \dim X$, we note that $\det(\mc{A}^*)$ is a line bundle. Consequently, given cohomologous forms $\omega,\omega' \in \Omega^n(\mc{A}) = \Gamma(\det(\mc{A}^*))$ we have $\omega = f \omega'$ for some nonvanishing function $f \in C^\infty(X)$. This function must be strictly positive as $\omega$ and $\omega'$ give rise to the same $\mc{A}$-orientation, so that $(1-t)\omega + t \omega' = ((1-t) + t f) \omega$ is nondegenerate for all $t \in [0,1]$. With more knowledge about the Lie algebroid $\mc{A}$, the assumption on induced $\mc{A}$-orientations can sometimes be dropped. This discussion summarizes as follows.
\begin{prop} Let $\mc{A} \to X$ be a Lie algebroid of divisor-type for which $Z_\mc{A}$ is smooth. An $\mc{A}$-Nambu structure $\Pi$ inducing a given $\mc{A}$-orientation is classified up to $\mc{A}$-orientation-preserving isomorphism by its $\mc{A}$-cohomology class.
\end{prop}
See also \cite{MartinezTorres04} and \cite{MirandaPlanas18} for the cases of log- and $b^k$-Nambu structures of top degree. Let us reiterate this classifies such $\mc{A}$-symplectic structures on surfaces using $H^2(\mc{A})$. Below we give some examples of local Darboux theorems for Poisson structures of divisor-type.
\begin{exa}[Nondegenerate] Let $\pi \in {\rm Poiss}(X^{2n})$ be nondegenerate. Then $\pi$ is dual to a symplectic structure. In other words, it is of nondegenerate $TX$-type. The classical Darboux theorem (e.g.\ proven using Moser methods) states that around any point in $X$ there are coordinates in which $\pi = \omega_0^{-1}$, with $\omega_0 = \sum_i dx_i \wedge dy_i$ the standard symplectic structure in $\R^{2n}$.
\end{exa}
\begin{exa}[Log-Poisson] Let $\pi \in {\rm Poiss}(X)$ be log-Poisson, i.e.\ of log divisor-type with divisor ideal $I_Z$. Then $\pi$ is of nondegenerate $\mc{A}_Z$-type by \autoref{exa:logpoissonlift}. Either directly by applying Weinstein's splitting theorem, or by Moser methods for $\mc{A}_Z$ (see \cite{GuilleminMirandaPires14}), around points in $Z$ there are coordinates $(z,x_i)$ with $I_Z = \langle z \rangle$ in which $\pi = z \partial_z \wedge \partial_{x_2} + \omega_0^{-1}$.
\end{exa}
The normal-crossing log case is more involved, as there can be local cohomology of the Lie algebroid $\mc{A}_{\underline{Z}}$ in contractible neighbourhoods. See the discussions found in \cite{Lanius16two,Lanius17,MirandaScott18,Radko02two}.
\begin{exa}[$b^k$-Poisson]\label{exa:bkpoissondarboux} Let $\pi \in {\rm Poiss}(X)$ be $b^k$-Poisson as in \autoref{exa:bkpoissonlift}. Then $\pi$ is of nondegenerate $\mc{A}_Z^k$-type for the associated $b^k$-tangent bundle. By Moser methods for $\mc{A}_Z^k$ (see \cite{GuilleminMirandaWeitsman17,Scott16}), around points in $Z$ there are coordinates $(z,x_i)$ with $z \in j_{k-1}$ such that
\be
	\pi = z^k \partial_z \wedge \partial_{x_1} + \omega_0^{-1}.
\ee
\end{exa}
\begin{exa}[Scattering Poisson] Let $\pi \in {\rm Poiss}(X)$ be scattering Poisson, so that $\pi$ is of nondegenerate $\mc{C}_Z$-type (as in \autoref{exa:scatteringpoissonlift}). By applying Moser methods for $\mc{C}_Z$ (see \cite{Lanius16}), around points in $Z$ there are coordinates $(z,x_i)$ with $I_Z = \langle z \rangle$ in which (for $\alpha_0$ the standard contact structure, so that $d\alpha_0 = \omega_0$)
\be
	\pi = z^3 \partial_z \wedge \partial_{x_1} + z^2 \partial_{z} \wedge \alpha_0^{-1} + z^2 \omega_0^{-1}.
\ee
\end{exa}
\begin{exa}[Elliptic Poisson]\label{exa:ellpoissondarboux} Let $\pi \in {\rm Poiss}(X)$ be elliptic Poisson, i.e.\ of elliptic divisor-type with divisor ideal $I_{|D|}$. Then $\pi$ is of nondegenerate $\mc{A}_{|D|}$-type by \autoref{exa:ellpoissonlift}. The local (and global) behavior of $\pi$ depends on the value(s) of its elliptic residue ${\rm Res}_q(\pi)$. This is a residue map associated to the elliptic tangent bundle in \cite{CavalcantiGualtieri18}, which we recall in \autoref{exa:elltangentresidues}. There are two qualitatively different local Darboux models for $\pi$, depending on the vanishing of ${\rm Res}_q(\pi)$ (see \cite{KlaasseLanius18}). Around points in $D$ there are coordinates $(r,\theta,x_i)$ with $I_{|D|} = \langle r^2 \rangle$ with:
\be
	\pi = \begin{cases} r \partial_r \wedge \partial_{x_1} + \partial_{\theta} \wedge \partial_{x_2} + \omega_0^{-1} \qquad &\text{if } {\rm Res}_q(\pi) = 0,\\ \lambda r \partial_r \wedge \partial_\theta + \omega_0^{-1} \qquad &\text{if } {\rm Res}_q(\pi) = \lambda \neq 0.
\end{cases}
\ee
This recovers the classification results of such Poisson structures on surfaces which can be found in \cite{Radko02two} (where ${\rm Res}_q(\pi) \neq 0$ is forced by nondegeneracy). See also \cite{KlaasseLi18}.
\end{exa}
The local behavior of elliptic-log Poisson structures (\autoref{exa:elllogpoissonlift}) is described in \cite{KlaasseLanius18}.

%% file: residues.tex
\section{Residue maps and the modular foliation}
\label{sec:residues}
In this section we discuss two distinct but related concepts. The first of these is the residue maps one can associate to a Lie algebroid $\mc{A}$ if one is given an appropriate $\mc{A}$-invariant submanifold (\autoref{defn:ainvariantresidues}). After this we will discuss Lie algebroid modules and consequently the Lie algebroid cohomologies they give rise to. We then discuss when the residue maps, which are initially defined on the levels of Lie algebroid forms, descend to the level of cohomology. After this we proceed to discuss the Lie algebroid modules in the presence of Poisson structures, and in particular consider the modular foliation. We finish by relating the residue maps to symplectic Lie algebroids and the modular foliation in concrete examples.
\subsection{Residue maps}
\label{sec:residuemaps}
In this section we discuss residues of Lie algebroid forms. These are a way to extract relevant information of $\mc{A}$-forms over $\mc{A}$-invariant submanifolds where the restriction of $\mc{A}$ is projective. We further discuss how such residues interact with Lie algebroids morphisms. Residues are used when dealing with geometric structures on $\mc{A}$, in order to extract information along the submanifold (see e.g.\ \cite{CavalcantiGualtieri18, CavalcantiKlaasse18}). In \cite{CavalcantiKlaasse18} the following is discussed for transitive $\mc{A}$-invariant submanifolds, but here we expose a generalization allowing for projective ones, recovering the transitive case when we can take $\mc{B} = TD$ below. As a special case it further covers the case when the restriction of $\mc{A}$ is regular.

Let $\mc{A} \to X$ be a Lie algebroid and let $D \subseteq X$ be a projective $\mc{A}$-invariant submanifold (as in \autoref{defn:lalgebroidadjective}), with $\mc{B} \to D$ and $(\rho_A|_D)(\Gamma(\mc{A}|_D)) \cong \Gamma(\mc{B})$. We obtain a short exact sequence
\be
0 \to \ker \wt{\rho}_\mc{A}|_D \to \mc{A}|_D \to \mc{B} \to 0,
\ee
where the $\wt{\rho}_\mc{A}|_D\colon \mc{A}|_D \to \mc{B}$ is the induced bundle map coming from the map on sections. Note that in the transitive case we have that $\wt{\rho}_{\mc{A}}|_D = \rho_{A}|_D$.
We see that $\mc{A}|_D$ is an extension of $\mc{B}$ by $\ker \wt{\rho}_\mc{A}|_D$, the latter being the germinal isotropy of \autoref{sec:singularliealgebroids}. The case in which we are mainly interested in is when $D = X \backslash X_\mc{A} = Z_\mc{A}$, the complement of the isomorphism locus of $\mc{A}$ (although in general the restriction need not be transitive nor projective). This is an $\mc{A}$-invariant submanifold of $X$ if it is smooth by \autoref{prop:degenlocus}, as $X \backslash X_\mc{A} = X_{\mc{A},n-1}$ where $\dim X = {\rm rank}(\mc{A}) = n$ (see \autoref{exa:denseisolocus}). Note further that $X$ itself is always $\mc{A}$-invariant, as is sometimes useful (say, for regular Lie algebroids). Dualizing the above sequence, we obtain
\be
0 \to \mc{B}^* \to \mc{A}^*|_D \to (\ker \wt{\rho}_\mc{A}|_D)^* \to 0.	
\ee
We are now in the following more general situation: given a short exact sequence $\mc{S}\colon 0 \to E \to W \to V \to 0$ of vector spaces, there is an associated dual sequence $\mc{S}^*\colon 0 \to V^* \to W^* \to E^* \to 0$. For a given $k \in \N$, by taking $k$th exterior powers we obtain a filtration of spaces $\mc{F}^i := \{\rho \in \wedge^k W^* \, | \, \iota_x \rho = 0$ for all $x \in \wedge^i E \}$, for $i = 0,\dots,k+1$. These spaces satisfy $\mc{F}^0 = 0$, $\mc{F}^1 = \wedge^k V^*$, $\mc{F}^i \subset \mc{F}^{i+1}$, and $\mc{F}^{i+1} / \mc{F}^i \cong \wedge^{k-i} V^* \otimes \wedge^i E^*$. Setting $\ell := \dim E$, we have $\mc{F}^{\ell + 1} = \wedge^k W^*$. The \emph{residue} of an element $\rho \in \wedge^k W^*$ is defined as its equivalence class ${\rm Res}(\rho) = [\rho] \in \mc{F}^{\ell + 1} / \mc{F}^{\ell} \cong \wedge^{k-\ell} V^* \otimes \wedge^\ell E^*$. Upon a choice of trivialization of $\wedge^\ell E^*$, i.e.\ a choice of volume element for $E$, one can view the residue ${\rm Res}(\rho)$ as an element of $\wedge^{k-\ell} V^*$.

We now define the top residue of a Lie algebroid along a projective $\mc{A}$-invariant submanifold, by performing the above operation to a smoothly varying section of $\wedge^\bullet \mc{A}^*$ restricted to $D$.
\begin{defn}\label{defn:ainvariantresidues} Let $\mc{A} \to X$ be a Lie algebroid and $D \subseteq X$ a projective $\mc{A}$-invariant submanifold with image $\mc{B}$ and germinal isotropy $E := \ker \wt{\rho}_\mc{A}|_D$. The \emph{residue map} of $\mc{A}$ along $D$ is the map ${\rm Res}_D\colon \Omega^\bullet(\mc{A}) \to \Omega^{\bullet-\ell}(\mc{B}; \wedge^\ell E^*)$, where $\ell = \dim \ker \rho_\mc{A}|_D = {\rm rank}(\mc{A}) - {\rm rank}(\mc{B})$.
\end{defn}
There are also lower residues ${\rm Res}_{-m}\colon \wedge^k W^* \to \mc{F}^{\ell+1} / \mc{F}^{\ell-m}$ for $m > 0$. These are always defined but have a better description for forms $\rho \in \wedge^k W^*$ whose higher residues vanish, so that ${\rm Res}_{-m}(\rho) \in \mc{F}^{\ell-m +1} / \mc{F}^{\ell -m} \cong \wedge^{k-\ell+m} V^* \otimes \wedge^{\ell-m} E^*$. This leads to the following.
\begin{defn} Let $\mc{A} \to X$ be a Lie algebroid and $D \subseteq X$ a projective $\mc{A}$-invariant submanifold with image $\mc{B}$ and germinal isotropy $E := \ker \wt{\rho}_\mc{A}|_D$. The \emph{lower residue maps} for $m \geq 0$ of $\mc{A}$ along $D$ are ${\rm Res}_{D,-m}\colon \Omega^\bullet_{1-m}(\mc{A}) \to \Omega^{\bullet - \ell + m}(\mc{B}; \wedge^{\ell - m} E^*)$, where $\Omega^\bullet_{1-m}(\mc{A})$ is inductively the space of $\mc{A}$-forms all of whose $(1-m)$th or higher residues along $D$ vanish.
\end{defn}
\begin{rem} In case $D = Z_\mc{A}$ for a Lie algebroid $\mc{A}$ with dense isomorphism locus, one should think of these residues as extracting the coefficients in front of the singular parts of an $\mc{A}$-form along $Z_\mc{A}$, as then $(\ker \wt{\rho}_\mc{A})^*|_D$ consists of ``singular'' generators.
\end{rem}
\begin{rem} Residue maps arise out of any extension of Lie algebroids. This extension need not be induced by the anchor map of the Lie algebroid, even though this is typical.
\end{rem}
Given a map of short exact sequences $\Psi\colon \wt{\mc{S}} \to \mc{S}$ with dual map $\Psi^*\colon \mc{S}^* \to \wt{\mc{S}}^*$, there is a corresponding map of filtrations $\Psi^*\colon \mc{F}^i \to \wt{\mc{F}}^i$. Setting $\wt{\ell} := \dim \wt{E}$, we have the following.
\begin{lem}\label{lem:residues} In the above setting, assume that $\wt{\ell} > \ell$. Then $\wt{\rm Res}(\Psi^* \rho) = 0$ for all $\rho \in \wedge^k W^*$.
\end{lem}
\bp We have $\rho \in \mc{F}^{\ell +1}$ so that $\Psi^* \rho \in \wt{\mc{F}}^{\ell + 1}$. As $\wt{\ell} > \ell$, we have $\wt{\mc{F}}^{\ell+1} \subset \wt{\mc{F}}^{\wt{\ell}} \subset \wt{\mc{F}}^{\wt{\ell}+1}$, so that $\wt{\rm Res}(\Psi^* \rho) = [\Psi^* \rho] \in \wt{\mc{F}}^{\wt{\ell}+1} / \wt{\mc{F}}^{\wt{\ell}}$ vanishes by degree reasons as desired.
\ep
Assuming $\wt{\ell} > \ell$, we see from \autoref{lem:residues} that all lower residues automatically vanish by degree reasons until considering $\Psi^* \rho \in \wt{\mc{F}}^{\ell +1}$. Hence the first possibly nonzero residue is $\wt{\rm Res}_{\ell-\wt{\ell}}(\Psi^* \rho) = [\Psi^* \rho] \in \wt{\mc{F}}^{\ell+1} / \wt{\mc{F}}^{\ell} \cong \wedge^{k-\ell} \wt{V}^* \otimes \wedge^\ell \wt{E}^*$. As in \autoref{lem:residues} we obtain the following.
\begin{lem}\label{lem:residuecommute} In the above setting, assuming $\wt{\ell} \geq \ell$, we have $\Psi^* \circ {\rm Res} = \wt{\rm Res}_{\ell - \wt{\ell}} \circ \Psi^*$.
\end{lem}
Out of the above discussion the following result on residue maps is immediate.
\begin{prop}\label{prop:ainvsubmfdresidue} Let $(\varphi,f)\colon (\mc{A},X) \to (\mc{A}',X')$ be a Lie algebroid morphism and $D \subseteq X^n$, $D' \subseteq X'$ transitive $\mc{A}$- respectively $\mc{A}'$-invariant submanifolds such that $f\colon (X,D) \to (X',D')$ is a strong map of pairs. Define $\ell = {\rm rank}(\mc{A}) - \dim D$ and $\ell' = {\rm rank}(\mc{A}') - \dim D'$, and assume that $\ell' \geq \ell$. Then ${\rm Res}_{\ell - \ell'} \circ \varphi^* = f^* \circ {\rm Res}$. Moreover, ${\rm Res}_{-m} \circ \varphi^* = 0$ for $m > \ell' - \ell$.
\end{prop}
\bp As $f$ is a strong map of pairs, we have $T f\colon TD \to TD'$. Moreover, as noted below \autoref{defn:laisotropy}, $\varphi$ restricts to $\varphi\colon \ker \rho_\mc{A}|_D \to \ker \rho_{\mc{A}'}|_{D'}$. Consequently, $\varphi$ induces a map of the relevant short exact sequences defining the residues. The result follows from \autoref{lem:residuecommute}.
\ep
For Lie algebroids of divisor-type this specifies to the following.
\begin{cor}\label{cor:residuemapsdenseiso} Let $(\varphi,f)\colon (\mc{A},X) \to (\mc{A}',X')$ be a Lie algebroid morphism between Lie algebroids of divisor-type with transitive degeneracy loci. Assume that $f\colon (X,Z_\mc{A}) \to (X',Z_{\mc{A}'})$ is a strong map of pairs and that ${\rm codim}\, D' \geq {\rm codim}\, D$. Then ${\rm Res}_{\ell-\ell'} \circ \varphi^* = f^* \circ {\rm Res}$ and ${\rm Res}_{-m} \circ \varphi^* = 0$ for $m > \ell' - \ell$, where $\ell = {\rm codim}\, D$ and $\ell' = {\rm codim}\, D'$.
\end{cor}
\bp This follows from \autoref{prop:ainvsubmfdresidue}, where we use that degeneracy loci are $\mc{A}$-invariant due to \autoref{prop:degenlocus}, and that ${\rm rank}(\mc{A}) = \dim X$ because $X_\mc{A}$ is nonempty.
\ep
\begin{rem} When we speak of projective or transitive degeneracy loci, it is implied that we assume that each degeneracy locus has components being smooth submanifolds.
\end{rem}
\subsection{Lie algebroid modules and cohomology}
\label{sec:aconnareps}
We next discuss Lie algebroid connections and Lie algebroid modules, which are vector bundles equipped with a flat Lie algebroid connection. We further introduce the associated Lie algebroid cohomologies and characteristic classes, including the modular class. For more information, see e.g.\ \cite[Chapter 7]{Mackenzie05} and also \cite{EvensLuWeinstein99, KosmannSchwarzbachLaurentGengouxWeinstein08,GualtieriPym13}. Let $\mc{A} \to X$ be a Lie algebroid and $E \to X$ a vector bundle.
\begin{defn} An \emph{$\mc{A}$-connection} on $E$ is a bilinear map $\nabla\colon \Gamma(\mc{A}) \times \Gamma(E) \to \Gamma(E)$, $(v,\sigma) \mapsto \nabla_v \sigma$, such that $\nabla_{f v} \sigma = f \nabla_v \sigma$ and $\nabla_v(f\sigma) = f \nabla_v \sigma + (\rho_\mc{A}(v) f) \cdot \sigma$ for all $v \in \Gamma(\mc{A})$, $\sigma \in \Gamma(E)$ and $f \in C^\infty(X)$. We will also consider $\nabla$ as a map $\nabla\colon \Gamma(E) \to \Omega^1(\mc{A};E)$ via $(\nabla s)(v) := \nabla_v s$.
\end{defn}
Any $\mc{A}$-connection $\nabla$ on $E$ has a \emph{curvature tensor} $F_\nabla \in \Gamma(\wedge^2 \mc{A}^* \otimes {\rm End}(E))$, given by $F_\nabla(v,w) := \nabla_v \circ \nabla_w - \nabla_w \circ \nabla_v - \nabla_{[v,w]_\mc{A}}$ for $v, w \in \Gamma(\mc{A})$. An $\mc{A}$-connection $\nabla$ is \emph{flat} if $F_\nabla \equiv 0$.
\begin{defn}\label{defn:amodule} An \emph{$\mc{A}$-module} is a vector bundle $E$ equipped with a flat $\mc{A}$-connection $\nabla$.
\end{defn}
Given an $\mc{A}$-module $(E,\nabla)$, there is a differential $d_{\mc{A},\nabla}$ on $\Omega^\bullet(\mc{A};E)$, the space of $\mc{A}$-forms with values in $E$. Recalling that $\Omega^\bullet(\mc{A};E) = \Gamma(\wedge^\bullet \mc{A}^*) \otimes \Gamma(E)$, we set $d_{\mc{A},\nabla}(\eta \otimes s) := d_{\mc{A}} \eta \otimes s + (-1)^{|\eta|} \eta \otimes \nabla s$, with $|\eta|$ the degree of $\eta$. This differential satisfies $d_{\mc{A},\nabla}(\eta \wedge \xi \otimes s) = d_{\mc{A}} \eta \wedge \xi \otimes s + (-1)^{|\eta|} \eta \wedge dt_{\mc{A},\nabla} (\xi \otimes s)$, and $d_{\mc{A},\nabla}$ squares to zero if and only if $\nabla$ is flat. In fact, there is a bijective correspondence between such operators on $\Omega^\bullet(\mc{A};E)$, and flat $\mc{A}$-connections on $E$.
\begin{defn} Let $\mc{A} \to X$ be a Lie algebroid and $(E,\nabla)$ an $\mc{A}$-module. The \emph{$\mc{A}$-cohomology with values in $E$} is given by $H^k(\mc{A}; E) := H^k(\Omega^\bullet(\mc{A};E), d_{\mc{A},\nabla})$ for $k \in \N \cup \{0\}$.
\end{defn}
Note that $\mc{A}$-modules are also called \emph{$\mc{A}$-representations}, and that they are Lie algebra representations when $X = \{{\rm pt}\}$. Denote the space of all $\mc{A}$-representations by $\mc{A}$-${\rm Rep}(X)$.
\begin{exa} Let $\underline{\R} \to X$ be the trivial line bundle, carrying a trivial $\mc{A}$-representation structure given by $\nabla_v f := \mc{L}_{\rho_{\mc{A}}(v)} f$ for $v \in \Gamma(\mc{A})$ and $f \in \Gamma(\underline{\R}) = C^\infty(X)$. We then have that $H^\bullet(\mc{A};\underline{\R}) =: H^\bullet(\mc{A})$, obtaining what is usually called the \emph{Lie algebroid cohomology} of $\mc{A}$.
\end{exa}
Lie algebroid cohomology is generally hard to compute and need not be finite-dimensional. However, the cohomology of most of the Lie algebroids of divisor-type can be determined (see \cite{KlaasseLanius18} for an overview). Any Lie algebroid $\mc{A} \to X$ always has a canonical $\mc{A}$-module (see \cite{EvensLuWeinstein99}).
\begin{defn}\label{defn:canamodule} Let $\mc{A} \to X$ be a Lie algebroid and set $Q_\mc{A} := \det(\mc{A}) \otimes \det(T^*X)$. Then $Q_\mc{A}$ is  the \emph{canonical $\mc{A}$-module}, using the $\mc{A}$-connection on $Q_\mc{A}$ that is defined for $v \in \Gamma(\mc{A})$, $V \in \Gamma(\det(\mc{A}))$ and $\mu \in \Gamma(\det(T^*X))$ by the formula $\nabla_v(V \otimes \mu) := \mc{L}_{v} V \otimes \mu + V \otimes \mc{L}_{\rho_\mc{A}(v)} \mu$.
\end{defn}
One readily checks this formula indeed defines a flat $\mc{A}$-connection on the line bundle $Q_\mc{A}$. Several operations can be performed on the class of Lie algebroid representations, such as pullbacks along Lie algebroid morphisms, duals, tensor products, and exterior products.
\begin{rem} For Lie algebroids $\mc{A}^n \to X$ of divisor-type we have $Q_\mc{A}^* \cong {\rm div}(\mc{A})$. There should be a relation between this and the Evens--Lu--Weinstein pairing $H^k(\mc{A}) \otimes H^{n-k}(\mc{A};Q_\mc{A}) \to \R$ (see \cite{EvensLuWeinstein99}), especially when $\mc{A} = TX_I$ for a standard projective divisor ideal $I \subseteq C^\infty(X)$.
\end{rem}
Any $\mc{A}$-module $(L,\nabla)$ of rank $1$ induces a cohomology class in $H^1(\mc{A})$. Given a nonvanishing section $s \in \Gamma(L^2)$ and $v \in \Gamma(\mc{A})$ (noting that $L^2 = L \otimes L$ is always trivial) we can write $\langle v, d_{\mc{A},\nabla} s\rangle = \theta_s(v) s$ for some element $\theta_s \in \Gamma(\mc{A}^*)$. As $d_{\mc{A},\nabla}$ squares to zero, the section $\theta_s$ is $d_\mc{A}$-closed, and one verifies that the cohomology class $\theta_L := \frac12 [\theta_s] \in H^1(\mc{A})$ is independent of $s$. This is the \emph{characteristic class} of $(L,\nabla)$. Characteristic classes are natural with respect to tensor products and pullbacks of modules, i.e.\ they satisfy $\theta_{L \otimes L'} = \theta_L + \theta_{L'}$ and $\varphi^* \theta_L = \theta_{\varphi^*L} \in H^1(\mc{A}')$ given $\mc{A}$-modules $L$ and $L'$ of rank $1$ and a Lie algebroid morphism $\varphi\colon \mc{A}' \to \mc{A}$.

The characteristic class of the trivial $\mc{A}$-module is zero, and the characteristic class of $Q_\mc{A}$ is called the \emph{modular class} ${\rm mod}(\mc{A})$ of $\mc{A}$. One says that $\mc{A}$ is \emph{unimodular} if ${\rm mod}(\mc{A}) = 0$. Similarly an $\mc{A}$-module $(L,\nabla)$ of rank $1$ is \emph{umimodular} if its characteristic class vanishes. In this case it admits a global nonvanishing section which is $d_{\mc{A},\nabla}$-closed.
\subsubsection{Representations up to homotopy}
In studying the residue map of \autoref{sec:residuemaps} we will have use for a powerful extension of the notion of an $\mc{A}$-representation, introduced in \cite{AriasAbadCrainic12}. To start, recall that a \emph{$\Z$-graded vector bundle} $\mathbb{E} \to X$ is a direct sum $\mathbb{E} = \bigoplus_{i \in \Z} E_i$ of vector bundles, where $E_i$ has \emph{degree} $i$. Further, ${\rm End}_k(\mathbb{E}) \to X$ is the associated bundle of endomorphisms of $\mathbb{E}$ of degree $k$. Let $\mc{A} \to X$ be a Lie algebroid. Then $\Omega^\bullet(\mc{A};\mathbb{E}) := \Gamma(\wedge^\bullet \mc{A}^* \otimes \mathbb{E})$ becomes a right graded $\Omega^\bullet(\mc{A})$-module in the standard way. With this we have:
\begin{defn}[\cite{AriasAbadCrainic12}] A $\Z$-graded vector bundle $\mathbb{E} \to X$ is an \emph{$\mc{A}$-representation up to homotopy} if $\Omega^\bullet(\mc{A};\mathbb{E})$ carries a differential $D$ such that it is a right graded $(\Omega^\bullet(\mc{A}),d_\mc{A})$-module.
\end{defn}
As for $\mc{A}$-representations, these lead to cohomology theories $H^\bullet(\mc{A};\mathbb{E})$. Moreover, there is a similar bijective correspondence between such differentials $D$, and $\mc{A}$-connections on $\mathbb{E}$ coming with maps $\omega_i\colon \Gamma(\wedge^i \mc{A}) \to {\rm End}_{1-i}(\mathbb{E})$ for $i \geq 2$ satisfying certain coherence conditions.

Denote by $\mc{A}$-${\rm RepHtpy}$ the space of $\mc{A}$-representations up to homotopy, with $\mc{A}$-${\rm RepHtpy}_{k,k'}$ for $k,k' \in \Z$ the subclasses of those for which the associated $\Z$-graded vector bundle $\mathbb{E}$ satisfies $E_i = \{0\}$ for all $i$ not satisfying $k \leq i < k'$. Finally, set $\mc{A}$-${\rm RepHtpy}_{k,k} := \mc{A}$-${\rm RepHtpy}_{k,k}$ for those concentrated in degree $k$. Note then that there is a natural inclusion $\mc{A}$-${\rm Rep} \hookrightarrow \mc{A}$-${\rm RepHtpy}_0$ of $E \mapsto \mathbb{E}$ with $\mathbb{E}_0 = E$. The class of $\mc{A}$-representations up to homotopy admits duals, tensor products, and exterior products, similar to what is true for $\mc{A}$-representations.

We now explain that $\mc{A}$-representations up to homotopy give rise to $\mc{A}$-representations. Given a $\Z$-graded vector bundle $\mathbb{E} \to X$, let $\mathbb{E}_{\rm red}$ be the resulting $\Z_2$-graded vector bundle obtained from the map $\Z \to \Z_2$ of reduction modulo two. Given such a $\Z_2$-graded vector bundle $\mathbb{F} \to X$ with summands $\mathbb{F} = F_0 \oplus F_1$, we can consider its \emph{Berizinian}
\be
	{\rm Ber}(F) := \det(F_0) \otimes \det(F_1^*) \to X.
\ee
With this, given a $\Z$-graded vector bundle $\mathbb{E} \to X$, define ${\rm Ber}(\mathbb{E}) := {\rm Ber}(\mathbb{E}_{\rm red}) \to X$.
\begin{prop}[\cite{Mehta15}]\label{prop:berrephtpy} Let $\mathbb{E} \in \mc{A}$-${\rm RepHtpy}$. Then assigning $\mathbb{E} \mapsto {\rm Ber}(\mathbb{E})$ defines a map
\be
	{\rm Ber}\colon \mc{A}\text{-}{\rm RepHtpy} \to \mc{A}\text{-}{\rm Rep}.
\ee
\end{prop}
\begin{rem}\label{rem:adjointrep} We noted in \autoref{defn:canamodule} that there is a canonical modular $\mc{A}$-module $Q_\mc{A} = \det(\mc{A}) \otimes \det(T^*X)$. We can understand this better using the above framework. Namely, every Lie algebroid $\mc{A} \to X$ has a canonical $\mc{A}$-representation up to homotopy, namely its \emph{adjoint representation (up to homotopy)} ${\rm Adj}(\mc{A})$, given by the anchor sequence $\mc{A} \to TX$ with its natural differential $D_\nabla$ \cite{AriasAbadCrainic12}. Note now that ${\rm Ber}({\rm Adj}(\mc{A})) = Q_\mc{A}$, and one readily checks that the resulting induced $\mc{A}$-module structure on $Q_\mc{A}$ is the one of \autoref{defn:canamodule} (see \cite{Mehta15}).
\end{rem}
\begin{rem} The cohomology groups $H^\bullet(\mc{A}; {\rm Adj}(\mc{A}))$ of the adjoint representation of \autoref{rem:adjointrep} form the \emph{deformation cohomology} $H^\bullet_{\rm def}(\mc{A})$ of the Lie algebroid $\mc{A} \to X$ \cite{CrainicMoerdijk08}. This important fact motivated the very definition of an $\mc{A}$-representation up to homotopy (see \cite{AriasAbadCrainic12}).
\end{rem}
There is a class of representations up to homotopy called \emph{Serre representations} that is of main importance to us (\cite[Example 4.15]{AriasAbadCrainic12}). Consider an extension of Lie algebroids over a manifold $X$, i.e.\ a short exact sequence of Lie algebroid morphisms as follows:
\be
	0 \to E \stackrel{i}{\to} \mc{A} \stackrel{p}{\to} \mc{A}' \to 0.
\ee
Let $\sigma\colon \Gamma(\mc{A}') \to \Gamma(\mc{A})$ be an Ehresmann connection splitting the above sequence, i.e.\ such that $p \circ \sigma = {\rm id}_{\mc{A}'}$. This gives a $\sigma$-dependent identification $\Gamma(\mc{A}) \cong \Gamma(\mc{A}') \oplus \Gamma(E)$. Consider the $\Z$-graded complex $(\Gamma(\wedge^\bullet E^*), d_E)$. Given $v \in \Gamma(\mc{A}')$ there is a degree-zero operator $\nabla_v^\sigma := {\rm ad}^*_{\sigma(v)}\colon \Gamma(\wedge^\bullet E^*) \to \Gamma(\wedge^\bullet E^*)$ using the bracket $[\cdot,\cdot]_\mc{A}$, which is thus an $\mc{A}'$-connection $\nabla^\sigma$ on $\Gamma(\wedge^\bullet E^*)$. Further, the \emph{curvature} $R_\sigma \in \Gamma(\wedge^2 \mc{A}'^* \otimes E)$ of $\sigma$ is given for $v,v' \in \Gamma(\mc{A}')$ by
\be
	R^\sigma(v,v') := [\sigma(v),\sigma(v')]_\mc{A} - \sigma([v,v']_{\mc{A}'}).
\ee
This naturally turns into a map $i(R^\sigma)\colon \Gamma(\wedge^2 \mc{A}') \to {\rm End}_{-1}(\wedge^\bullet E^*)$ after using the contraction $i\colon E \to {\rm End}_{-1}(\wedge^\bullet E^*)$. Now define the operator $D := d_E + \nabla^\sigma + i(R^\sigma)$ on $\mathbb{E} := \wedge^\bullet E^*$.
\begin{prop}[{\cite[Example 4.15]{AriasAbadCrainic12}}]\label{prop:rephtpyextension} The triple $(\mathbb{E},\sigma, D)$ defines the structure of an $\mc{A}'$-representation up to homotopy on the exterior algebra $\mathbb{E} = \wedge^\bullet E^*$.
\end{prop}
\begin{rem} When $E$ is abelian, i.e.\ has trivial bracket, the above structure in fact defines an $\mc{A}$-module structure on $E$. More generally this is true for the center and abelianization $Z(E)$ and $E/[E,E]$ of the bundle Lie algebroid $E$ (see \cite[Proposition 3.3.20]{Mackenzie05}).
\end{rem}
\subsubsection{Residue maps and representations}
We now discuss when the residue maps of \autoref{sec:residuemaps} induce maps in cohomology. Let $\mc{A} \to X$ be a Lie algebroid and let $D \subseteq X$ a projective $\mc{A}$-invariant submanifold, so that we obtain an extension of Lie algebroids
\be
	0 \to \ker\wt{\rho}_{\mc{A}}|_D \to \mc{A}|_D \to \mc{B} \to 0.
\ee
By applying \autoref{prop:rephtpyextension} we obtain the structure of a $\mc{B}$-representation up to homotopy on $E := \ker \wt{\rho}_\mc{A}|_D \to D$. Using \autoref{prop:berrephtpy} we obtain a $\mc{B}$-module structure on $\det(E) \to D$, and consequently also on its dual $\det(E^*)$. We can now prove the following (\autoref{thm:introresidues}).
\begin{thm}\label{thm:residues} Let $\mc{A} \to X$ be a Lie algebroid and $(\mc{B},D)$ a projective $\mc{A}$-invariant submanifold with $\ell$-dimensional germinal isotropy $E \to D$. Then ${\rm Res}_D\colon \Omega^\bullet(X,\mc{A}) \to \Omega^{\bullet-\ell}(D, \mc{B}; \det(E^*))$ is a cochain morphism, hence induces a map $[{\rm Res}_D]\colon H^\bullet(X,\mc{A}) \to H^{\bullet-\ell}(D,\mc{B};\det(E^*))$.
\end{thm}
\bp The $\mc{B}$-representation on $\det(E)$ is induced from the $\mc{B}$-connection $\nabla^\sigma = [\sigma(\cdot),\cdot]_\mc{A}$. The condition that ${\rm Res}_D$ is a cochain morphism is equivalent to the fact that $\det(E^*) \subseteq \wedge^\bullet \mc{A}^*|_D$ has closed local generating sections. This is implied by \autoref{prop:berrephtpy}, giving the result.
\ep
\begin{rem} The above theorem is not true in general for the lower residue maps discussed in \autoref{sec:residuemaps}. Counterexamples arise from the local generators of $\wedge^{\ell-k} E^*$ not being closed. However, by analyzing the proof of \autoref{thm:residues} we see that it is true for ${\rm Res}_{D,-k}$ whenever the representation up to homotopy $\wedge^{\ell-k} E$ is a $\mc{B}$-module. This property is inherited (i.e.\ if it holds for some $k$, then also for all $k' \leq k$), and the above shows that this is true for $k = 0$.
\end{rem}
To better understand the Lie algebroid cohomology acting as the target for the residue map, it is clear that it is important to determine when, or whether, the $\mc{B}$-representation on $\det(E)$ is unimodular, so that one can suitably trivialize $\det(E)$ as a representation.
\begin{rem}\label{rem:residuesymplecticliealgebroid} Let $(\mc{A},\omega_{\mc{A}}) \to X$ be a symplectic Lie algebroid and let $D \subseteq X$ be a transitive $\mc{A}$-invariant submanifold, possibly the degeneracy locus of $\mc{A}$ (see \autoref{prop:degenlocus}). The residue map gives a form ${\rm Res}_D(\omega_\mc{A}) \in \Omega^{2-\ell}(D; \wedge^\ell E^*)$, where $E = \ker \rho_{\mc{A}}|_D$ and $\ell = \dim(E)$. Note however that it can only be nonzero for $\ell \leq 2$. One further has the residue ${\rm Res}_D(\omega_{\mc{A}}^n) \in \Omega^{2n-\ell}(D;\wedge^\ell E^*)$, if ${\rm rank}(\mc{A}) = 2n$. These describe important features of the induced geometry on $D$, and the latter has a Poisson-geometric interpretation which we discuss in \autoref{sec:residuesympla}.
\end{rem}
\begin{rem} While \autoref{thm:residues} shows the residue map descends to cohomology, it often does not respect the ring structure induced by the wedge product. See also \autoref{sec:residuesympla}.
\end{rem}
\subsubsection{Examples of residue maps}
\label{sec:residueexamples}
In this section we discuss several examples of residue maps, for some of the Lie algebroids of divisor-type that we have discussed in this paper.

\begin{exa}[$\mc{A}_Z^k$]\label{exa:bkersiduemaps} For the $b^k$-bundles $\mc{A}_Z^k \to X$ of \autoref{exa:bktangentbundle}, the hypersurface $Z \subseteq X$ is a transitive invariant submanifold. The resulting extension sequence of Lie algebroids
\be
	0 \to \mathbb{L}_{Z,k} \to \mc{A}_Z^k|_Z \to TZ \to 0,
\ee
is such that $\mathbb{L}_{Z,k} \to Z$ is trivial, and has a canonical trivializing section (it is also called the \emph{$b^k$-normal bundle}, see \cite[Proposition 4]{GuilleminMirandaPires14} and \cite[Proposition 3.2]{Scott16}). As such the residue map is a cochain morphism, given by ${\rm Res}_{Z,k}\colon \Omega^\bullet(\mc{A}_Z^k) \to \Omega^{\bullet-1}(Z)$, and locally $\frac{dx}{x^k} \wedge \alpha + \beta \mapsto \iota_Z^*(\alpha)$, where $x \in j_{k-1}$ is a local $(k-1)$-jet generator, $\alpha,\beta \in \Omega^\bullet(X)$ and $\iota_Z\colon Z \hookrightarrow X$ is the inclusion. These residue maps along with the dual $\mc{A}_Z^{k-1}$-anchor $\varphi_{A_Z^k}^*$ fit in the short exact sequence
\be
	0 \to \Omega^\bullet(\mc{A}_Z^{k-1}) \to \Omega^\bullet(\mc{A}_Z) \to \Omega^{\bullet-1}(Z) \to 0,
\ee
which is called the \emph{residue sequence} for $\mc{A}_Z^k$. These sequences split and imply the cohomological consequences $H^\bullet(\mc{A}_Z^k) \cong H^\bullet(\mc{A}_Z^{k-1}) \oplus H^{\bullet-1}(Z)$ (c.f.\ \cite{GuilleminMirandaPires14,MarcutOsornoTorres14,Melrose93,Scott16}; see \cite{KlaasseLanius18} for more).
\end{exa}
 There are no lower residue maps for $\mc{A}_Z^k \to X$, as $\mathbb{L}_{Z,k}$ is one-dimensional. While \autoref{thm:residues} asserts the existence of the above residue map as stated, it is special in this situation that $\mathbb{L}_{Z,k}$ is not only (canonically) trivial, but is trivial as a $TZ$-representation, so that the image of the residue map is smooth forms on $Z$, i.e.\ $\Omega^{\bullet-1}(Z)$, instead of $\Omega^{\bullet-1}(Z;\mathbb{L}_{Z,k}^*)$.
\begin{exa}[$\mc{A}_{|D|}$]\label{exa:elltangentresidues} For the elliptic tangent bundle $\mc{A}_{|D|} \to X$ of \autoref{exa:ellvx}, the degeneracy locus $D \subseteq X$ is a transitive invariant submanifold, hence leads to an extension (\cite{CavalcantiGualtieri18})
\be
	0 \to \underline{\mathbb{R}} \oplus \mf{k} \to \mc{A}_{|D|}|_D \to TD \to 0.
\ee
Note that $\ker(\rho_{\mc{A}_{|D|}}) \cong \underline{\mathbb{R}} \oplus \mf{k}$, so that $\det(\ker(\rho_{\mc{A}_{|D|}})) \cong \mf{k}$. This results in an \emph{elliptic residue}
\be	
	{\rm Res}_q\colon \Omega^\bullet(\mc{A}_{|D|}) \to \Omega^{\bullet-2}(D;\mf{k}^*), \qquad d\log r \wedge d\theta \wedge \alpha + d\log r \wedge \beta + d\theta \wedge \gamma + \eta \mapsto \iota_D^*(\alpha),
\ee
for $\alpha,\beta,\gamma,\eta \in \Omega^\bullet(X)$ and $\iota_D\colon D \hookrightarrow X$ the inclusion. Letting $\Omega^\bullet_0(\mc{A}_{|D|}) = \ker({\rm Res}_q)$ we have
\be
	{\rm Res}_r\colon \Omega^\bullet_0(\mc{A}_{|D|}) \to \Omega^{\bullet-1}(D), \qquad d\log r \wedge \beta + d\theta \wedge \gamma + \eta \mapsto \iota_D^*(\beta),
\ee
which is called the \emph{radial residue}. A coorientation for $D$ trivializes $\mf{k}$ and gives further
\be
	{\rm Res}_\theta\colon \Omega^\bullet_0(\mc{A}_{|D|}) \to \Omega^{\bullet-1}(D), \qquad d\log r \wedge \beta + d\theta \wedge \gamma + \eta \mapsto \iota_D^*(\gamma),
\ee
which is the \emph{$\theta$-residue}. In general ${\rm Res}_r\colon \Omega^\bullet(\mc{A}_{|D|}) \to \Omega^{\bullet-1}(D,{\rm At}(S^1 ND))$ and ${\rm Res}_q = \iota_{\partial_{\theta}} \circ {\rm Res}_r$. This is because $\mc{A}_{|D|}|_D$ is an extension by $\underline{\R}$ of the Atiyah algebroid ${\rm At}(S^1 ND) \to D$, i.e.\
\be
	0 \to \underline{\R} \to \mc{A}_{|D|}|_D \to {\rm At}(S^1 ND) \to 0,
\ee
so that by the general recipe for residue maps we obtain the map ${\rm Res}_r$ above, as ${\rm det}(\underline{\R}^*)$ is the trivial representation. Finally ${\rm At}(S^1 ND)$ has $D$ as transitive invariant submanifold, with
\be
	0 \to \mf{k} \to {\rm At}(S^1 ND) \to TD \to 0,
\ee
whose residue map ${\rm Res}_{\rm At}\colon \Omega^\bullet(D, {\rm At}(S^1 ND)) \to \Omega^{\bullet-1}(D; \mf{k}^*)$ combines as ${\rm Res}_q = {\rm Res}_{\rm At} \circ {\rm Res}_r$, and ${\rm Res}_{\rm At} = \iota_{\partial_{\theta}}$ in local coordinates.
As in \autoref{exa:bkersiduemaps} there are cohomological consequences one can draw from this (see \cite{CavalcantiGualtieri18,CavalcantiKlaasse18}), for example $H^\bullet(\mc{A}_{|D|}) \cong H^\bullet(X\backslash D) \oplus H^{\bullet-1}(S^1 ND)$.
\end{exa}
The previous example is special as the isotropy bundle is two-dimensional, but splits as a sum of two line bundles. As such we have ${\rm Res}_{D,0} = {\rm Res}_q$ and ${\rm Res}_{D,-1} =  {\rm Res_r} + {\rm Res}_\theta$ here. There are also residue maps for complex Lie algebroids, for example the complex log-tangent bundle $\mc{A}_D$ of \autoref{exa:complexlogtgnt} with its complex log residue ${\rm Res}_D\colon \Omega^\bullet(\mc{A}_D) \to \Omega^{\bullet-1}(D;\C)$, which is related to the \emph{complex residue} ${\rm Res}_\C := {\rm Res}_{D,-1}$ of the elliptic tangent bundle $\mc{A}_{|D|}$ (see \cite{CavalcantiGualtieri18}).
\begin{exa}[$\mc{A}_{\underline{Z}}$]  Let $\underline{Z} = \cup_{j \in I} Z_j$ be a normal-crossing log divisor with associated log-tangent bundle $\mc{A}_{\underline{Z}}$ as below \autoref{exa:zvx}. Set $Z_\tau = \cap_{i \in \tau} Z_i$ for $\tau \subseteq I$ and let $\underline{Z}_k = \cup_{|\tau| = k} Z_\tau$, so that $\underline{Z}_k$ consists of all $k$-fold intersections of hypersurfaces. By transversality, these are all smooth submanifolds, with induced normal-crossing log divisors. We have for $\mc{A}_{\underline{Z}_k} \to \underline{Z}_{k}$ that $\underline{Z}_{k+1}$ is $\mc{A}_{\underline{Z}_k}$-invariant, with its associated Lie algebroid extension sequence over $\underline{Z}_{k+1}$ being
\be
	0 \to \mathbb{L}_{\underline{Z},k} \to \mc{A}_{\underline{Z}_k}|_{\underline{Z}_{k+1}} \to \mc{A}_{\underline{Z}_{k+1}} \to 0.
\ee
Again $\mathbb{L}_{\underline{Z},k} \to \underline{Z}_{k+1}$ is canonically trivial, leading to ${\rm Res}_{Z,k}\colon \Omega^\bullet(\underline{Z}_k, \mc{A}_{\underline{Z}_k}) \to \Omega^{\bullet-1}(\underline{Z}_{k+1},\mc{A}_{\underline{Z}_{k+1}})$. This leads to a sequence of composable residue maps, and splitting residue sequences
\be
	0 \to \Omega^\bullet(\underline{Z}_k) \to \Omega^\bullet(\underline{Z}_k, \mc{A}_{\underline{Z}_k}) \to \Omega^{\bullet-1}(\mc{A}_{\underline{Z}_{k+1}}) \to 0,
\ee
which leads to the cohomological consequence $H^\bullet(\underline{Z}_k,\mc{A}_{\underline{Z}_k}) \cong H^\bullet(\underline{Z}_k) \oplus H^{\bullet-1}(\underline{Z}_{k+1}, \mc{A}_{\underline{Z}_{k+1}})$. These are straightforward analogues of the case of the log-tangent bundle $\mc{A}_Z$ of \autoref{exa:bkersiduemaps}. By applying induction to $k$ and noting that $X = Z_0$, this gives (c.f.\ \cite[Appendix A.24]{GualtieriLiPelayoRatiu17})
\be
	H^\bullet(\mc{A}_{\underline{Z}}) \cong H^\bullet(X) \oplus \bigoplus_{k \geq 1} H^{\bullet-k}(\underline{Z}_k).
\ee
\end{exa}
The above can also be done in similar fashion for self-crossing log divisors (see \cite{MirandaScott18}) instead mapping onto its $k$-strata, but this requires more notation, so that we will not do so here. The other Lie algebroids found in this paper, such as the scattering tangent bundle $\mc{C}_Z \to X$, also admit residue maps, but these are less useful, because $Z$ is $\mc{C}_Z$-invariant with $\mc{B} = 0_{TZ}$.
\begin{exa}[$\mc{A}_{Z,F}$] Consider the Lie algebroid $\mc{A}_{Z,F} := [TX\:TF] \to X$ given a codimension-one involutive distribution $TF \subseteq TZ$ on a hypersurface $Z \subseteq X$ (see \autoref{sec:examplesmodification}). Then $\mc{A}_{Z,F}$ has $Z$ as a projective invariant submanifold with $\mc{B} = TF$ (i.e.\ $\mc{A}_{Z,F}|_Z$ is regular):
\be
	0 \to \ker(\rho_{\mc{A}_{Z,F}}|_Z) \to \mc{A}_{Z,F}|_Z \to TF \to 0.
\ee
Locally we have $\ker(\rho_{\mc{A}_{Z,F}}|_Z) = \langle x \partial_x, x \partial_y\rangle$ for $x \in I_Z$ a local generator, and $\partial_y \subseteq \Gamma(TZ)$ normal to $TF$. The dual generator $dx/x$ is closed, but the generator $dy/x$ is not. Thus, while the top residue map is a cochain morphism by \autoref{thm:residues}, the lower residue map is not. Note that $\mc{A}_{Z,F}$ is an $\mc{A}_Z$-Lie algebroid, and that its $\mc{A}_Z$-anchor is a quasi-isomorphism (c.f.\ \cite{GualtieriLiPelayoRatiu17,KlaasseLanius18}).
\end{exa}
Similar ideas apply to invariantly capture what is done in \cite{Lanius16}, but we will not do so here.
\begin{exa}[$\mc{A}_W$]\label{exa:ellipticlogresidue} Let $W = Z \otimes |D|$ be an elliptic-log divisor with associated Lie algebroid $\mc{A}_W \to X$ as in \autoref{prop:elllogalgebroid}. The degeneracy locus $Z \subseteq X$ is $\mc{A}_W$-invariant, and $D \subseteq Z$ is a hypersurface with log-tangent bundle $\mc{A}_{Z,D} \to Z$. There is an extension sequence
\be
	0 \to E \to \mc{A}_W|_Z \to \mc{A}_{Z,D} \to 0,
\ee
with $E$ the rank-$1$ germinal isotropy bundle. This follows from the local description of $\Gamma(\mc{A}_W) = \langle x \partial_x + y \partial_y, x(y \partial_x - x \partial_y) \rangle \oplus \Gamma(TD)$. This results in a residue map
\be
	{\rm Res}_Z\colon \Omega^\bullet(\mc{A}_W) \to \Omega^{\bullet-1}(Z,\mc{A}_{Z,D};E^*).
\ee
The submanifold $D\subseteq X$ is also $\mc{A}_W$-invariant, with a transitive extension sequence
\be
	0 \to \ker(\rho_{\mc{A}_W}|_D) \to \mc{A}_W|_D \to TD \to 0.
\ee
This gives another residue map onto $D$, namely
\be
	{\rm Res}_D\colon \Omega^\bullet(\mc{A}_W) \to \Omega^{\bullet-2}(D;\det(F^*)),
\ee
with isotropy plane bundle $F = \ker(\rho_{\mc{A}_W}|_D) \to D$. These residue maps are related, because for the log-tangent bundle $\mc{A}_{Z,D} \to Z$, the submanifold $D \subseteq Z$ is invariant with sequence
\be
	0 \to \wt{E} \to \mc{A}_{Z,D}|_D \to TD \to 0,
\ee
which in turn results in the logarithmic residue map
\be
	{\rm Res}_{Z,D}\colon \Omega^\bullet(Z,\mc{A}_{Z,D}) \to \Omega^{\bullet-1}(D;\wt{E}^*).
\ee
We have $\det(F^*) \cong E^* \otimes \wt{E}^*$ and through this we have that ${\rm Res}_D = {\rm Res}_{Z,D} \circ {\rm Res}_Z$. Note that the line bundle $\wt{E} = \mathbb{L}_D \to D$ is canonically trivial as in \autoref{exa:bkersiduemaps}. The cohomological consequences of the existence of these residue maps is discussed in \cite{KlaasseLanius18}.
\end{exa}
Further examples of residue maps for other Lie algebroids can be found in \cite{KlaasseLanius18}.
\subsection{Poisson modules and modular residue maps}
\label{sec:poissonmodulesmodfoliation}
In this section we discuss modules associated to $\mc{A}$-Poisson structures. For simplicity we will mainly focus on the case where $\mc{A} = TX$ \cite{GualtieriLi14,GualtieriPym13,Polishchuk97}. In essence, given a Lie algebroid $\mc{A} \to X$, $\mc{A}$-Poisson modules are Lie algebroid modules for the $\mc{A}$-Poisson algebroid $\mc{A}^*_{\pi_\mc{A}} \to X$ of an $\mc{A}$-Poisson structure $\pi_\mc{A} \in {\rm Poiss}(\mc{A})$. Unravelling \autoref{defn:amodule} gives the following definition of an $\mc{A}$-Poisson module.
\begin{defn} Let $\pi_\mc{A} \in {\rm Poiss}(\mc{A})$. A \emph{$\pi_\mc{A}$-module} is a bundle $E \to X$ equipped with a flat $\mc{A}$-Poisson connection, i.e.\ a linear morphism $\nabla\colon \Gamma(E) \to \Gamma(\mc{A} \otimes E)$ satisfying the Leibniz rule
	\be
	\nabla(f s) = - \pi_{\mc{A}}^\sharp(d_\mc{A} f) \otimes s + f \nabla s,
	\ee
	for all $f \in C^\infty(X)$ and $s \in \Gamma(E)$, and which for all $\alpha,\beta \in \Gamma(\mc{A}^*_{\pi_{A}})$ satisfies (with $\nabla_\alpha s = (\nabla s)(\alpha)$)
	\be
		\nabla_{[\alpha,\beta]_{\pi_\mc{A}}} = \nabla_\alpha \circ \nabla_\beta - \nabla_\beta \circ \nabla_\alpha.
	\ee
\end{defn}
When $\mc{A} = TX$, the bracket on $T^*_\pi X$ is described directly by the relation $[df,dg]_\pi = d\{f,g\}_\pi$. As such, in this case flatness is given by $\nabla_{d\{f,g\}_\pi} = \nabla_{df} \circ \nabla_{dg} - \nabla_{dg} \circ \nabla_{df}$ as operators on $\Gamma(E)$.
\begin{exa}\label{exa:modularreppoisson} Given $\pi \in {\rm Poiss}(X)$, the Lie algebroid $T^*_\pi X$ has a representation on $\det(T^*X)$ which is given, for $\alpha \in \Gamma(T^*_\pi X)$ and $\mu \in \Gamma(\det(T^*X))$, by
\be
	\nabla_\alpha \mu = \mc{L}_{\pi^\sharp(\alpha)} \mu + (\pi,d\alpha) \mu = [\alpha,\mu]_\pi - (\pi,d\alpha)\mu = \alpha \wedge d(\iota_\pi \mu).
\ee
Note that this is essentially the modular representation $Q_{T^*_\pi X} \cong \det(T^*_\pi X)^2$ of $T^*_\pi X$.
\end{exa}
The \emph{modular class} ${\rm Mod}(\pi)$ of $\pi \in {\rm Poiss}(X)$ is given by $2\, {\rm Mod}(\pi) = {\rm Mod}(T^*_\pi X) = \theta_{Q_{T^*_\pi X}}$.
Given a $\pi_{\mc{A}}$-line bundle $(L,\nabla)$, a local trivialization $s \in \Gamma(L)$ specifies a unique $\mc{A}$-Poisson section $v \in {\rm Poiss}(\pi_{\mc{A}}) \subseteq \Gamma(\mc{A})$ by $\nabla s = v \otimes s$. In fact, $v$ is $\mc{A}$-Poisson if and only if $\nabla$ is flat.
\begin{rem} For $\pi \in {\rm Poiss}(X)$, in case $L = \det(T^*X)$ with local covolume form $\mu$, the above section $v = V_\mu \in \Gamma(TX)$ is called the \emph{$\pi$-modular vector field} associated to $\mu$. For general Lie algebroids $\mc{A}$, the above section then gives the \emph{$\pi_{A}$-modular section} $V_{\mu_\mc{A}}$ associated to $\mu_\mc{A}$.
\end{rem}
\begin{exa} Given $\pi \in {\rm Poiss}(X^n)$ locally expressed as $\pi = \sum_{ij} \pi^{ij} \partial_{x_i} \wedge \partial_{x_j}$, the modular vector field associated to $\mu = \partial_{x_1} \wedge \dots \wedge \partial_{x_n}$ is given by $V_\mu = \sum_i \partial_k(\pi^{ik}) \partial_{x_i}$.
\end{exa}
We defined $\mc{A}$-Poisson submanifolds as those equipped with Lie subalgebroids $(Y,\mc{A}',\pi_{\mc{A}'})$ for which the inclusion is an $\mc{A}$-Poisson map. We now specify to the case $\mc{A} = TX$ to discuss restriction of Poisson modules. Recall that $Y \subseteq (X,\pi)$ is $\pi$-Poisson if $I_Y$ is a $\pi$-Poisson ideal.
\begin{defn}[\cite{GualtieriPym13,Polishchuk97}] Let $\pi \in {\rm Poiss}(\mc{A})$ be given. A submanifold $Y \subseteq X$ is a \emph{strong $\pi$-Poisson submanifold} if $\mc{L}_V I_Y \subseteq I_Y$ for all $V \in {\rm Poiss}(\pi)$.
\end{defn}
As every $\pi$-Hamiltonian vector field is $\pi$-Poisson, it is clear that strong $\pi$-Poisson submanifolds are also $\pi$-Poisson submanifolds. In general, a \emph{strong $\pi$-Poisson ideal} is an ideal $I \subseteq C^\infty(X)$ that is preserved by all $\pi$-Poisson vector fields. \autoref{rem:preservingideals} shows that strong $\pi$-Poisson ideals are closed under sums and intersections, and \autoref{rem:preservingradical} shows the same for taking radicals. As such, we see that (strong) $\pi$-Poisson submanifolds are closed under unions and intersections. There moreover is a strengthening of \autoref{prop:pidegenpipoisson} (c.f.\ \cite{GualtieriPym13,Polishchuk97}). 
\begin{prop}[\cite{Polishchuk97}]\label{prop:pidegpistrongpoisson} Let $\pi \in {\rm Poiss}(X)$. Then the degeneracy ideals $I_{\pi,2k}$ are strong $\pi$-Poisson ideals, hence the degeneracy loci $X_{\pi,2k}$ are strong $\pi$-Poisson submanifolds if smooth.
\end{prop}
\bp Unravelling definitions, for $k \geq 0$, the ideal $I_{2k,\pi}$ is defined as the image $\pi^k(\Gamma(\wedge^{2k} T^*X))$. Let $V \in {\rm Poiss}(\pi)$ and $\alpha \in \Gamma(\wedge^{2k} T^*X)$. Then we readily compute that
\be
	\mc{L}_V(\pi^k(\alpha)) = (\mc{L}_V \pi^k)(\alpha) + \pi^k(\mc{L}_V \alpha) = \pi^k(\mc{L}_V \alpha) \in I_{2k,\pi}.\qedhere
\ee
\ep
Essentially, the above follows because the degeneracy loci are defined directly in terms of the (rank of the) Poisson bivector itself: all $\pi$-Poisson vector fields preserve all powers of $\pi$. Turning this on its head, the above proposition says that given $\pi \in {\rm Poiss}(X)$, all $\pi$-Poisson vector fields must be tangent to all the degeneracy loci of $\pi$.

Consider a submanifold $Y \subseteq X$ and the sequence of its conormal bundle:
\be
	0 \to N^*Y \to T^*X|_Y \to T^*Y \to 0.
\ee
When $Y$ is a $\pi$-Poisson submanifold, this becomes an extension sequence of Lie algebroids, using the Poisson algebroids $T^*X_\pi$ and $T^*Y_{\pi_Y}$, where $\pi_Y \in {\rm Poiss}(Y)$ is the induced Poisson structure on $Y$. One naturally wonders about restriction of $\pi$-Poisson line bundles. We have:
\begin{prop}[{\cite[Lem 4.11]{GualtieriPym13}}]\label{prop:degenlocusstrongpoisson} Let $\pi \in {\rm Poiss}(X)$ and let $(L,\nabla)$ be a $\pi$-Poisson line bundle, and let $Y \subseteq X$ be a strong $\pi$-Poisson submanifold with induced Poisson structure $\pi_Y$. Then the restriction $(L|_Y, \nabla|_Y)$ is a $\pi_Y$-Poisson module.
\end{prop}
\bp Let $s \in \Gamma(L)$ be a local trivialization. Then $\nabla s = v \otimes s$ for a unique $\pi$-Poisson vector field $v \in {\rm Poiss}(\pi)$. To be able to restrict $\nabla s$ to $Y$, the image $[v \otimes s] \in \Gamma(Y,{\rm Hom}(N^*Y \otimes L|_Y))$ is zero. As $Y$ is a strong $\pi$-Poisson submanifold we have $\mc{L}_v I_Y \subseteq I_Y$. Thus $v$ is tangent to $Y$, so that $[v \otimes s] \equiv 0$. We conclude that $v \otimes s|_Y \in \Gamma(Y; T^*Y \otimes L|_Y)$ as required.
\ep
As explained in \cite{GualtieriPym13}, a $\pi$-Poisson line bundle $(L,\nabla)$ has several \emph{residue maps}. Given $k \geq 0$, with $Y := X_{\pi,2k}$, these are sections ${\rm Res}_k(\nabla) \in \Gamma(Y,\wedge^{2k+1} TY)$. They are locally given using a trivialization $s \in \Gamma(L)$, so that $\nabla s = v \otimes s$ for $v$ a local $\pi$-Poisson vector field, by the expression
\be
	{\rm Res}_k(\nabla) = v \wedge \pi^k_Y.
\ee
Invariantly, they are obtained by the composition
\be
	(\cdot \wedge \sigma^k) \circ \nabla\colon \Gamma(L) \to \Gamma(TX \otimes L) \to \Gamma(\wedge^{2k+1}TX \otimes L),
\ee
which is then restricted to $Y$ using \autoref{prop:degenlocusstrongpoisson} to give a morphism
\be
	\Gamma(Y, L|_Y) \to \Gamma(Y, \wedge^{2k+1}TY \otimes L|_Y),
\ee
from which the multiderivation ${\rm Res}_k(\nabla) \in \Gamma(Y, \wedge^{2k+1} TY)$ is extracted. During this process it is used that $\pi^{k+1}|_Y \equiv 0$ because $Y$ is the $2k$th-degeneracy locus of $\pi$. Applying the above to $L = \det(T^*X)$ we obtain the \emph{$k$th modular residue} ${\rm Res}_{k,{\rm mod}}(\pi)$.

In \autoref{sec:residuesympla} we briefly discuss these residue maps and their relation to those in \autoref{sec:residuemaps} in the context of Poisson structures of nondegenerate $\mc{A}$-type for a Lie algebroid $\mc{A} \to X$.
\subsection{The modular foliation}
\label{sec:modularfoliation}
In this section we discuss the modular foliation of a Poisson manifold, which was introduced in \cite{GualtieriPym13}. Our main goal here is to prove \autoref{thm:modfoliation} which discusses how the modular foliation interacts with the lifting procedure of \autoref{sec:liftingpoisson}.

Let $\pi \in {\rm Poiss}(X)$ be given and let $(L,\nabla) \to X$ be a $\pi$-Poisson line bundle. Then ${\rm At}(L) \to X$, its \emph{Atiyah algebroid}, is a transitive Lie algebroid. There is a natural section $\sigma_\nabla \in \Gamma(\wedge^2 {\rm At}(L))$, which is in fact an ${\rm At}(L)$-Poisson structure (\cite[Corollary 5.3]{Polishchuk97}). It can be viewed as the natural Poisson structure on ${\rm Tot}(L)$ induced by $\pi$ and $\nabla$. In our language, $\sigma_\nabla$ is a Poisson lift of the underlying Poisson structure $\pi$. Due to this, $({\rm At}(L), {\rm At}^*(L), \sigma_\nabla)$ is a triangular Lie bialgebroid (c.f.\ \autoref{rem:apoissontriangbialgebroid}). In \autoref{sec:poissondivtype} we studied cases where these were of divisor-type, but here instead the Lie algebroid ${\rm At}(L)$ is transitive. We have the associated module $\mc{F}_{\sigma_\nabla} = \sigma_\nabla^\sharp(\Gamma({\rm At}(L)^*) \subseteq \Gamma({\rm At}(L))$, which can be pushed down using $\rho\colon {\rm At}(L) \to TX$.
\begin{defn}[\cite{GualtieriPym13}] Let $\pi \in {\rm Poiss}(X)$ and $(L,\nabla) \to X$ a $\pi$-Poisson line bundle. Then $\mc{F}_{\pi,(L,\nabla)} := \rho(\mc{F}_{\sigma_\nabla}) \subseteq \Gamma(TX)$ is the singular foliation associated to $(L,\nabla)$.
\end{defn}
It is immediate that the module $\mc{F}_{\pi,(L,\nabla)}$ is involutive as $\rho \circ \sigma_\nabla^\sharp$ is a Lie algebroid morphism.
\begin{defn}[\cite{GualtieriPym13}] Let $\pi \in {\rm Poiss}(X)$. Then the \emph{modular foliation} $\mc{F}_{\pi,{\rm mod}}$ is the singular foliation associated to the modular representation $L = \det(T^*X)$ of $\pi$ (see \autoref{exa:modularreppoisson}).
\end{defn}
Let us describe the modular foliation locally, as in \cite{PymSchedler18}. Let $\pi \in {\rm Poiss}(X)$ and let $\mu$ be a local volume form. Then the associated modular vector field $V_{\pi,\mu}$ of $\pi$ is characterized by
\be
	\nabla(\mu) = V_{\pi,\mu} \otimes \mu, \qquad \text{or} \qquad \mc{L}_{\pi^\sharp(df)} \mu = (\mc{L}_{V_{\pi,\mu}} f) \cdot \mu.
\ee
Moreover, it is immediate that $V_{\pi,\mu} \in {\rm Poiss}(\pi)$. Consequently, $V_{\pi,\mu}$ is always tangent to the degeneracy loci of $\pi$ by \autoref{prop:pidegenpipoisson}. The same computation showing that the characteristic class $\theta_L$ is independent of local generators shows that if we change volume element to $\mu' = f \mu$ for some nonvanishing function $f$, we obtain that $V_{\pi,\mu'} = V_{\pi,\mu} - V_{\log f}$, i.e.\ the associated modular vector fields differ by a Hamiltonian vector field. As a result we have
\be
	\mc{F}_{\pi,{\rm mod}} = \mc{F}_\pi \cup \{V_{\pi,\mu} \, | \, \mu \in {\rm Vol}_{\rm loc}(X)\}.
\ee
The involutivity of $\mc{F}_{\pi,{\rm mod}}$ locally becomes the computation that, given $f \in C^\infty(X)$, we have:
\be
	[V_{\pi,\mu}, V_f] = \mc{L}_{V_{\pi,\mu}}(\pi^\sharp(df)) = \pi^\sharp(d \mc{L}_{V_{\pi,\mu}} f) = V_{\mc{L}_{V_{\pi,\mu}} f}.
\ee
We see that $\mc{F}_\pi \subseteq \mc{F}_{\pi,{\rm mod}}$. The leaves of the former are the $\pi$-symplectic leaves, which are thus necessarily even-dimensional. The leaves of the modular folation come in two flavors:
\bi
	\item The even-dimensional leaves, where $V_{\pi,{\rm mod}}$ is tangent to the symplectic leaves;
	\item The odd-dimensional leaves, where $V_{\pi,{\rm mod}}$ is transverse to the symplectic leaves.
\ei
Let us record the consequence of \autoref{prop:pidegpistrongpoisson} regarding strong Poisson submanifolds.
\begin{prop}\label{prop:modvfielddegenloci} Let $\pi \in {\rm Poiss}(X)$. Then any local modular vector field for $\pi$ is tangent to all strong Poisson submanifolds of $\pi$, hence in particular to all of its degeneracy loci.
\end{prop}
The following (combined with \autoref{prop:atypepoissonsubmfd} showing $\mc{F}_\pi \subseteq \mc{F}_\mc{A}$) is \autoref{thm:intromodfoliation}.
\begin{thm}\label{thm:modfoliation} Let $\pi \in {\rm Poiss}(X)$ be of $\mc{A}$-type. Then $\mc{F}_\pi \subseteq \mc{F}_{\pi,{\rm mod}} \subseteq \mc{F}_\mc{A}$.
\end{thm}
\bp As noted, \autoref{prop:atypepoissonsubmfd} shows that $\mc{F}_\pi \subseteq \mc{F}_\mc{A}$, while we always have $\mc{F}_\pi \subseteq \mc{F}_{\pi,{\rm mod}}$. We are left with showing that all local modular vector fields are contained in $\mc{F}_\mc{A}$. Let $\pi_{A} \in {\rm Poiss}(\mc{A})$ be an $\mc{A}$-lift of $\pi$. Given local trivializing sections $\nu \in \Gamma(\det(\mc{A}^*))$ and $\mu \in \Gamma(\det(T^*X))$, note that $\nu \otimes \mu \in  \Gamma(\det(\mc{A}^*) \otimes \det(T^*X))$ can serve as a nonvanishing section for the modular representation $Q_{\mc{A}^*_{\pi_\mc{A}}}$ of the Poisson algebroid $\mc{A}^*_{\pi_\mc{A}}$. These now give rise to the modular section $v_{\nu \otimes \mu} \in \Gamma(\mc{A})$, but also to the $\pi_\mc{A}$-modular section $V_{\pi_\mc{A},\nu} \in \Gamma(\mc{A})$ and the $\pi$-modular vector field $V_{\pi,\mu} \in \Gamma(TX)$. By \cite[Equation (64)]{KosmannSchwarzbach00}, these are related by the equality
\be
	\rho_\mc{A}(v_{\nu \otimes \mu} - V_{\pi_\mc{A},\nu}) = V_{\pi,\mu}.
\ee
This shows that $V_{\pi,\mu} \in \mc{F}_\mc{A}$, from which $\mc{F}_{\pi,{\rm mod}} \subseteq \mc{F}_\mc{A}$ follows.
\ep
In \autoref{cor:poissonhamiltonian} we showed that if $\pi$ is of $I$-divisor-type, then any $\pi$-Poisson vector field admits a $TX_I$-lift. This shows that $\mc{F}_{\pi,{\rm mod}} \subseteq \Gamma(TX)_I = \rho_{TX_I}(\Gamma(TX_I))$ because any $\pi$-modular vector field is $\pi$-Poisson. Hence, it is not always necessary that $\pi$ is of $\mc{A}$-type.
\begin{rem} It is an interesting question to ask when the modular foliation $\mc{F}_{\pi,{\rm mod}}$ is projective. Note that $\mc{F}_\pi$ is projective if and only if $\pi$ is of $m$-divisor-type, i.e.\ is almost-regular (\autoref{cor:divtypealmostreg}). How is projectivity of $\mc{F}_{\pi,{\rm mod}}$ related to projectivity of $\mc{F}_\pi$?
\end{rem}
\subsubsection{Relations between modular classes}
There are extensions of the above concepts to $\mc{A}$-Poisson structures, using ideas contained in \cite{CaseiroFernandes13,KosmannSchwarzbachLaurentGengouxWeinstein08,KosmannSchwarzbach00, KosmannSchwarzbach08,KosmannSchwarzbachLaurentGengoux05,KosmannSchwarzbachWeinstein05}. Recall from \autoref{sec:aconnareps} that any Lie algebroid $\mc{A} \to X$ has a modular class ${\rm Mod}(\mc{A})$. As an $\mc{A}$-Poisson structure turns $\mc{A}^*$ into a Lie algebroid (\autoref{defn:apoissonalgebroid}), there are various modular classes around. In particular, the modular class ${\rm Mod}(\pi_{\mathcal{A}})$ of $\pi_{\mathcal{A}}$ is related to ${\rm Mod}(\pi)$ through the modular class of $\mc{A}$: an $\mc{A}$-Poisson structure $\pi_{\mathcal{A}}$ comes equipped with:
\bi
	\item A Lie algebroid morphism $\pi_{\mathcal{A}}^\sharp\colon \mc{A}^*_{\pi_{\mathcal{A}}} \to \mc{A}$;
	\item A Lie algebroid morphism $\rho_{A}\colon \mc{A} \to TX$;
	\item A Lie algebroid morphism $\pi^\sharp\colon T^*X_\pi \to TX$;
	\item A Lie algebroid comorphism $\rho_{A}^*\colon \mc{A}^*_{\pi_{\mathcal{A}}} \dashrightarrow T^*X_\pi$.
\ei
These are of course related by the lifting relation $\pi^\sharp = \rho_{A} \circ \pi_\mc{A}^\sharp \circ \rho_{A}^*$. As $\rho_{A}^*$ is base-preserving, it can also be considered as a Lie algebroid morphism from $T^*X_\pi$ to $\mc{A}^*_{\pi_\mc{A}}$. Each of these four Lie algebroids has a modular class, and each morphism has a relative modular class \cite{KosmannSchwarzbachWeinstein05}. These are related to (or define) the modular classes ${\rm Mod}(\pi_\mc{A})$ and ${\rm Mod}(\pi)$ as follows:
\bi
	\item $2\, {\rm Mod}(\pi_{\mathcal{A}}) = {\rm Mod}^{\pi_{\mathcal{A}}^\sharp}(\mc{A}^*_{\pi_\mc{A}},\mc{A}) = {\rm Mod}(\mc{A}^*_{\pi_\mc{A}}) - (\pi_\mc{A}^\sharp)^*({\rm Mod}(\mc{A}))$;
	\item $2\, {\rm Mod}(\pi) = {\rm Mod}^{\pi^\sharp}(T^*X_\pi, TX) = {\rm Mod}(T^*X_\pi) - (\pi^\sharp)^*({\rm Mod}(TX)) = {\rm Mod}(T^*X_\pi)$;
	\item ${\rm Mod}^{\rho_{A}}(\mc{A},TX) = {\rm Mod}(\mc{A}) - \rho_{A}^*({\rm Mod}(TX)) = {\rm Mod}(\mc{A})$, because ${\rm Mod}(TX) = 0$;
	\item ${\rm Mod}^{\rho_{A}^*}(T^*X_\pi, \mc{A}^*_{\pi_\mc{A}}) = {\rm Mod}(T^*X_\pi) - \rho_{A}({\rm Mod}(\mc{A}^*_{\pi_{\mathcal{A}}}))$, as $(\rho_\mc{A}^*)^* = \rho_\mc{A}$.
\ei
The equality $\pi^\sharp = \rho_{A} \circ \pi_\mc{A}^\sharp \circ \rho_{A}^*$ then implies the following relations:
	\begin{align*}
		2\, {\rm Mod}(\pi) &= {\rm Mod}^{\pi^\sharp}(T^*X_\pi, TX) = {\rm Mod}^{\rho_{A} \circ \pi_\mc{A}^\sharp \circ \rho_{A}^*}(T^*X_\pi, TX)\\
		&= {\rm Mod}^{\rho_{A}^*}(T^*X_\pi, \mc{A}^*_{\pi_\mc{A}}) + \rho_{A}({\rm Mod}^{\rho_\mc{A} \circ \pi_\mc{A}^\sharp}(\mc{A}^*_{\pi_\mc{A}}, TX))\\
		&= {\rm Mod}^{\rho_{A}^*}(T^*X_\pi, \mc{A}^*_{\pi_\mc{A}}) + \rho_{A}({\rm Mod}^{\pi_\mc{A}^\sharp}(\mc{A}_{\pi_\mc{A}}^*,\mc{A}) + (\pi_{A}^\sharp)^*({\rm Mod}^{\rho_\mc{A}}(\mc{A},TX)))\\
		&= {\rm Mod}^{\rho_{A}^*}(T^*X_\pi, \mc{A}^*_{\pi_\mc{A}}) + \rho_{A}(2\,{\rm Mod}(\pi_\mc{A}) + (\pi_{A}^\sharp)^*({\rm Mod}(\mc{A}))).
	\end{align*}
This type of description can possibly be used in the discussion of liftability of unimodular Poisson structures to unimodular Lie algebroids.
\subsubsection{Examples of modular foliations}
We finish by describing several examples of modular foliations for Poisson structures of divisor-type, starting locally in dimension two.
\begin{exa} As in \cite[Example 3.1]{PymSchedler18}, given $\pi = f \partial_x \wedge \partial_y$ on $(\R^2,(x,y))$, its modular vector field associated to $\mu = \partial_x \wedge \partial_y$ is given by $V_{\pi,\mu} = (\partial_x f) \partial_y - (\partial_y f) \partial_x$. The symplectic foliation is given by the subset where $f$ is nonvanishing, and the points where $f$ vanishes. The modular vector field is tangent to the zero set of $f$, and vanishes exactly at its critical points.
\end{exa}
\begin{exa}[Nondegenerate] Let $\pi \in {\rm Poiss}(X)$ be nondegenerate. Then $X$ is a symplectic leaf (or: its connected components). This shows that $\mc{F}_\pi = \Gamma(TX)$, as is also immediate as $\pi^\sharp\colon T^*X \to TX$ is an isomorphism. The only degeneracy locus of $\pi$ is $X$ itself, and any modular vector field must be tangent to it by \autoref{prop:modvfielddegenloci}. Thus $\mc{F}_{\pi,{\rm mod}} = \mc{F}_{\pi} = \Gamma(TX)$. In fact, nondegenerate Poisson manifolds are unimodular, so that their modular class vanishes and there exists a global volume form (the Liouville form) which is preserved by Hamiltonian flows. For this volume form the modular vector field vanishes, which again shows that $\mc{F}_\pi = \mc{F}_{\pi, {\rm mod}}$.
\end{exa}
\begin{exa}[Log-Poisson] Let $\pi \in {\rm Poiss}(X)$ be a log-Poisson structure. Then $X$ and $Z$ are the only degeneraci loci of $\pi$. From this it is immediate that any local modular vector field must be tangent to $Z$. However, in fact any local modular vector field is transverse to the symplectic leaves of the induced Poisson structure $\pi_Z$ \cite{GualtieriLi14,GuilleminMirandaPires14}. As a consequence we have that $\mc{F}_{\pi_Z,{\rm mod}} = \Gamma(TZ)$, while $\mc{F}_{\pi,{\rm mod}} = \Gamma(TX,TZ) = \Gamma(\mc{A}_Z)$. By \autoref{exa:logpoissonlift} we know that $\pi$ is of $\mc{A}_Z$-type, so that \autoref{thm:modfoliation} tells us that $\mc{F}_\pi \subseteq \mc{F}_{\pi,{\rm mod}} \subseteq \mc{F}_{\mc{A}_Z} = \Gamma(\mc{A}_Z)$. We see that the second inclusion is an equality, and that $\mc{F}_\pi$ and $\mc{F}_{\pi, {\rm mod}}$ are both projective. The symplectic foliation of $\pi_Z$ is regular of corank-1, with tangencies $TF \subseteq TZ$, so that $\mc{F}_\pi = \Gamma(TX,TF)$.
\end{exa}
\begin{exa}[$b^k$-Poisson] Let $\pi \in {\rm Poiss}(X)$ be a $b^k$-Poisson structure. Then the behavior of $\pi$ is similar for all $k \geq 1$, and in particular to log-Poisson structures (where $k = 1$). Again the only degeneraci loci of $\pi$ are $X$ and $Z$, but if $k > 1$, then any modular vector field is tangent to the symplectic leaves of $\pi_Z$. If locally $\pi = z^k \partial_z \wedge \partial_{x_1} + \omega_0^{-1}$ as in \autoref{exa:bkpoissondarboux}, then $V_{\pi,\mu} = k z^{k-1} \partial_{x_1}$ with respect to $\mu = dz \wedge dx_i$. In terms of this local description, we then have
\be
	\mc{F}_\pi = \langle z^k \partial_z, z^k \partial_{x_1} \rangle \oplus \Gamma(TF), \quad \mc{F}_{\pi,{\rm mod}} = \langle z^k \partial_z, z^{k-1} \partial_{x_1} \rangle \oplus \Gamma(TF), \quad \mc{F}_{\mc{A}_Z^k} = \langle z^k \partial_z, \partial_{x_1} \rangle \oplus \Gamma(TF),
\ee
 where $TF \subseteq TX$ are the tangencies to the corank-one symplectic foliation of $\pi_Z \in {\rm Poiss}(X)$.
\end{exa}
The modular foliations of elliptic and elliptic-log Poisson manifolds are described in \cite{KlaasseLanius18,KlaasseLi18}.
\begin{exa}[$m$-Log-Poisson]\label{exa:modfolmlog} Let $\pi \in {\rm Poiss}(X)$ be $m$-log-Poisson with associated regular distribution $D$. Then the degeneracy loci of $\pi$ are $X$ and $Z$, so that any local modular vector field must be tangent to $Z$. However, as $\pi$ is of $\mc{D}$-type by \autoref{prop:mdivtypelift}, it must further be tangent to $D$ due to \autoref{thm:modfoliation}. This suggests that the modular vector fields can detect how the orbits of $\mc{D}$ meet $Z$. In the local examples of \autoref{rem:mlogpoissonlocal} on $(\R^3,(x,y,z))$, with respect to the volume form $\mu = dx \wedge dy \wedge dz$, the modular vector fields of $\pi_1$, $\pi_2$ and $\pi_3$ are:
\be
	V_{\pi_1,\mu} = \partial_y, \quad V_{\pi_2,\mu} = 0, \quad V_{\pi_3,\mu} = - 2x \partial_y.
\ee
From this we see that the associated modules $\mc{F}_{\pi_i}$ and $\mc{F}_{\pi_i, {\rm mod}}$ are given by
\begin{align*}
	&\mc{F}_{\pi_1} = \langle x \partial_x, x \partial_y \rangle, \quad \mc{F}_{\pi_2} = \langle z \partial_x, z \partial_y \rangle, \quad \mc{F}_{\pi_3} = \langle (z-x^2) \partial_x, (z-x^2)\partial_y \rangle,\\
	&\mc{F}_{\pi_1,{\rm mod}} = \langle x \partial_x, \partial_y \rangle, \quad \mc{F}_{\pi_2,{\rm mod}} = \langle z \partial_x, z \partial_y \rangle \quad \mc{F}_{\pi_3,{\rm mod}} = \langle (z-x^2) \partial_x, z \partial_y, x \partial_y. \rangle
\end{align*}
\end{exa}
In general, the modular foliation of an almost-regular Poisson manifold behaves as follows: in \autoref{sec:liftingalmostreg} we saw that if $\pi \in {\rm Poiss}(X)$ is of $m$-$I$-divisor-type with associated involutive distribution $D = D_\pi$, then $\pi$ is of both $D$- and $TX_I$-type (assuming projectivity of $I$). Then, \autoref{thm:modfoliation} implies that $\mc{F}_{\pi, {\rm mod}} \subseteq \Gamma(\mc{D}_\pi)$ and $\mc{F}_{\pi, {\rm mod}} \subseteq \Gamma(TX_I) = \Gamma(TX)_I$. In particular any modular vector field must be tangent to $D$ and to the subset $Z_I \subseteq X$ (if it is smooth).
\subsection{Residues for symplectic Lie algebroids}
\label{sec:residuesympla}
In this section we discuss some of the consequences of the fact that symplectic Lie algebroids $\mc{A}$ have residue maps, which can be applied to an $\mc{A}$-symplectic form $\omega_{\mc{A}}$ or its powers. We recall again that residue maps, while often chain maps, typically do not respect the ring structure coming from the wedge product. We will discuss this process in more detail in \cite{Klaasse18} in the context of Dirac geometry, yet show here what happens in the case of $b^k$-Poisson or elliptic Poisson structures.

We will need the following definitions (\cite[Definitions 3.2.20 and 3.2.21]{OsornoTorres15}).
\begin{defn} Let $M^m$ be a manifold equipped with a $(m-2n = 2\ell)$-regular foliation $\mc{F}$.
	\bi
	\item an \emph{$n$-cosymplectic tuple} on $M$ is a tuple $(\alpha_1, \dots, \alpha_n, \beta)$ of forms, where $\beta \in \Omega^2_{\rm cl}(M)$ has constant rank $2\ell$, $\alpha_i \in \Omega^1_{\rm cl}(M)$ for all $i$, and $\alpha_1 \wedge \dots \wedge \alpha_n \wedge \beta^{\ell} \neq 0$;
	\item an \emph{$n$-Poisson tuple} on $M$ is a tuple $(V_1, \dots V_n, \pi)$, where $\pi \in {\rm Poiss}(M)$ is regular of rank $2\ell$, $V_i \in \Gamma(TM)$ are pairwise commuting $\pi$-Poisson vector fields, and $V_1 \wedge \dots V_n \wedge \pi^\ell \neq 0$.
	\ei
	We say the tuple is \emph{adapted to $\mc{F}$} if $T\mc{F}$ is spanned by (the kernel of) the $v_i$ (respectively $\alpha_i$).
\end{defn}
When $n=1$ the above is simply referred to as a \emph{cosymplectic} or \emph{corank-$1$ Poisson tuple}. The following correspondence is in full analogy with the case $n=1$ contained in \cite{GuilleminMirandaPires11}.
\begin{lem}[{\cite[Lemma 3.2.29]{OsornoTorres15}}]\label{lem:oneonecorrespondence} Let $(M,\mc{F})$ be a manifold with a corank-$2n$ foliation. There is a one-to-one correspondence between $\mc{F}$-adapted $n$-cosymplectic and $n$-Poisson tuples.
\end{lem}
We wish to reinterpet these structures from a Dirac geometric point of view. Recall that a \emph{Dirac structure} $E \subseteq \mathbb{T}M$ is an involutive Lagrangian subbundle of the standard Courant algebroid structure on $\mathbb{T}M = TM \oplus T^*M$. A Dirac structure can equivalently be described by a \emph{pure spinor line}, as follows (see \cite{AlekseevBursztynMeinrenken09,Gualtieri11}). Sections $v = V + \xi \in \Gamma(\mathbb{T}M)$ act on differential forms via Clifford multiplication, given by $v \cdot \rho = \iota_V \rho + \xi \wedge \rho$ for $\rho \in \Omega^\bullet(M)$. With this in mind, a pure spinor line is a line bundle $K_E \subseteq \wedge^\bullet T^*M$ that is pointwise generated by a pure spinor $\rho = e^B \wedge \Omega$ with $B$ a two-form and $\Omega$ a decomposable $k$-form, such that for any nonzero local section $\rho \in \Gamma(K_E)$ we have $d \rho = v \cdot \rho$ for some $v \in \Gamma(\mathbb{T}M)$. The correspondence  is then given by the fact that $E = {\rm Ann}(K_E)$ is the annihilator of $K_E$ under the Clifford action.
\begin{prop} Let $M^m$ be a manifold. Then an $n$-cosymplectic tuple $(\alpha_1,\dots,\alpha_n,\beta)$ defines a global pure spinor $\rho := e^\beta \wedge (\alpha_1 \wedge \dots \wedge \alpha_n)$ for which $\rho_m \neq 0 \in \Gamma(\wedge^m T^*M)$. 
\end{prop}
We call the associated pure spinor line $K_\rho = \langle \rho \rangle$ an \emph{$n$-cosymplectic structure} on $M$.
\bp This is almost immediate. Note that $\rho_m = \alpha_1 \wedge \dots \wedge \alpha_n \wedge \beta^\ell$, which is nonvanishing. Letting $\Omega := \alpha_1 \wedge \dots \wedge \alpha_n$, note that because $\beta$ and each $\alpha_i$ is closed we have
\be
	d \rho = d(e^\beta \wedge \Omega) = d\beta \wedge \rho + e^\beta \wedge d\Omega = 0 = 0 \cdot \rho.
\ee
The Dirac structure associated to $K_\rho = \langle \rho \rangle$ is given by $e^\beta( T\mc{F}) = \{v + \beta^\flat(v) \, | \, v \in \Gamma(T\mc{F})\}$.
\ep
From the Dirac perspective, the intrinsic object is the pure spinor line $K_\rho$, instead of $\rho$ itself. However, $n$-cosymplectic structures have a global generating pure spinor.  In general, this still means that the one-forms $\alpha_i$ need not exist separately, but only their associated volume form $\alpha_1 \wedge \dots \wedge \alpha_n \in \Gamma(\det(T\mc{F})^*)$ (making $\mc{F}$ coorientable). Their separate existence distinguishes $n$-cosymplectic structures from other symplectic foliations (such foliations are called \emph{unimodular}). Similarly, the two-form $\beta$ should be viewed as an extension to $M$ of a symplectic form on $\mc{F}$. However, $n$-cosymplectic structures are such that the symplectic form on $\mc{F}$ admits a closed global extension (such symplectic foliations are called \emph{strong}). Thus, an $n$-cosymplectic structure determines a coorientable and unimodular corank-$2n$ strong symplectic foliation, whose associated pure spinor line has a global generating pure spinor.

We now discuss how these structures arise in the context of log-Poisson structures, and of elliptic Poisson structures with zero elliptic residue. To start, recall that log-Poisson structures admit nondegenerate lifts to the log-tangent bundles that they define (\autoref{exa:logpoissonlift}). In other words, $\pi \in {\rm Poiss}(X)$ of $I_Z$-divisor-type has an associated dual $\mc{A}_Z$-symplectic structure $\omega$.
\begin{prop}\label{prop:logcosymplectic} Let $(X^{2n},Z,\pi,\omega)$ be a log-Poisson manifold. Then
\be
	\rho := {\rm Res}_Z(e^\omega) = {\rm Res}_Z\left(\sum_{k=1}^n \frac{1}{k!}  \omega^k\right) \in \Omega^{\bullet}_{\rm cl}(Z)
\ee
forms a $1$-cosymplectic global pure spinor on $Z$, hence defines a $1$-cosymplectic structure. In particular, $Z$ inherits a unimodular corank-$1$ symplectic foliation.
\end{prop}
We are asserting here that the residue map should be applied to the pure spinor $e^\omega$ associated to the $\mc{A}_Z$-Dirac structure defined by $\omega$ (see \cite{Klaasse18} for more discussion on this process).
\bp Because ${\rm Res}_Z$ is a cochain morphism (due to \autoref{thm:residues}), we have
\be
	d\rho = d \, {\rm Res}_Z\left(\sum_{k=1}^n \frac{1}{k!} \omega^k\right) = {\rm Res}_Z\left(\sum_{k=1}^n \frac{1}{k!} d (\omega^k)\right) = 0.
\ee
 Let $z$ be a local defining function for $\iota_Z\colon Z \hookrightarrow X$, so that locally $\omega = d\log z \wedge \wt{\alpha} + \wt{\beta}$. Then 
\be
	\omega^k = k (d\log z \wedge \wt{\alpha} + \wt{\beta}) \wedge \wt{\beta}^{k-1}, \qquad 1 \leq k \leq n.
\ee
Thus ${\rm Res}_Z(\omega^k) = k \alpha \wedge \beta^{k-1} \neq 0$ for $\alpha := {\rm Res}_Z(\omega) = \iota_Z^*(\wt{\alpha})$ and $\beta := \iota_Z^*(\wt{\beta})$ (see \autoref{exa:bkersiduemaps}). Both $\rho$ and $\alpha$ are invariantly defined closed forms on $Z$, and $\rho$ is a global pure spinor. We see that $\rho$ is locally given by
\be
	\rho = \alpha \wedge e^\beta.
\ee
Moreover, $\rho_{2n-1} = \frac{1}{(n-1)!} \alpha \wedge \beta^{n-1}$, which is nonzero due to nondegeneracy of $\omega$ (i.e.\ $\omega^n \neq 0$).
\ep
\begin{rem} \autoref{prop:logcosymplectic} is the intrinsic reformulation of \cite[Proposition 10]{GuilleminMirandaPires14}, that log-Poisson structures carry an equivalence class of cosymplectic structures $(\alpha,\beta)$, where $(\alpha,\beta) \sim (\alpha',\beta')$ if $\alpha = \alpha'$ and $\beta' = \beta + \alpha \wedge df$ for $f \in C^\infty(Z)$. There they note that the one-form $\alpha$ is determined invariantly, and here this is immediate because it arises from the logarithmic residue map ${\rm Res}_Z$. Any choice of (local) defining function for $Z$ picks out a representative for $\beta$, but the associated pure spinor $\rho$ is the invariant object, globally on $Z$. This is again because it is obtained by (repeated) application of the logarithmic residue.
\end{rem}
By repeating the proof of \autoref{prop:logcosymplectic} essentially verbatim, we obtain the following invariant statement about the induced geometry on the degeneracy loci of $b^k$-Poisson structures, which lead to symplectic Lie algebroids by \autoref{exa:bkpoissonlift} (compare this with \cite{Scott16,GuilleminMirandaWeitsman17}).
\begin{prop}\label{prop:bkpoissoncosymplectic} Let $(X,Z,\pi)$ be a $b^k$-Poisson manifold. Then $Z$ carries a $1$-cosymplectic global pure spinor, so that it inherits a $1$-cosymplectic structure.
\end{prop}
There is a similar result for elliptic Poisson manifolds, as follows. Due to \autoref{exa:ellpoissonlift}, an elliptic Poisson structure $\pi \in {\rm Poiss}(X)$ (with divisor ideal $I_{|D|}$) is of nondegenerate $\mc{A}_{|D|}$-type, with dual $\mc{A}_{|D|}$-symplectic structure $\omega$. We will focus on the case where ${\rm Res}_q(\omega) = 0$, which are those with zero elliptic residue, or \emph{zero elliptic Poisson structures}. We will moreover require that $D$ is cooriented (which is immediate if $\pi$ stems from a stable generalized complex structure \cite{CavalcantiGualtieri18}). As mentioned in \autoref{exa:elltangentresidues}, this trivializes the isotropy line $\mf{k} \to D$, so that the residue one-forms ${\rm Res}_r(\omega) \in \Omega^1(D)$ and ${\rm Res}_\theta(\omega) \in \Omega^1(D)$ are both closed.
\begin{prop}\label{prop:elliptic2cosymplectic} Let $(X^{2n},|D|,\pi,\omega)$ be a cooriented zero elliptic Poisson manifold. Then
\be
	\rho := {\rm Res}_q(e^\omega) = -{\rm Res}_r(\omega) \wedge {\rm Res}_{\theta}(\omega) + {\rm Res}_q \left(\sum_{k=3}^{n} \frac1{k!} \omega^k \right) \in \Omega^\bullet_{\rm cl}(D)
\ee
forms a $2$-cosymplectic global pure spinor on $D$, hence defines a $2$-cosymplectic structure. In particular, $D$ inherits a unimodular corank-$2$ symplectic foliation.
\end{prop}
There are extensions of this statement when $D$ is not cooriented (see \cite{KlaasseLanius18,KlaasseLi18}). Moreover, this result is known in the context of stable generalized complex geometry due to \cite{BaileyCavalcantiGualtieri17,CavalcantiGualtieri18}. There are three residue maps, namely ${\rm Res}_r$, ${\rm Res}_\theta$ and ${\rm Res}_q$, so that the pure spinor above differs slightly from \autoref{prop:logcosymplectic}. We have that ${\rm Res}_r(\omega) \wedge {\rm Res}_\theta(\omega) = -{\rm Res}_q(\omega^2)$. As for \autoref{prop:logcosymplectic} we apply ${\rm Res}_q$ to the pure spinor $e^\omega$ defining the $\mc{A}_{|D|}$-symplectic structure.
\bp Let $(r,\theta)$ be local normal polar coordinates to $\iota_D\colon D \hookrightarrow X$ in which $I_{|D|} = \langle r^2 \rangle$. Then locally $\omega = d\log r \wedge \wt{\alpha}_1 + d\theta \wedge \wt{\alpha}_2 + \wt{\beta}$. There is no term involving $d\log r \wedge d\theta$ because ${\rm Res}_q(\omega) = 0$. We have $\alpha_1 := {\rm Res}_r(\omega) = \iota_D^*(\wt{\alpha}_1)$ and $\alpha_2 := {\rm Res}_\theta(\omega) = \iota_D^*(\wt{\alpha}_2)$. These are both closed because the residue maps ${\rm Res}_r$ and ${\rm Res}_\theta$ are cochain morphisms. We have for $2 \leq k \leq n$ that
\be
	\omega^k = (k d\log r \wedge \wt{\alpha}_1 \wedge \beta + k d\theta \wedge \wt{\alpha}_2 \wedge \beta + k(k-1) d \log r \wedge \wt{\alpha}_1 \wedge d\theta \wedge \wt{\alpha}_2 + \wt{\beta}^2) \wedge \wt{\beta}^{k-2}.
\ee
Note here that $n \geq 2$ is forced by nondegeneracy and the condition that ${\rm Res}_q(\omega) = 0$. Thus
\be
	{\rm Res}_q(\omega^k) = -k (k-1) \alpha_1 \wedge \alpha_2 \wedge \beta^{k-2}, \qquad 2 \leq k \leq n,
\ee
where $\beta := \iota_D^*(\wt{\beta})$. In particular this shows that ${\rm Res}_q(\omega^2) = -{\rm Res}_r(\omega) \wedge {\rm Res}_\theta(\omega)$. Further, $\rho$ is closed because ${\rm Res}_q$ is a cochain morphism. Hence $\rho$, $\alpha_1$ and $\alpha_2$ are invariantly defined forms on $D$, and $\rho$ is a $2$-cosymplectic global pure spinor. We see that $\rho$ is locally given by
\be
	\rho = -\alpha_1 \wedge \alpha_2 \wedge e^\beta.
\ee
Moreover, we have that $\rho_{2n-2} = -\frac{1}{(n-2)!} \alpha_1 \wedge \alpha_2 \wedge \beta^{n-2} \neq 0$, due to nondegeneracy of $\omega$.
\ep
There is also a result for elliptic Poisson structures $\pi$ whose elliptic residue is nonzero, in other words, for \emph{nonzero elliptic Poisson structures} (for which $D$ is cooriented). Namely, in this case $D$ is a symplectic leaf of $\pi$. See \cite{KlaasseLanius18,KlaasseLi18} for more information.

Note that in the presence of $n$ rational one-forms $\alpha_i \in \Omega^1(M)$ whose cohomology classes are linearly independent, Tischler's theorem \cite{Tischler70} implies that $M$ naturally fibers over $T^n$. An $n$-cosymplectic tuple whose one-forms satisfy this condition is called \emph{proper}. Any $n$-cosymplectic tuple on a compact manifold can be perturbed slightly to be proper. Call a Poisson structure \emph{proper} if its induced geometric structure on its degeneracy locus is proper. We obtain the following result for cooriented zero elliptic Poisson manifolds. This is the direct Poisson geometric analogue results in \cite{BaileyCavalcantiGualtieri17} and a preliminary version of \cite{CavalcantiGualtieri18} in the setting of generalized complex geometry.
\begin{prop} Let $(X,|D|,\pi)$ be a cooriented zero elliptic Poisson manifold. Then $D$ carries a fibration over $T^2$, which can be made into a symplectic fibration after slight perturbation. Moreover, $\pi$ can be perturbed through zero elliptic Poisson structures to be proper.
\end{prop}
\begin{rem} A similar result holds for log-Poisson manifolds (or $b^k$-Poisson) (\cite{GuilleminMirandaPires14,MarcutOsornoTorres14,Cavalcanti17}): by \autoref{prop:logcosymplectic} the degeneracy locus $Z$ carries a $1$-cosymplectic structure, hence fibers over $S^1$. This can be perturbed to be proper, making $Z$ a symplectic mapping torus. In \cite{Cavalcanti17,MarcutOsornoTorres14} it is shown that $\pi$ can be perturbed through log-Poisson structures so that $Z$ is proper. The argument used there proves essentially verbatim that the same holds for $b^k$-Poisson structure.
\end{rem}
\subsubsection{Relation with the Poisson residue}
There is a relation between the Lie algebroid residue maps that were used above, and the modular Poisson residues of \autoref{sec:poissonmodulesmodfoliation} due to \cite{GualtieriLi14,GualtieriPym13}.

Let $\mc{A} \to X^{2n}$ be a Lie algebroid of smooth divisor-type and let $\pi \in {\rm Poiss}(X)$ be of nondegenerate $\mc{A}$-type with dual $\mc{A}$-symplectic structure $\omega$. This forces $I_\pi = I_\mc{A}$ by \autoref{prop:apoissondivtypelifts}, and it ensures that $Z_\mc{A} = Z_{I_\mc{A}}$ is both the degeneraci locus of $\mc{A}$, and the $(2n-2)$nd degeneracy locus of $\pi$. Assuming that the restriction of $\mc{A}$ to $Z_\mc{A}$ is transitive with unimodular isotropy, and that $Z_\mc{A}$ is of codimension-one, the residues ${\rm Res}_{Z_\mc{A}}(\omega^n) \in \Omega^{2n-1}(Z_\mc{A})$ and ${\rm Res}_{n-1,{\rm mod}}(\pi) \in \mf{X}^{2n-1}(Z_\mc{A})$ are dual to each other in some sense. In particular, this happens for $b^k$-Poisson structures and their associated $\mc{A}_Z^k$-symplectic structure. From \autoref{prop:bkpoissoncosymplectic} we see that ${\rm Res}_Z(\omega^n) = \alpha \wedge \beta^{n-1}$ if locally $\omega = dz/z^k \wedge \wt{\alpha} + \wt{\beta}$ for $z \in j_{k-1}$. On the other hand, we have that ${\rm Res}_{n-1,{\rm mod}}(\pi) = v \wedge \pi_Z^{n-1}$ for $v$ a local modular vector field. Note that $\beta$ is not globally defined, and nor is $v$. However, their pure spinors make sense, and define the same Dirac structure (the forms we used above as spinors are also called contravariant spinors, with multivector fields being covariant spinors, see \cite{AlekseevBursztynMeinrenken09,Klaasse18}). See also the discussions in \cite[Proposition 1.8]{GualtieriLi14} and \cite[Remark 7.3]{GualtieriPym13}. In general the relation between the Poisson modular residue and the residue maps of the associated symplectic Lie algebroid can be more involved. 
%
%END OF TEX FILE